\documentclass[a4paper,11pt, leqno]{amsart}
\usepackage[T1]{fontenc}
\usepackage{lmodern}
\usepackage{todonotes}
\usepackage{amssymb}
\usepackage{amsthm}
\usepackage{mathrsfs, amsfonts, amsmath}
\usepackage{xcolor}
\usepackage{graphicx}
\usepackage{tikz}
\usetikzlibrary{
  lindenmayersystems,
  arrows.meta,
  decorations.markings
}
\usetikzlibrary{arrows.meta}
\usepackage{tikz,graphicx}
\usetikzlibrary{lindenmayersystems}
\usetikzlibrary{patterns}
\usepackage{physics}
\usepackage[final]{pdfpages}
\usepackage{comment}
\usepackage{nicefrac}
\allowdisplaybreaks
\usepackage[alphabetic,initials,nobysame]{amsrefs}

\AtBeginDocument{%
  \def\MR#1{}
}

\newcommand{\diam}{\operatorname{diam}}
\newcommand{\dist}{\operatorname{dist}}

\newcommand{\s}{{\mathbb S}}

\newcommand{\modu}{\operatorname{mod}}

\newcommand{\hatc}{\hat{\mathbb{C}}}

\numberwithin{equation}{section} \theoremstyle{plain}
\newtheorem{theorem}{\textbf{THEOREM}}[section]
\newtheorem{lemma}[theorem]{\textsc{Lemma}}
\newtheorem{proposition}[theorem]{\textsc{Proposition}}
\newtheorem{corollary}[theorem]{\textsc{Corollary}}
\theoremstyle{definition}
\newtheorem{definition}[theorem]{\textsc{Definition}}
{\theoremstyle{remark} }
\def\charfn_#1{{\raise1.2pt\hbox{$\chi_{\kern-1pt\lower3pt\hbox{{$\scriptstyle#1$}}}$}}}
\def\leq{\leqslant }
\def\geq{\geqslant }

\def\XXint#1#2#3{{\setbox0=\hbox{$#1{#2#3}{\int}$}
\vcenter{\hbox{$#2#3$}}\kern-.5\wd0}}

\def\le {\leqslant}
\def\ge {\geqslant}

\usepackage{latexsym,amsthm,amsmath,amssymb,mathrsfs,mathtools}
\usepackage[T1]{fontenc}
\usepackage[utf8]{inputenc}
\usepackage[mathscr]{eucal}
\usepackage{enumerate}
\usepackage{enumitem}
\usepackage{color}
\usepackage{xcolor}

\usepackage{ esint }
\usepackage[colorlinks	=	true,						
linkcolor	=	red,
urlcolor	=	red,
citecolor	=	red]%
{hyperref}
\usepackage{cleveref}
\usepackage{graphicx}
\usepackage{tikz}
\usetikzlibrary{decorations.pathreplacing}

\def\diam{\operatorname{diam}}
\def\card{\operatorname{card}}
\def\dist{\operatorname{dist}}
\def\area{\operatorname{area}}

\DeclareMathOperator{\Mod}{mod}

\DeclareMathOperator{\D}{\mathbb{D}}
\DeclareMathOperator{\C}{\mathbb{C}}

\title[Uniformization of uniformly non-collinear domains]{Conformal uniformization of uniformly non-collinear domains}
\author{Chengxi Li and Kai Rajala} 
\address{Department of Mathematics and Statistics, University of Jyv\"askyl\"a, P.O. Box 35, FI-40014, University of Jyv\"askyl\"a, Finland. }
\email { cheli@jyu.fi, kai.i.rajala@jyu.fi}%

\date{\today}

\begin{document}

	\thanks{  
		\newline C.L. was supported by the Research Council of Finland, project numbers 364210 and 360505, and China Scholarship Council Fellowship, project number 202306200022. K.R. was supported by the Research Council of Finland, project number 360505. }

\begin{abstract} 
We prove a strong version of \emph{Koebe's conjecture} for domains $\Omega \subset \hatc$ whose complementary components satisfy a \emph{uniform non-collinearity condition} and a \emph{packing condition}: there exists a conformal map $f:\Omega \to D$, which preserves the classes of point components and non-trivial components, onto a \emph{circle domain} $D$. 

We demonstrate the sharpness of the conditions by presenting a countably connected slit domain $\Omega \subset \hatc$ for which the packing condition holds, but every conformal map $f:\Omega \to D$ onto a circle domain maps a point component to a disk. The non-collinearity condition is satisfied by Schramm's \emph{cofat domains} and by Esmayli-Rajala's \emph{cospread domains}. 
\end{abstract}
	
\maketitle

\section{Introduction}

\subsection{Background} \label{sec:background}

This paper is concerned with the conformal uniformization of multiply connected subdomains of the Riemann sphere $\hatc$. We seek conformal maps from domains satisfying geometric conditions onto \emph{circle domains}, i.e., domains whose boundary components are circles or points. 

In 1908, Koebe \cite{koebe1907uniformisierung} conjectured that every subdomain of $\hatc$ admits a conformal map onto a circle domain. The simply connected case is the Riemann mapping theorem. The conjecture was established for finitely connected domains by Koebe \cite{koebe1920abhandlungen}, and more than 70 years later for countably connected domains  by He and Schramm \cite{he1993fixed}. A complete solution of Koebe's conjecture was recently announced in \cite{OpenAI26}; see also \cites{KarNt26b,NtRa25,Nta26Survey,NtaYou26,Raj25} and the references below. 

Conformal maps $f:\Omega \to D$ onto circle domains may map point components to non-trivial components and vice versa, even when $\Omega$ is countably connected; see \cite{ntalampekos2023rigidity}*{Section 6}, 
\cite{LiRaj26}*{Theorems 1.3 and 1.4}. Much of the recent research on Koebe's conjecture is aimed at proving the existence of an $f$ as above under \emph{geometric conditions} in $\Omega$, such that \emph{the types (point or non-trivial) of complementary components are preserved by $f$}. 

The first significant result in this direction is due to Schramm \cite{schramm1995transboundary}, whose \emph{cofat uniformization theorem} gives the above conclusion under the assumption that the complementary components of $\Omega$ are \emph{uniformly Ahlfors $2$-regular}. He also introduced the \emph{transboundary modulus} of path families, which has become the main tool for various uniformization theories; see \cites{bonk2011uniformization,bonk2013quasisymmetric,HaLi23,HakReh25,LiRaj26,Mer25,Nta20} and the discussion below. 

Recently, Esmayli and the second author \cite{esmayli2024conformal} proved the above conclusion for \emph{cospread domains} whose complementary components contain uniformly ``spread'' \emph{tripods}; see Section \ref{sec:additional}. Another line of research was developed by Karafyllia and Ntalampekos \cites{karafyllia2024uniformization,KarNt26b}. Their assumptions, such as the \emph{John condition}, concern hyperbolic-type metrics on $\Omega$ instead of the geometry of individual complementary components.  

In this work, we introduce the \emph{uniformly non-collinear (UNC) domains} whose complementary components are uniformly far from being segments. We show in Theorem \ref{thm:mainthm} that the above conclusion holds for UNC-domains that satisfy a \emph{packing condition} that was also applied in \cite{esmayli2024conformal}. The UNC-condition is simple and geometrically appealing. It is satisfied by all cofat and cospread domains; see Section \ref{sec:additional} and Section \ref{sec:cofatspread}. Theorem \ref{thm:counterexample} illustrates the sharpness of the conditions in Theorem \ref{thm:mainthm}. 

\subsection{Main results} \label{sec:mainresults}

We denote the collection of complementary components of a domain $G \subset \hatc$ by $\mathcal{C}(G):=\mathcal{C}_N(G) \cup \mathcal{C}_P(G)$; here $p \in \mathcal{C}_N(G)$ are the \emph{non-trivial} components of positive diameter, and $p \in \mathcal{C}_P(G)$ the \emph{point components}. Every homeomorphism $f:G \to G'$ has a unique extension to a homeomorphism $\hat f: \hat G \to \hat G'$ between the \emph{quotient spaces} $\hat G:= G \cup \mathcal{C}(G)$ and $\hat G':= G' \cup \mathcal{C}(G')$, see Section \ref{sec:preliminaries}. 

We denote the solid triangle in $\mathbb{C}$ with vertices $a,b,c$ by $\triangle(a,b,c)$. 

\begin{definition}
	We call a domain $\Omega \subset \hatc$, $\infty \in \Omega$, \emph{uniformly non-collinear (UNC or $\varepsilon$-UNC for short)} if there exists $\varepsilon>0$ such that for every $p \in \mathcal{C}_N(\Omega)$ there are points $a_p, b_p, c_p\in p$ for which  
	\begin{equation}\label{eq:non-collineardef}
		\area(\triangle(a_p,b_p,c_p))\ge \varepsilon\diam(p)^2.  
	\end{equation}
 \end{definition}

Note that $\triangle(a_p,b_p,c_p)$ is not required to be a subset of $p$. Our first main result is a uniformization theorem for a subclass of UNC-domains. 

\begin{theorem}\label{thm:mainthm}
	Let $\Omega\subset \hatc$, $\infty \in \Omega$, be a UNC-domain that satisfies the following packing condition: there exists $N\ge 1$ such that for every $z_0 \in \mathbb{C}$ and $r>0$, we have
\begin{equation}\label{eq:packingcondition}
	\card \{p \in \mathcal{C}_N(\Omega): \,\diam(p) \geq r, \,  p \cap \mathbb{D}(z_0,r) \neq \emptyset \} \leq N.
\end{equation}
Then there exists a conformal map $f\colon\Omega \to D$ onto a circle domain $D \subset \hatc$. Moreover, 
	\begin{equation} \label{eq:claimi2}
		\hat{f}(\mathcal{C}_N(\Omega))=\mathcal{C}_N(D) \quad \text{and} \quad  \hat{f}(\mathcal{C}_P(\Omega))=\mathcal{C}_P(D). 
	\end{equation}
\end{theorem}

By \cite{esmayli2024conformal}*{Proposition 6.1}, the UNC-condition alone does not guarantee \eqref{eq:claimi2}. Our second main result shows that the packing  condition \eqref{eq:packingcondition} does not guarantee \eqref{eq:claimi2} either. 

\begin{theorem}\label{thm:counterexample}
There exists a domain $\Omega'\subset \hatc$, $\infty \in \Omega'$, such that $\mathcal{C}(\Omega')$ is countable and \eqref{eq:packingcondition} holds, but 
$$ 
\{0\}\in \mathcal{C}_P(\Omega') \quad \text{and} \quad \hat{f}(\{0\}) \in \mathcal{C}_N(D) 
$$ 
for every conformal map $f: \Omega' \to D$ onto a circle domain $D \subset \hatc$.
\end{theorem}

Theorem \ref{thm:counterexample} was stated in \cite{esmayli2024conformal}*{Remark 6.2} without proof. We give a proof of Theorem \ref{thm:counterexample} with complete details in Section \ref{sec:counterexample}. 

\subsection{Additional comments} \label{sec:additional}

Schramm \cite{schramm1995transboundary} proved Theorem \ref{thm:mainthm} for \emph{cofat domains} $\Omega$, that is, assuming that there exists $\tau>0$ such that every $p \in \mathcal{C}_N(\Omega)$ satisfies 
\[
\operatorname{area}(p\cap \mathbb{D}(z,r))\geq \tau\,r^2 \quad 
\text{for every } z \in p \text{ and } 0<r \leq \diam(p). 
\]
Recently, Esmayli and the second author \cite{esmayli2024conformal} proved a version of Theorem \ref{thm:mainthm} where the UNC-condition is replaced with a stronger condition involving \emph{quasitripods}, i.e., images of the standard tripod under \emph{quasisymmetric maps} (see Section \ref{sec:cofatspread}): there exists a homeomorphism $\eta:[0,\infty)\to [0,\infty)$ such that every $p \in \mathcal{C}_N(\Omega)$ contains an $\eta$-quasitripod  $T_p$ such that 
\begin{equation}\label{eq:quasitripod} 
\diam(T_p) \geq \diam(p)/\eta(1). 
\end{equation}

Theorem \ref{thm:mainthm} unifies the above results, answering an inquiry in  \cite{esmayli2024conformal}*{page 559}; the class of domains $\Omega$ that satisfy the assumptions of Theorem \ref{thm:mainthm} contains all cofat domains, as well as the domains for which \eqref{eq:quasitripod} and the packing condition \eqref{eq:packingcondition} hold; see Section \ref{sec:cofatspread}. 

Cofatness is preserved under Möbius transformations (see \cite{schramm1995transboundary}*{Theorem 2.1}) and cofat domains satisfy the packing condition \eqref{eq:packingcondition}. Similarly, assuming \eqref{eq:quasitripod} at all scales leads to the Möbius invariant class of \emph{cospread domains} (\cite{esmayli2024conformal}*{Definition 1.4}), which satisfy \eqref{eq:packingcondition}. One could hope that analogous localizations of the non-invariant UNC-condition lead to a Möbius invariant class of domains for which the packing condition \eqref{eq:packingcondition} holds. 

\begin{figure}[htbp]
\centering
\begingroup

\def\UNCArcColor{blue!65!black}

\def\UNCseedone{%
  0.00000000/0.00000000,0.10000000/0.00000000,0.36666667/0.00000000,
  0.50000000/-0.23094011,0.63333333/0.00000000,0.90000000/0.00000000,
  1.00000000/0.00000000
}

\def\UNCseedtwo{%
  0.00000000/0.00000000,0.10000000/0.00000000,0.18888889/0.00000000,
  0.23333333/-0.07698004,0.27777778/0.00000000,0.36666667/0.00000000,
  0.41111111/-0.07698004,0.36666667/-0.15396007,0.45555556/-0.15396007,
  0.50000000/-0.23094011,0.54444444/-0.15396007,0.63333333/-0.15396007,
  0.58888889/-0.07698004,0.63333333/0.00000000,0.72222222/0.00000000,
  0.76666667/-0.07698004,0.81111111/0.00000000,0.90000000/0.00000000,
  1.00000000/0.00000000
}

\newcommand{\UNCTrace}[3]{%
  \ifnum#2=1\relax
    \let\UNCseed\UNCseedone
  \else
    \let\UNCseed\UNCseedtwo
  \fi
  \foreach \UNCx/\UNCy [count=\UNCi,
    remember=\UNCx as \UNCprevx,
    remember=\UNCy as \UNCprevy] in {#3}{%
    \ifnum\UNCi>1\relax
      \pgfmathsetmacro{\UNClen}{
        veclen(\UNCx-\UNCprevx,\UNCy-\UNCprevy)
      }
      \pgfmathsetmacro{\UNCangle}{
        atan2(\UNCy-\UNCprevy,\UNCx-\UNCprevx)
      }
      \pgfmathtruncatemacro{\UNCn}{
        max(1,ceil(\UNClen/#1))
      }
      \pgfmathtruncatemacro{\UNClast}{\UNCn-1}
      \pgfmathsetmacro{\UNCtile}{\UNClen/\UNCn}
      \begin{scope}[
        shift={(\UNCprevx,\UNCprevy)},
        rotate=\UNCangle
      ]
        \foreach \UNCk in {0,...,\UNClast}{%
          \draw[unc approximant] ({\UNCk*\UNCtile},0)
            \foreach \UNCtx/\UNCty in \UNCseed {
              -- ({(\UNCk+\UNCtx)*\UNCtile},
                  {\UNCty*\UNCtile})
            };
        }
      \end{scope}
    \fi
  }
}

\resizebox{.90\linewidth}{!}{%
\begin{tikzpicture}[
  scale=5.2,
  line cap=butt,
  line join=miter,
  miter limit=10,
  every node/.style={font=\normalsize},
  unc limit/.style={
    draw=black,
    line width=.85pt
  },
  unc approximant/.style={
    draw=\UNCArcColor,
    line width=.45pt
  },
  unc boundary/.style={
    draw=black!65,
    line width=.40pt,
    dash pattern=on 4pt off 3pt
  },
  unc leader/.style={
    draw=black!65,
    line width=.30pt,
    line cap=round,
    densely dotted,
    shorten <=1pt
  },
  unc label/.style={
    anchor=west,
    inner sep=.7pt,
    font=\scriptsize
  }
]

\UNCTrace{0.090}{2}{%
  -0.79000000/-0.30000000,-0.53361541/-0.30000000,-0.52666667/-0.28796443,
  -0.51971793/-0.30000000,-0.18250430/-0.30000000,-0.01389748/-0.00796443,
  -0.02084622/0.00407114,-0.00694874/0.00407114,0.00000000/0.01610671,
  0.00694874/0.00407114,0.02084622/0.00407114,0.01389748/-0.00796443,
  0.18250430/-0.30000000,0.51971793/-0.30000000,0.52666667/-0.28796443,
  0.53361541/-0.30000000,0.79000000/-0.30000000
}

\UNCTrace{0.060}{2}{%
  -0.79000000/-0.26500000,-0.55382267/-0.26500000,-0.52666667/-0.21796443,
  -0.49951067/-0.26500000,-0.20271156/-0.26500000,-0.05431200/-0.00796443,
  -0.08146800/0.03907114,-0.02715600/0.03907114,0.00000000/0.08610671,
  0.02715600/0.03907114,0.08146800/0.03907114,0.05431200/-0.00796443,
  0.20271156/-0.26500000,0.49951067/-0.26500000,0.52666667/-0.21796443,
  0.55382267/-0.26500000,0.79000000/-0.26500000
}

\UNCTrace{0.050}{2}{%
  -0.79000000/-0.23500000,-0.57114317/-0.23500000,-0.52666667/-0.15796443,
  -0.48219016/-0.23500000,-0.22003206/-0.23500000,-0.08895302/-0.00796443,
  -0.13342952/0.06907114,-0.04447651/0.06907114,0.00000000/0.14610671,
  0.04447651/0.06907114,0.13342952/0.06907114,0.08895302/-0.00796443,
  0.22003206/-0.23500000,0.48219016/-0.23500000,0.52666667/-0.15796443,
  0.57114317/-0.23500000,0.79000000/-0.23500000
}

\UNCTrace{0.040}{1}{%
  -0.79000000/-0.21000000,-0.70261397/-0.21000000,-0.70222222/-0.20932148,
  -0.70183048/-0.21000000,-0.58557693/-0.21000000,-0.52745016/-0.10932148,
  -0.52784190/-0.10864295,-0.52705841/-0.10864295,-0.52666667/-0.10796443,
  -0.52627492/-0.10864295,-0.52549143/-0.10864295,-0.52588318/-0.10932148,
  -0.46775640/-0.21000000,-0.35150286/-0.21000000,-0.35111111/-0.20932148,
  -0.35071937/-0.21000000,-0.23446582/-0.21000000,-0.17633905/-0.10932148,
  -0.17673079/-0.10864295,-0.17594730/-0.10864295,-0.11782053/-0.00796443,
  -0.17594730/0.09271409,-0.17673079/0.09271409,-0.17633905/0.09339262,
  -0.17673079/0.09407114,-0.17594730/0.09407114,-0.17555556/0.09474967,
  -0.17516381/0.09407114,-0.05891026/0.09407114,-0.00078349/0.19474967,
  -0.00117524/0.19542819,-0.00039175/0.19542819,0.00000000/0.19610671,
  0.00039175/0.19542819,0.00117524/0.19542819,0.00078349/0.19474967,
  0.05891026/0.09407114,0.17516381/0.09407114,0.17555556/0.09474967,
  0.17594730/0.09407114,0.17673079/0.09407114,0.17633905/0.09339262,
  0.17673079/0.09271409,0.17594730/0.09271409,0.11782053/-0.00796443,
  0.17594730/-0.10864295,0.17673079/-0.10864295,0.17633905/-0.10932148,
  0.23446582/-0.21000000,0.35071937/-0.21000000,0.35111111/-0.20932148,
  0.35150286/-0.21000000,0.46775640/-0.21000000,0.52588318/-0.10932148,
  0.52549143/-0.10864295,0.52627492/-0.10864295,0.52666667/-0.10796443,
  0.52705841/-0.10864295,0.52784190/-0.10864295,0.52745016/-0.10932148,
  0.58557693/-0.21000000,0.70183048/-0.21000000,0.70222222/-0.20932148,
  0.70261397/-0.21000000,0.79000000/-0.21000000
}

\UNCTrace{0.028}{1}{%
  -0.79000000/-0.19000000,-0.71416097/-0.19000000,-0.70222222/-0.16932148,
  -0.69028347/-0.19000000,-0.59712394/-0.19000000,-0.55054417/-0.10932148,
  -0.56248292/-0.08864295,-0.53860542/-0.08864295,-0.52666667/-0.06796443,
  -0.51472792/-0.08864295,-0.49085041/-0.08864295,-0.50278916/-0.10932148,
  -0.45620940/-0.19000000,-0.36304986/-0.19000000,-0.35111111/-0.16932148,
  -0.33917236/-0.19000000,-0.24601283/-0.19000000,-0.19943306/-0.10932148,
  -0.21137181/-0.08864295,-0.18749431/-0.08864295,-0.14091454/-0.00796443,
  -0.18749431/0.07271409,-0.21137181/0.07271409,-0.19943306/0.09339262,
  -0.21137181/0.11407114,-0.18749431/0.11407114,-0.17555556/0.13474967,
  -0.16361680/0.11407114,-0.07045727/0.11407114,-0.02387750/0.19474967,
  -0.03581625/0.21542819,-0.01193875/0.21542819,0.00000000/0.23610671,
  0.01193875/0.21542819,0.03581625/0.21542819,0.02387750/0.19474967,
  0.07045727/0.11407114,0.16361680/0.11407114,0.17555556/0.13474967,
  0.18749431/0.11407114,0.21137181/0.11407114,0.19943306/0.09339262,
  0.21137181/0.07271409,0.18749431/0.07271409,0.14091454/-0.00796443,
  0.18749431/-0.08864295,0.21137181/-0.08864295,0.19943306/-0.10932148,
  0.24601283/-0.19000000,0.33917236/-0.19000000,0.35111111/-0.16932148,
  0.36304986/-0.19000000,0.45620940/-0.19000000,0.50278916/-0.10932148,
  0.49085041/-0.08864295,0.51472792/-0.08864295,0.52666667/-0.06796443,
  0.53860542/-0.08864295,0.56248292/-0.08864295,0.55054417/-0.10932148,
  0.59712394/-0.19000000,0.69028347/-0.19000000,0.70222222/-0.16932148,
  0.71416097/-0.19000000,0.79000000/-0.19000000
}

\UNCTrace{0.016}{1}{%
  -0.79000000/-0.17400000,-0.76241092/-0.17400000,-0.76074074/-0.17110716,
  -0.75907056/-0.17400000,-0.72339858/-0.17400000,-0.70556259/-0.14310716,
  -0.70723277/-0.14021432,-0.70389240/-0.14021432,-0.70222222/-0.13732148,
  -0.70055204/-0.14021432,-0.69721167/-0.14021432,-0.69888186/-0.14310716,
  -0.68104587/-0.17400000,-0.64537389/-0.17400000,-0.64370370/-0.17110716,
  -0.64203352/-0.17400000,-0.60636154/-0.17400000,-0.58852555/-0.14310716,
  -0.59019573/-0.14021432,-0.58685537/-0.14021432,-0.56901938/-0.10932148,
  -0.58685537/-0.07842864,-0.59019573/-0.07842864,-0.58852555/-0.07553579,
  -0.59019573/-0.07264295,-0.58685537/-0.07264295,-0.58518519/-0.06975011,
  -0.58351500/-0.07264295,-0.54784302/-0.07264295,-0.53000703/-0.04175011,
  -0.53167721/-0.03885727,-0.52833685/-0.03885727,-0.52666667/-0.03596443,
  -0.52499648/-0.03885727,-0.52165612/-0.03885727,-0.52332630/-0.04175011,
  -0.50549031/-0.07264295,-0.46981833/-0.07264295,-0.46814815/-0.06975011,
  -0.46647797/-0.07264295,-0.46313760/-0.07264295,-0.46480778/-0.07553579,
  -0.46313760/-0.07842864,-0.46647797/-0.07842864,-0.48431396/-0.10932148,
  -0.46647797/-0.14021432,-0.46313760/-0.14021432,-0.46480778/-0.14310716,
  -0.44697179/-0.17400000,-0.41129981/-0.17400000,-0.40962963/-0.17110716,
  -0.40795945/-0.17400000,-0.37228747/-0.17400000,-0.35445148/-0.14310716,
  -0.35612166/-0.14021432,-0.35278129/-0.14021432,-0.35111111/-0.13732148,
  -0.34944093/-0.14021432,-0.34610056/-0.14021432,-0.34777075/-0.14310716,
  -0.32993476/-0.17400000,-0.29426278/-0.17400000,-0.29259259/-0.17110716,
  -0.29092241/-0.17400000,-0.25525043/-0.17400000,-0.23741444/-0.14310716,
  -0.23908462/-0.14021432,-0.23574426/-0.14021432,-0.21790827/-0.10932148,
  -0.23574426/-0.07842864,-0.23908462/-0.07842864,-0.23741444/-0.07553579,
  -0.23908462/-0.07264295,-0.23574426/-0.07264295,-0.23407407/-0.06975011,
  -0.23240389/-0.07264295,-0.19673191/-0.07264295,-0.17889592/-0.04175011,
  -0.18056610/-0.03885727,-0.17722574/-0.03885727,-0.15938975/-0.00796443,
  -0.17722574/0.02292841,-0.18056610/0.02292841,-0.17889592/0.02582125,
  -0.19673191/0.05671409,-0.23240389/0.05671409,-0.23407407/0.05382125,
  -0.23574426/0.05671409,-0.23908462/0.05671409,-0.23741444/0.05960694,
  -0.23908462/0.06249978,-0.23574426/0.06249978,-0.21790827/0.09339262,
  -0.23574426/0.12428546,-0.23908462/0.12428546,-0.23741444/0.12717830,
  -0.23908462/0.13007114,-0.23574426/0.13007114,-0.23407407/0.13296398,
  -0.23240389/0.13007114,-0.19673191/0.13007114,-0.17889592/0.16096398,
  -0.18056610/0.16385682,-0.17722574/0.16385682,-0.17555556/0.16674967,
  -0.17388537/0.16385682,-0.17054501/0.16385682,-0.17221519/0.16096398,
  -0.15437920/0.13007114,-0.11870722/0.13007114,-0.11703704/0.13296398,
  -0.11536685/0.13007114,-0.07969487/0.13007114,-0.06185888/0.16096398,
  -0.06352907/0.16385682,-0.06018870/0.16385682,-0.04235271/0.19474967,
  -0.06018870/0.22564251,-0.06352907/0.22564251,-0.06185888/0.22853535,
  -0.06352907/0.23142819,-0.06018870/0.23142819,-0.05851852/0.23432103,
  -0.05684834/0.23142819,-0.02117636/0.23142819,-0.00334037/0.26232103,
  -0.00501055/0.26521387,-0.00167018/0.26521387,0.00000000/0.26810671,
  0.00167018/0.26521387,0.00501055/0.26521387,0.00334037/0.26232103,
  0.02117636/0.23142819,0.05684834/0.23142819,0.05851852/0.23432103,
  0.06018870/0.23142819,0.06352907/0.23142819,0.06185888/0.22853535,
  0.06352907/0.22564251,0.06018870/0.22564251,0.04235271/0.19474967,
  0.06018870/0.16385682,0.06352907/0.16385682,0.06185888/0.16096398,
  0.07969487/0.13007114,0.11536685/0.13007114,0.11703704/0.13296398,
  0.11870722/0.13007114,0.15437920/0.13007114,0.17221519/0.16096398,
  0.17054501/0.16385682,0.17388537/0.16385682,0.17555556/0.16674967,
  0.17722574/0.16385682,0.18056610/0.16385682,0.17889592/0.16096398,
  0.19673191/0.13007114,0.23240389/0.13007114,0.23407407/0.13296398,
  0.23574426/0.13007114,0.23908462/0.13007114,0.23741444/0.12717830,
  0.23908462/0.12428546,0.23574426/0.12428546,0.21790827/0.09339262,
  0.23574426/0.06249978,0.23908462/0.06249978,0.23741444/0.05960694,
  0.23908462/0.05671409,0.23574426/0.05671409,0.23407407/0.05382125,
  0.23240389/0.05671409,0.19673191/0.05671409,0.17889592/0.02582125,
  0.18056610/0.02292841,0.17722574/0.02292841,0.15938975/-0.00796443,
  0.17722574/-0.03885727,0.18056610/-0.03885727,0.17889592/-0.04175011,
  0.19673191/-0.07264295,0.23240389/-0.07264295,0.23407407/-0.06975011,
  0.23574426/-0.07264295,0.23908462/-0.07264295,0.23741444/-0.07553579,
  0.23908462/-0.07842864,0.23574426/-0.07842864,0.21790827/-0.10932148,
  0.23574426/-0.14021432,0.23908462/-0.14021432,0.23741444/-0.14310716,
  0.25525043/-0.17400000,0.29092241/-0.17400000,0.29259259/-0.17110716,
  0.29426278/-0.17400000,0.32993476/-0.17400000,0.34777075/-0.14310716,
  0.34610056/-0.14021432,0.34944093/-0.14021432,0.35111111/-0.13732148,
  0.35278129/-0.14021432,0.35612166/-0.14021432,0.35445148/-0.14310716,
  0.37228747/-0.17400000,0.40795945/-0.17400000,0.40962963/-0.17110716,
  0.41129981/-0.17400000,0.44697179/-0.17400000,0.46480778/-0.14310716,
  0.46313760/-0.14021432,0.46647797/-0.14021432,0.48431396/-0.10932148,
  0.46647797/-0.07842864,0.46313760/-0.07842864,0.46480778/-0.07553579,
  0.46313760/-0.07264295,0.46647797/-0.07264295,0.46814815/-0.06975011,
  0.46981833/-0.07264295,0.50549031/-0.07264295,0.52332630/-0.04175011,
  0.52165612/-0.03885727,0.52499648/-0.03885727,0.52666667/-0.03596443,
  0.52833685/-0.03885727,0.53167721/-0.03885727,0.53000703/-0.04175011,
  0.54784302/-0.07264295,0.58351500/-0.07264295,0.58518519/-0.06975011,
  0.58685537/-0.07264295,0.59019573/-0.07264295,0.58852555/-0.07553579,
  0.59019573/-0.07842864,0.58685537/-0.07842864,0.56901938/-0.10932148,
  0.58685537/-0.14021432,0.59019573/-0.14021432,0.58852555/-0.14310716,
  0.60636154/-0.17400000,0.64203352/-0.17400000,0.64370370/-0.17110716,
  0.64537389/-0.17400000,0.68104587/-0.17400000,0.69888186/-0.14310716,
  0.69721167/-0.14021432,0.70055204/-0.14021432,0.70222222/-0.13732148,
  0.70389240/-0.14021432,0.70723277/-0.14021432,0.70556259/-0.14310716,
  0.72339858/-0.17400000,0.75907056/-0.17400000,0.76074074/-0.17110716,
  0.76241092/-0.17400000,0.79000000/-0.17400000
}

\pgfmathsetlengthmacro{\UNCstep}{1.58cm/3^4}

\begin{scope}[shift={(-.79,-.16)}]
  \draw[unc limit] (0,0)
    lindenmayer system[
      l-system={
        rule set={F -> F+F--F+F},
        axiom=F,
        order=4,
        angle=60,
        step=\UNCstep
      }
    ];
\end{scope}

\draw[unc leader] (.79,-.160) -- (.815,-.105);
\node[unc label] at (.830,-.105) {$S_0$};

\draw[unc leader] (.79,-.235) -- (.815,-.245);
\node[unc label] at (.830,-.245) {$S_3$};

\draw[unc leader] (.79,-.265) -- (.815,-.325);
\node[unc label] at (.830,-.325) {$S_2$};

\draw[unc leader] (.79,-.300) -- (.815,-.405);
\node[unc label] at (.830,-.405) {$S_1$};

\node[inner sep=0pt,font=\scriptsize]
  at (.850,-.150) {$\vdots$};

\end{tikzpicture}%
}
\endgroup
\caption{$S_j$, $j=1,2,\cdots$ are pairwise disjoint von Koch snowflake-type arcs that converge to $S_0$ in the Hausdorff sense. The local $\varepsilon$-UNC-condition \eqref{eq:thana} holds for all $S_j$, but  $\Omega:=\hatc\setminus \cup_{j=0}^\infty S_j$ fails the packing condition \eqref{eq:packingcondition}. }
\label{fig:snowflake-arcs}
\end{figure}
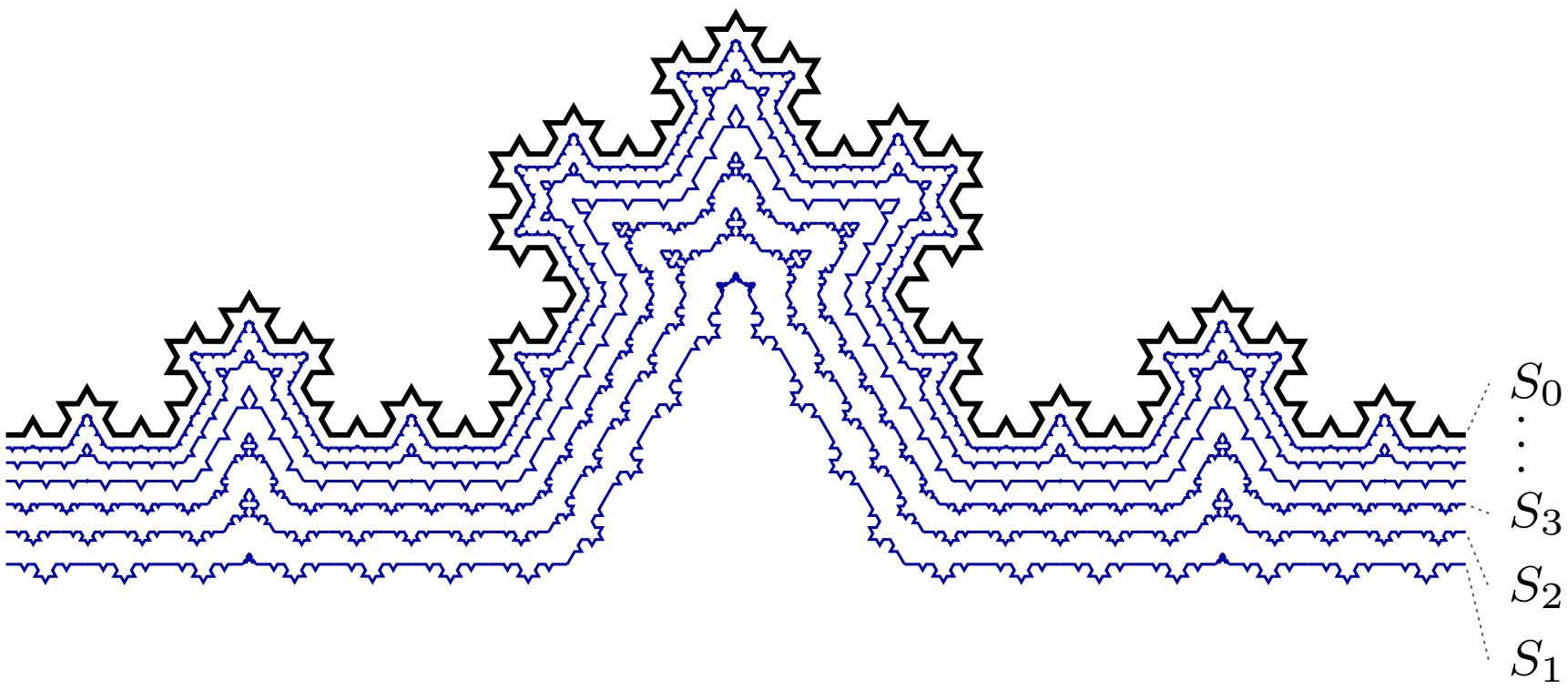

Unfortunately, this is not the case: there exist $\varepsilon>0$ and a sequence of pairwise disjoint von Koch snowflake-type arcs $S_j \subset \mathbb{D}(0,1)$, $j=0,1,2\ldots$, whose union is closed such that the local $\varepsilon$-UNC-condition
\begin{equation}\label{eq:thana}
S_j \cap \mathbb D(z_0,r) \text{ satisfies \eqref{eq:non-collineardef} for all } z_0 \in S_j, \, r<1, 
\end{equation}
holds and $\diam(S_j)\geq 1$ for all $j$; see Figure \ref{fig:snowflake-arcs}. In particular, the packing condition \eqref{eq:packingcondition} fails for domain $\Omega:=\hatc\setminus \cup_{j=0}^\infty S_j$. We are not aware of a localized version of the UNC-condition that implies \eqref{eq:packingcondition}. 

The proof of Theorem \ref{thm:mainthm} follows a standard approach: approximate $\Omega$ by a decreasing sequence of finitely connected domains $\Omega_j$, apply Koebe's theorem to find conformal maps $f_j:\Omega_j \to D_j$ onto circle domains, and apply transboundary modulus estimates to show that a conformal limit map preserves the types of complementary components and maps onto a circle domain. We discuss transboundary modulus in Section \ref{sec:preliminaries} and the outline of the proof of Theorem \ref{thm:mainthm} in Section \ref{sec:outline}. 

The challenging part is establishing the necessary transboundary modulus estimates. Our new contribution is Theorem \ref{thm:keyestimate}, an intricate modulus upper bound on finitely connected UNC-domains, which can be combined with known results to prove Theorem \ref{thm:mainthm}. 

The proof of Theorem \ref{thm:keyestimate} involves geometric and topological properties of triangles, and a construction of admissible functions by grouping together complementary components. Substantial difficulties arise from the fact that the standard admissible function, whose weight in $p \in \mathcal C_N(\Omega)$ equals $\diam(p)$, may have infinite energy. Our solution is to let the weight equal zero at many ``buried'' components, and use the UNC-condition to show that paths traveling through them must gain weight elsewhere. We prove Theorem \ref{thm:keyestimate} in Section \ref{sec:key}, which is divided into several subsections. 

To prove Theorem \ref{thm:counterexample}, in Section \ref{sec:counterexample} we construct the complement of a domain $\Omega$ as the union of increasingly dense horizontal slit configurations approaching $\infty$. We show that slits can be chosen such that the packing condition \eqref{eq:packingcondition} is satisfied, and prove transboundary modulus estimates showing that $\hat f(\{\infty\})$ is a disk for every conformal map $f$ onto a circle domain. The domain $\Omega'$ in Theorem \ref{thm:counterexample} is the inversion of $\Omega$. 

In Section \ref{sec:cofatspread}, we prove that the cofat domains and the domains with the quasitripod condition \eqref{eq:quasitripod} are UNC-domains. Thus, Theorem \ref{thm:mainthm} contains the cofat and cospread uniformization theorems as special cases.


\section{Transboundary modulus} \label{sec:preliminaries} 

	The proofs of our main results, Theorems \ref{thm:mainthm} and \ref{thm:counterexample}, are based on the \emph{transboundary modulus}, introduced by Schramm \cite{schramm1995transboundary}. The idea is to generalize the classical conformal modulus, which concerns families of paths in a domain $G \subset \hatc$, to the setting where the paths are also allowed to pass through the collection of complementary components $\mathcal{C}(G)$. More precisely, the paths travel in the quotient space $\hat G$ associated to $G$.   

We now define the transboundary modulus and recall its main properties. Given a domain $G \subset \hatc$, let $\hat{G}:=\hatc / \sim$, where $z \sim w$ if either $z=w \in G$ or $z,w \in p$ for some $p \in \mathcal{C}(G)$. Thus, $\hat G=G \cup \mathcal{C}(G)$. We equip $\hat G$ with the quotient topology and denote the induced quotient map by $\pi_G:\hatc\to \hat{G}$. 

By Moore's theorem \cite{MR2341468}*{page 3}, $\hat{G}$ is homeomorphic to $\hatc$. Every homeomorphism $f:G\to D\subset \hatc$ can be extended uniquely to a homeomorphism $\hat{f}:\hat{G}\to \hat{D}$ between their quotient spaces. 
In what follows, we consider $z \in G$ as an element of both $\hatc$ and $\hat G$, and $p \in \mathcal{C}(G)$ as both a subset of $\hatc$ and an element of $\hat G$, when there is no danger of confusion. 

\begin{definition}\label{lem:transmoddef}
		The \emph{transboundary modulus} $\Mod(\Gamma)$ of a family $\Gamma$ of paths in $\hat G$ is 
		$$
		\Mod(\Gamma)=\inf_{\rho \in Adm(\Gamma)} \int_{G\cap \mathbb{C}} \rho^2 \, dA + \sum_{p \in \mathcal{C}(G)} \rho(p)^2, 
		$$
		where $Adm(\Gamma)$ is the collection of \emph{admissible functions for $\Gamma$}, i.e., Borel measurable functions $\rho\colon \hat{G}\to [0,\infty]$ for which 
		$$
		1 \leq \int_{\gamma} \rho \, ds +\sum_{p \in \mathcal{C}(\gamma)} \rho(p)=:\hat \ell(\rho,\gamma)  \quad \text{for all } \gamma \in \Gamma.  
		$$ 
		Here $\mathcal{C}(\gamma):=|\gamma| \cap \mathcal{C}(G)$, and $\int_{\gamma} \rho \, ds$ is the path integral defined as follows: 
        the relatively open set $I \cap \gamma^{-1}(G)$ is the union of pairwise disjoint intervals $I_j$. We denote $\gamma_j:=\gamma_{|I_j}$, and define 
 		$$
		\int_{\gamma} \rho \, ds := \sum_j \int_{\gamma_j} \rho \, ds.  
		$$
        If $\gamma_j$ is not locally rectifiable for some $j$, then the integral is infinite. 
	\end{definition}
	
	We will repeatedly use the conformal invariance of the transboundary modulus, a basic fact also due to Schramm \cite{schramm1995transboundary}*{Lemma 1.1}, whose proof relies on a change-of-variables argument analogous to the one used to establish the conformal invariance of the classical conformal modulus.
	
	\begin{lemma}[\cite{schramm1995transboundary}, Lemma 1.1]\label{lem:confinvariance}
		Suppose that $f\colon G \to D$ is conformal. Then for every path family $\Gamma$ in $\hat{G}$, we have $\Mod(\Gamma)=\Mod(\hat{f}(\Gamma))$. Here $\hat{f}(\Gamma):=\{\hat{f} \circ \gamma: \, \gamma \in \Gamma\}$. 
	\end{lemma}

	We will also apply the following basic properties of transboundary modulus; see \cite{HaLi23}*{Proposition 3.1} for the proofs.
	\begin{proposition}\label{prop:propoftrans}
		Let $\Gamma_1,\Gamma_2,\ldots$ be path families in $\hat{G}$. The following properties hold: 
		\begin{enumerate}
			\item If $\Gamma_1 \subset \Gamma_2$, then $\Mod(\Gamma_1) \leq \Mod(\Gamma_2)$. 
			\item If $\Gamma=\cup_j \Gamma_j$, then $\Mod(\Gamma) \leq \sum_j \Mod(\Gamma_j)$. 
			\item If every $\gamma_2\in \Gamma_2$ has a subpath $\gamma_1\in \Gamma_1$, then 
			$\Mod(\Gamma_1) \geq \Mod(\Gamma_2)$.
			\item  If $p \in G$ or if $p$ is a point component isolated in $\mathcal C(G)$, then the modulus of all non-constant paths $\gamma$ in $\hat{G}$ satisfying $p \in |\gamma|$ is zero. 
		\end{enumerate}
	\end{proposition}


\section{Proof of Theorem \ref{thm:mainthm}} \label{sec:outline}
The proof of Theorem \ref{thm:mainthm} is a combination of an approximation argument and an intricate modulus estimate, Theorem \ref{thm:keyestimate} below. In this section, we recall the approximation argument, which is fairly standard; it is also used in \cite{esmayli2024conformal} to prove a version of Theorem \ref{thm:mainthm} where the UNC-condition is replaced by the stronger 
quasitripod condition \eqref{eq:quasitripod}. We refer to  \cite{esmayli2024conformal}*{Section 3} for the details that are omitted here. 

Let $\Omega \subset \hatc$ be a domain that satisfies the assumptions of Theorem \ref{thm:mainthm}. We may assume that $\mathcal{C}_N(\Omega)$ is infinite, since otherwise the theorem follows by applying Koebe's theorem \cite{koebe1920abhandlungen} to the domain $\hatc \setminus \cup_{p \in \mathcal C_N(\Omega)}$ p.  

Condition \eqref{eq:packingcondition} implies that $\mathcal{C}_N(\Omega)$ is countable. We enumerate the elements $\mathcal{C}_N(\Omega)=\{p_0,p_1,\ldots\}$. By Koebe's theorem, there are conformal maps $g_k$, $k=0,1,2,\ldots$, from $\tilde{\Omega}_k:=\hatc\setminus \cup_{l=0}^k p_l$ onto circle domains $\tilde{D}_k$. We normalize the maps $g_k$ so that they  all fix the same three points, one of which is $\infty$ (recall that $\infty \in \Omega$ by assumption). 

By compactness, there is a subsequence $f_j:=g_{k_j}:\Omega_j \to D_j$, $\Omega_j=\tilde \Omega_{k_j}$,  such that 
\begin{itemize}
    \item[(i)] $f_j$ converges to a conformal map $f:\Omega \to D$ locally uniformly. 
    \item[(ii)] The sequence $(\hat{f}_j(p_k))_j$ of disks converges 
    to a (possibly degenerate) closed disk $q_k$ in the Hausdorff sense for every $k=0,1,2,\ldots$. 
\end{itemize} 
By Carath\'eodory's kernel convergence theorem (see \cite{Gol:69}*{Theorem V.5.1, p.~228}) we have $q_k \subset \hat{f}(p_k)$ for every $k=0,1,\ldots$. We claim that the following properties hold:   
\begin{eqnarray}\label{eq:cruuprop1}
& & q_k = \hat{f}(p_k) \quad \text{for every } k=0,1,\ldots, \\ 
& & \label{eq:cruuprop2} \text{if } p \in \mathcal{C}_P(\Omega), 
\text{ then } \hat f(p) \text{ is a point}, \\ 
& & \label{eq:cruuprop3} \text{if } p_k \in \mathcal{C}_N(\Omega), 
\text{ then } q_k \text{ is a disk.} 
\end{eqnarray}

Note that \eqref{eq:cruuprop1} and \eqref{eq:cruuprop2} imply that $D$ is a circle domain. Moreover, combining \eqref{eq:cruuprop1}-\eqref{eq:cruuprop3} shows that 
$$
\hat f(\mathcal{C}_P(\Omega))=\mathcal{C}_P(D) \quad \text{and} \quad 
\hat f(\mathcal{C}_N(\Omega))=\mathcal{C}_N(D). 
$$

We conclude that Theorem \ref{thm:mainthm} follows from \eqref{eq:cruuprop1}-\eqref{eq:cruuprop3}, whose proofs are based on basic transboundary modulus estimates on circle domains and the following key estimate applied to the domains $G=\Omega_j$.  

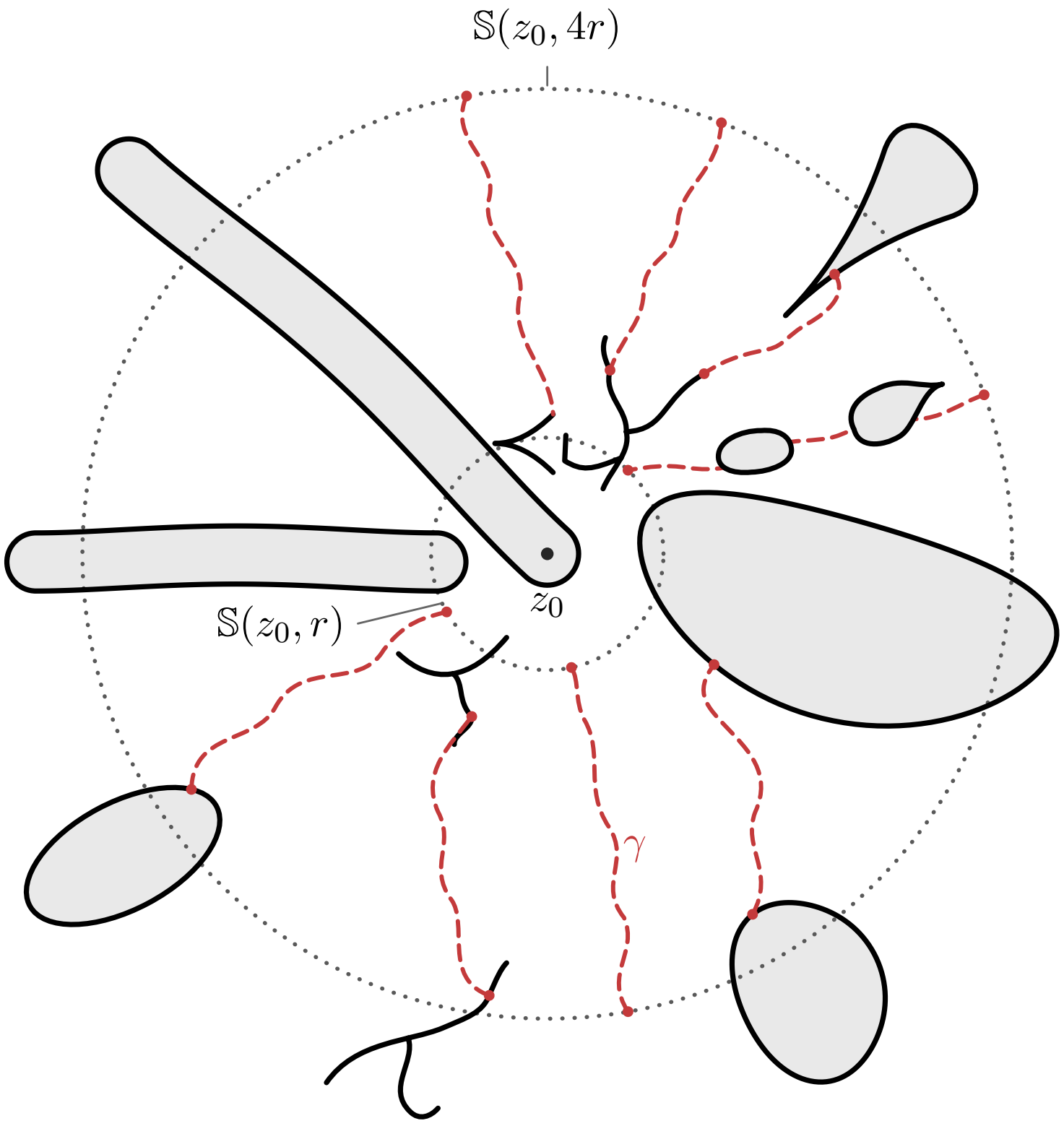
\begin{figure}[htbp]
\centering
\resizebox{0.82\linewidth}{!}{%
\begin{tikzpicture}[
  x=1.1cm,
  y=1.1cm,
  line cap=round,
  line join=round,
  component/.style={
    draw=black,
    solid,
    line width=1.4pt
  },
  boundary/.style={
    draw=black!65,
    line width=1.1pt,
    dash pattern=on 0pt off 3.2pt
  },
  pathcurve/.style={
    draw=pathRed,
    line width=1.15pt,
    dash pattern=on 4pt off 2.5pt,
    postaction={decorate},
    decoration={
      markings,
      mark=at position 0 with {
        \fill[pathRed] (0,0) circle[radius=1.45pt];
      },
      mark=at position 1 with {
        \fill[pathRed] (0,0) circle[radius=1.45pt];
      }
    }
  }
]

\definecolor{pathRed}{HTML}{C43A3B}
\definecolor{componentGray}{HTML}{E9E9E9}

\path[use as bounding box]
  (-5.05,-5.02) rectangle (4.70,4.65);

\def\annShapes#1{%
  \foreach \xx/\yy/\ang/\len/\wid/\bend in {
    -3.6/3.3/-42.510447/4.883646/.27/.14,
    -4.4/-.07/0/3.45/.25/.06
  }{
    \begin{scope}[shift={(\xx,\yy)},rotate=\ang]
      \path[#1]
        plot[domain=0:1,samples=65,variable=\t]
          ({\len*\t},{\wid+16*\bend*\t^2*(1-\t)^2})
        arc[start angle=90,end angle=-90,radius=\wid]
        -- plot[domain=1:0,samples=65,variable=\t]
          ({\len*\t},{-\wid+16*\bend*\t^2*(1-\t)^2})
        arc[start angle=-90,end angle=-270,radius=\wid]
        -- cycle;
    \end{scope}
  }

  \begin{scope}[shift={(2.55,-.42)},rotate=-16]
    \path[#1]
      plot[domain=0:360,samples=121,variable=\t]
        ({1.85*cos(\t)},{.86*sin(\t)+.15*cos(2*\t)})
      -- cycle;
  \end{scope}

  \path[#1,rotate around={28:(-3.65,-2.6)}]
    (-3.65,-2.6) ellipse[x radius=.91,y radius=.46];

  \begin{scope}[shift={(2.25,-3.75)},rotate=25]
    \path[#1]
      plot[domain=0:360,samples=121,variable=\t]
        ({.62*(1+.035*cos(3*\t))*cos(\t)},
         {.80*(1+.035*cos(3*\t))*sin(\t)})
      -- cycle;
  \end{scope}

  \begin{scope}[shift={(2.05,2.05)},rotate=45]
    \path[#1]
      (0,0)
      .. controls (.533333,0) and (1.066667,.133333)
      .. (1.6,.4)
      .. controls (1.9,.55) and (2.05,.32)
      .. (2.05,0)
      .. controls (2.05,-.32) and (1.9,-.55)
      .. (1.6,-.4)
      .. controls (1.066667,-.133333) and (.533333,0)
      .. (0,0)
      -- cycle;
  \end{scope}
}

\def\annInteriorCusp{%
  (3.40,1.46)
  .. controls (3.22,1.42) and (3.06,1.53) .. (2.84,1.42)
  .. controls (2.64,1.34) and (2.57,1.20) .. (2.64,1.07)
  .. controls (2.68,.91) and (2.88,.92) .. (3.03,1.01)
  .. controls (3.17,1.09) and (3.22,1.42) .. (3.40,1.46)
  -- cycle%
}
\def\annInteriorRounded{%
  (1.4724,0.8398)
  .. controls (1.4718,0.9713) and (1.6599,1.0619)
  .. (1.8424,1.0604)
  .. controls (1.9823,1.0724) and (2.0920,1.0190)
  .. (2.1009,0.8989)
  .. controls (2.1098,0.7788) and (1.9185,0.7080)
  .. (1.7573,0.7027)
  .. controls (1.5879,0.6860) and (1.4683,0.7379)
  .. (1.4724,0.8398)
  -- cycle%
}

\annShapes{fill=componentGray}
\path[fill=componentGray] \annInteriorCusp \annInteriorRounded;

\foreach \r in {1,4}
  \draw[boundary] (0,0) circle[radius=\r];

\annShapes{component}

\foreach \annSign in {-1,1}
  \draw[component]
    (-.45,.95)
    .. controls (-.283333,.95)
      and (-.116667,{.95+\annSign/12})
    .. (.05,{.95+\annSign/4});

\draw[component]
  (0.48,0.56)
  .. controls (0.55,0.73) and (0.72,0.87)
  .. (0.68,1.05)
  .. controls (0.69,1.25) and (0.48,1.36)
  .. (0.54,1.58)
  .. controls (0.48,1.7) and (0.47,1.78)
  .. (0.5,1.86)

  (0.68,1.05)
  .. controls (1.03,1.07) and (0.98,1.34)
  .. (1.35,1.55)

  (0.634674,0.824921)
  .. controls (0.38,0.72) and (0.3,0.71)
  .. (0.15,0.8)
  -- (0.16,1.02)

  (-1.28,-0.86)
  .. controls (-0.99,-1.13) and (-0.67,-1.1)
  .. (-0.35,-0.72)

  (-0.807482,-1.03081)
  .. controls (-0.71,-1.13) and (-0.78,-1.26)
  .. (-0.65,-1.4)
  .. controls (-0.64,-1.5) and (-0.81,-1.55)
  .. (-0.8,-1.64)

  (-1.9,-4.55)
  .. controls (-1.65,-4.19) and (-1.19,-4.2)
  .. (-0.9,-4.08)
  .. controls (-0.69,-3.99) and (-0.56,-3.95)
  .. (-0.5,-3.8)
  .. controls (-0.44,-3.66) and (-0.4,-3.57)
  .. (-0.35,-3.52)

  (-1.19664,-4.1664)
  .. controls (-1.13,-4.37) and (-1.39,-4.6)
  .. (-1.18,-4.8)
  .. controls (-1.1,-4.88) and (-1,-4.84)
  .. (-0.94,-4.77);

\draw[pathcurve]
  (1.498426,3.708735)
  .. controls (1.47,3.27) and (1.13,3.16)
  .. (1.12,2.82)
  .. controls (1.1,2.54) and (0.84,2.54)
  .. (0.83,2.28)
  .. controls (0.83,1.98) and (0.65,1.89)
  .. (0.54,1.58);

\draw[pathcurve]
  (-3.062142,-2.026939)
  .. controls (-3.05,-1.52) and (-2.51,-1.64)
  .. (-2.34,-1.32)
  .. controls (-2.11,-0.9) and (-1.62,-1.18)
  .. (-1.41,-0.83)
  .. controls (-1.22,-0.58) and (-1.06,-0.55)
  .. (-0.866025,-0.5);

\draw[pathcurve]
  (-0.5,-3.8)
  .. controls (-0.91,-3.65) and (-0.67,-3.18)
  .. (-0.83,-2.97)
  .. controls (-1.03,-2.69) and (-0.79,-2.5)
  .. (-0.94,-2.22)
  .. controls (-1.09,-1.88) and (-0.84,-1.65)
  .. (-0.65,-1.4);

\draw[pathcurve]
  (0.694593,-3.939231)
  .. controls (0.48,-3.61) and (0.79,-3.43)
  .. (0.66,-3.13)
  .. controls (0.43,-2.84) and (0.72,-2.64)
  .. (0.53,-2.35)
  .. controls (0.31,-2.02) and (0.49,-1.72)
  .. (0.32,-1.43)
  .. controls (0.23,-1.26) and (0.28,-1.1)
  .. (0.207912,-0.978148);

\draw[pathcurve]
  (1.766311,-3.101822)
  .. controls (1.97,-2.73) and (1.67,-2.53)
  .. (1.79,-2.22)
  .. controls (1.94,-1.89) and (1.72,-1.68)
  .. (1.55,-1.5)
  .. controls (1.38,-1.32) and (1.25,-1.2)
  .. (1.43487,-0.95306);

\draw[pathcurve]
  (2.472331,2.405487)
  .. controls (2.59,2.24) and (2.42,2.08)
  .. (2.24,1.98)
  .. controls (2.06,1.88) and (2.04,1.79)
  .. (1.9,1.76)
  .. controls (1.76,1.73) and (1.54,1.69)
  .. (1.35,1.55);

\draw[pathcurve]
  (100:4)
  .. controls (-.86,3.55) and (-.42,3.41)
  .. (-.50,3.06)
  .. controls (-.58,2.71) and (-.16,2.55)
  .. (-.24,2.20)
  .. controls (-.32,1.85) and (.03,1.58)
  .. (.05,1.20);

\begin{scope}[even odd rule]
  \clip
    (-5.05,-5.02) rectangle (4.70,4.65)
    \annInteriorCusp \annInteriorRounded;
  \draw[pathcurve]
    (20:4)
    .. controls (3.55,1.30) and (3.39,1.15)
    .. (3.14,1.16)
    .. controls (2.88,1.17) and (2.71,1.10)
    .. (2.51,1.02)
    .. controls (2.29,.94) and (2.09,1.02)
    .. (1.85,.88)
    .. controls (1.62,.73) and (1.43,.71)
    .. (1.21,.74)
    .. controls (1.02,.77) and (.88,.75)
    .. (46:1);
\end{scope}
\path[component] \annInteriorCusp \annInteriorRounded;

\fill[black!85] (0,0) circle[radius=1.7pt];

\node[anchor=north,inner sep=2pt]
  at (0,-.27) {$z_0$};

\node[anchor=south]
  at (0,4.22) {$\mathbb S(z_0,4r)$};

\draw[black!60,line width=.55pt]
  (0,4.19)--(0,4.03);

\node[anchor=east,inner sep=2pt]
  at (-1.67,-.62) {$\mathbb S(z_0,r)$};

\draw[black!60,line width=.55pt]
  (-1.63,-.60)--(-.915,-.427);

\node[pathRed,anchor=west]
  at (.53,-2.54) {$\gamma$};

\end{tikzpicture}%
}
\caption{The shaded regions and thick solid curves represent
complementary components of $G$. The dashed curves $\gamma$
represent paths in $\Gamma(z_0,r;G)$.}
\label{fig:annulus-paths}
\end{figure}

Given a domain $G \subset \hatc$, $z_0 \in \C$, and $r>0$, we denote 
$\mathbb{A}(z_0,r):=\D(z_0,4r) \setminus \overline \D(z_0,r)$ and let $\Gamma(z_0,r;G)$ be the family of paths $\gamma:[a,b] \to \hat G$, $a \leq b$, that satisfy $\gamma(a) \in \pi_G(\mathbb{S}(z_0,4r))$, $\gamma(b) \in \pi_G(\mathbb{S}(z_0,r))$, and $\gamma(t)\in\pi_G\bigl(\overline{\mathbb{A}}(z_0,r)\bigr)$ for all $a<t<b$; see Figure \ref{fig:annulus-paths}.

\begin{theorem} \label{thm:keyestimate} 
    	There exists $M=M(\varepsilon,N)>0$ such that if $G \subset \hatc$ is a \emph{finitely connected} $\varepsilon$-UNC-domain that satisfies the packing condition \eqref{eq:packingcondition} with constant $N$, then  
        $$
        \Mod \Gamma(z_0,r;G) \leq M \quad \text{for every } z_0 \in \C \text{ and } r>0.  
        $$
\end{theorem} 
It is important that $M$ does not depend on the number of complementary components of the finitely connected domain $G$ in Theorem \ref{thm:keyestimate}.  

We postpone the proof of Theorem \ref{thm:keyestimate} to Section \ref{sec:key}, and next show how it is applied to prove \eqref{eq:cruuprop1}-\eqref{eq:cruuprop3}. We start with \eqref{eq:cruuprop1} and \eqref{eq:cruuprop2}. We fix a Jordan curve $J \subset \Omega$ and $p \in \mathcal{C}(\Omega)$. 

If $p=p_k \in \mathcal{C}_N(\Omega)$, then we denote $q:=q_k$. If $p \in \mathcal C_P(\Omega)$, then $p \in \Omega_j$ for every $j$, and we may assume, after passing to a subsequence if necessary, that $f_j(p) \to q \in \C$ (our normalization above implies that $q$ is not $\infty$). 

In order to prove \eqref{eq:cruuprop1} and \eqref{eq:cruuprop2}, it suffices to show that if $(a_\ell)$ is a sequence of points in $\Omega$ such that $\dist(a_\ell,p) \to 0$ as $\ell \to \infty$, then 
\begin{equation} \label{eq:pasuo}
\lim_{\ell \to \infty} \limsup_{j \to \infty} \dist(f_j(a_\ell),q)=0. 
\end{equation}

Assume towards a contradiction that \eqref{eq:pasuo} does not hold for some $(a_\ell)$. After passing to a subsequence if necessary, we may assume that $a_\ell \to a \in \partial p$. We denote by $I_\ell$ the segment joining $a_\ell$ and $a$, and let $\Gamma_{j,\ell}$ be the family of paths that join $J$ and $\pi_{\Omega_j}(I_\ell)$ in $\hat \Omega_j \setminus \{p\}$. 

Notice that $\diam(I_\ell) \to 0$ as $\ell \to \infty$. It follows that we can apply Theorem \ref{thm:keyestimate} to $z_0=a$ and a sequence of radii $r_m$, $m=1,2,\ldots$, such that the projections $\pi_{\Omega_j}(\overline{\mathbb{A}}(a,r_m))\setminus\{p\}$ are pairwise disjoint, to conclude that 
\begin{equation} \label{eq:thickannulusestimate} 
\lim_{\ell \to \infty} \limsup_{j \to \infty} \Mod \Gamma_{j,\ell} =0, 
\end{equation}
see \cite{esmayli2024conformal}*{equation (5) on page 563}. 

On the other hand, since \eqref{eq:pasuo} does not hold, a transboundary modulus estimate on circle domains (see \cite{esmayli2024conformal}*{Proposition 3.2}) shows that 
there is $\delta'>0$ such that 
\begin{equation} \label{eq:pasuo3}
\limsup_{\ell \to \infty} \limsup_{j \to \infty} 
\Mod \hat f_j(\Gamma_{j,\ell}) \geq \delta'. 
\end{equation}
Since transboundary modulus is a conformal invariant (see Lemma \ref{lem:confinvariance}), \eqref{eq:pasuo3} contradicts \eqref{eq:thickannulusestimate}, proving \eqref{eq:cruuprop1} and \eqref{eq:cruuprop2}. 

It remains to prove \eqref{eq:cruuprop3}. Suppose towards a contradiction that $\hat{f}(p) \in \mathcal{C}_P(D)$ for some $p=p_k \in \mathcal{C}_N(\Omega)$. We fix a Jordan curve $J \subset \Omega$, and denote by $\Lambda_j$, $j \geq k$, the family of paths that separate $p$ and $J$ in $\hat \Omega_j$. 

Let $I \subset \C$ be a segment whose endpoints lie in $p$ and $J$, respectively, and fix $0<r< \min\{\diam(p),\diam(J)\}/8$. We cover $I$ with a finite number of disks $\mathbb D(z_n,r)$, $n=1,2,\ldots, n_I$, of radius $r$, and notice that if $j \geq k$ then for every $\gamma \in \Lambda_j$ there is at least one $1\leq n \leq n_I$ such that $\gamma$ enters $\mathbb D(z_n,r)$ and exits $\mathbb D(z_n,4r)$. Therefore, Theorem \ref{thm:keyestimate} and the subadditivity of the transboundary modulus (see Proposition \ref{prop:propoftrans}) show that  
\begin{equation} 
\label{eq:praavi1} 
\limsup_{j \to \infty} \Mod \Lambda_j <\infty.  
\end{equation} 

On the other hand, assumption $\hat{f}(p) \in \mathcal{C}_P(D)$ and a transboundary modulus estimate on circle domains (see \cite{esmayli2024conformal}*{Proposition 3.2}) show that 
\begin{equation} \label{eq:praavi2} 
\lim_{j \to \infty} \Mod \hat f_j(\Lambda_j)=\infty. 
\end{equation}
Since transboundary modulus is a conformal invariant by Lemma \ref{lem:confinvariance}, \eqref{eq:praavi1} contradicts \eqref{eq:praavi2}, proving \eqref{eq:cruuprop3}. We have shown how Theorem \ref{thm:mainthm} follows from the key estimate, Theorem \ref{thm:keyestimate}. 



\section{Proof of the key modulus estimate, Theorem \ref{thm:keyestimate}} \label{sec:key} 


\subsection{Quantitatively non-collinear triangles} \label{sec:non-collinear}
    In this section, we establish elementary trigonometric estimates for triangles. Let $\triangle\subset \mathbb{C}$ be a triangle with vertices $a,b,c$, interior angles $\theta_a,\theta_b,\theta_c$, incenter $o$, and inradius $r_\triangle:=\dist(o,\partial \triangle)$. We denote the segment with endpoints $z,w \in \C$ and its length by $zw$ and $|zw|$, respectively. The first lemma can be proved by applying the formula
    $$
    \area(\triangle)= \tfrac{1}{2}|ac|\cdot |bc| \sin \theta_c.  
    $$ 
    
    \begin{lemma}\label{lem:equivcon'}
		Suppose that there exists $0<\varepsilon<1/2$ such that 
		\begin{equation}\label{eq:non-collineardef'}
		\area(\triangle)\ge \varepsilon\, \diam(\triangle)^2. 
        \end{equation}
		 Then 
        \begin{eqnarray}\begin{aligned}
         2\varepsilon \leq \min\{\theta_a,\theta_b,\theta_c\} &\leq 
        \max\{\theta_a,\theta_b,\theta_c\} \leq \pi -2\varepsilon, 
        \quad \text{and } \\  
        \min\{|ab|,|ac|,|bc|\} &\geq 2 \varepsilon \cdot \max\{|ab|,|ac|,|bc|\}. 
        \end{aligned} \label{eq:alternativedef}
        \end{eqnarray}
		\end{lemma}

The next lemma follows from Lemma \ref{lem:equivcon'} and the inradius formula 
\begin{equation} \label{eq:inradiusformula}
r_\triangle = \frac{2\area(\triangle)}{|ab|+|ac|+|bc|}. 
\end{equation}

	\begin{lemma}\label{lem:incenterradius}
		Suppose that there exists $0<\varepsilon<1/2$ such that \eqref{eq:non-collineardef'} holds. Then
		\begin{equation*}
			\tfrac{2}{3}\varepsilon \diam(\triangle) \leq r_\triangle \leq \frac{\min\{|ab|,|ac|,|bc|\}}{2} \leq \diam(\triangle)/2.
		\end{equation*}
	\end{lemma}
The next lemma is important in the proof of key Proposition \ref{prop:blockingprop} below. 

    \begin{lemma}\label{lem:prepareforcosts'}
    		Suppose that there exists $0<\varepsilon<1/2$ such that \eqref{eq:non-collineardef'} holds, and let $a',c' \in ac$. Then 
    	$$
    		|oa'|+|oc'|-|a'c'| \geq \varepsilon^4  \,|ac|.  
        $$
    \end{lemma}

    \begin{proof}
    After interchanging $a'$ and $c'$ if necessary, we may assume
that $|aa'|\leq |ac'|$. By the triangle inequality, we have
    	\begin{eqnarray*} 
        \begin{aligned}
    		|oa'|+|oc'|-|a'c'| &\geq |oa|+|oc|-(|aa'|+|a'c'|+|c'c|)\\
            &=|oa|+|oc|-|ac|. 
        \end{aligned}
        \end{eqnarray*}
    	
    	Moreover, the law of cosines yields  
    	\begin{equation}\label{eq:rightanalysis'}
    		\frac{|oa|+|oc|-|ac|}{|ac|} \geq 
    	\Big(1+ \frac{2|oa||oc|(1+\cos \theta_0)}{|ac|^2} \Big)^{1/2}-1, 
    	\end{equation} 
        where $\theta_0$ is the angle between $oa$ and $oc$. Finally, 
        applying \eqref{eq:alternativedef} and elementary estimates shows that the right-hand side in \eqref{eq:rightanalysis'} is bounded from below by $\varepsilon^4$. The proof is complete. 
    \end{proof}


\subsection{Avoiding complementary components leads to detours} \label{sec:non-collinear2}

In this section, we establish Proposition \ref{prop:blockingprop}, a geometric property of paths on quotient spaces of finitely connected $\varepsilon$-UNC-domains, which will be important for analyzing admissible weights for the path family $\Gamma(z_0,r;G)$ in Theorem \ref{thm:keyestimate}. The packing condition \eqref{eq:packingcondition} will not be used in this section. 

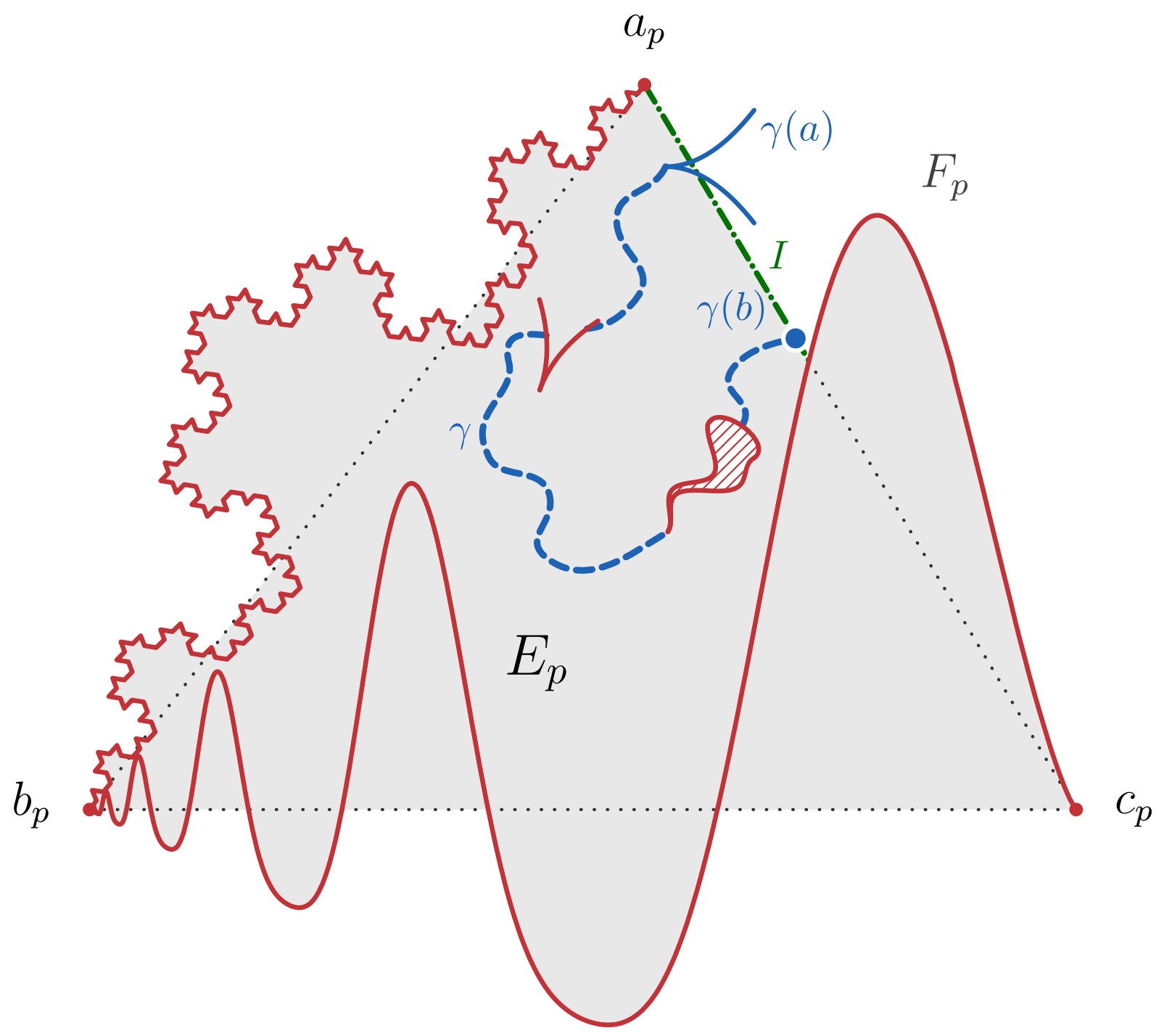
\begin{figure}[htbp]
\centering
\resizebox{0.95\linewidth}{!}{%
\begin{tikzpicture}[
  x=1.35cm,
  y=1.35cm,
  line cap=round,
  line join=round,
  epTravel/.style={
    draw=epBlue,
    line width=1.65pt,
    dash pattern=on 4pt off 2.5pt
  },
  epNonfatPlacement/.style={
    shift={(-0.28,2.72)},
    rotate=70,
    scale=1.2
  },
  epCuspPlacement/.style={
    shift={(0.55,1.80)},
    rotate=55
  },
  declare function={
    epT(\s)=max(0,min(1,1.006*exp(-\s/7)-0.006));
    epX(\s)=-3.2+6.4*epT(\s);
    epY(\s)=7.2*epT(\s)*max(0,1-epT(\s))^0.25*
      (0.725*sin(deg(\s))+0.275*sin(deg(\s))^2);
    epCuspCenter(\u)=0.15*sin(500*\u)*sin(250*\u);
  }
]

\definecolor{epRed}{HTML}{C33236}
\definecolor{epBlue}{HTML}{1D63B5}
\definecolor{epGray}{HTML}{E8E8E8}

\coordinate (A) at (0.4,4.7);
\coordinate (B) at (-3.2,0);
\coordinate (C) at (3.2,0);
\coordinate (Gb) at (1.38,3.055);
\coordinate (Ip) at (1.46620643,2.91029635);

\begin{scope}[epNonfatPlacement]
  \coordinate (NonfatIn) at (0.40,{-0.8*0.40^2});
  \coordinate (NonfatOut) at (0.29,{0.8*0.29^2});
\end{scope}

\begin{scope}[epCuspPlacement]
  \coordinate (CuspIn) at (0,0);
  \coordinate (CuspOut) at (0.84,0);
\end{scope}

\pgfmathsetmacro{\epPhase}{7*ln(1.006/0.006)}

\tikzset{
  epOscillation/.style={
    domain=\epPhase:0,
    samples=1001,
    variable=\s
  }
}

\path[use as bounding box]
  (-4.45,-1.9) rectangle (4.1,5.3);

\def\epSnowflake{%
  (-3.20000,0.00000) -- (-3.15522,0.05777) -- (-3.18283,0.12495) -- (-3.10978,0.11503) -- (-3.06367,0.17178)
  -- (-3.09062,0.23846) -- (-3.16417,0.24876) -- (-3.11742,0.30501) -- (-3.14405,0.37145) -- (-3.07002,0.36078)
  -- (-3.02201,0.41607) -- (-2.99584,0.34999) -- (-2.92151,0.33909) -- (-2.87290,0.39393) -- (-2.89863,0.45967)
  -- (-2.82372,0.44833) -- (-2.77400,0.50231) -- (-2.79919,0.56764) -- (-2.87450,0.57929) -- (-2.82426,0.63287)
  -- (-2.84921,0.69802) -- (-2.92476,0.70985) -- (-2.97500,0.65627) -- (-3.00020,0.72160) -- (-3.07550,0.73325)
  -- (-3.02527,0.78683) -- (-3.05022,0.85198) -- (-2.97455,0.84006) -- (-2.92340,0.89294) -- (-2.94792,0.95776)
  -- (-3.02390,0.96992) -- (-2.97235,1.02250) -- (-2.99669,1.08718) -- (-2.92044,1.07481) -- (-2.86822,1.12688)
  -- (-2.84411,1.06238) -- (-2.76772,1.04990) -- (-2.71523,1.10176) -- (-2.73916,1.16612) -- (-2.66253,1.15346)
  -- (-2.60966,1.20504) -- (-2.58586,1.14077) -- (-2.50916,1.12806) -- (-2.56203,1.07648) -- (-2.53815,1.01216)
  -- (-2.46153,0.99950) -- (-2.40866,1.05108) -- (-2.38486,0.98681) -- (-2.30816,0.97410) -- (-2.25517,1.02558)
  -- (-2.27890,1.08979) -- (-2.20212,1.07701) -- (-2.14906,1.12844) -- (-2.17279,1.19265) -- (-2.24956,1.20542)
  -- (-2.19654,1.25687) -- (-2.22031,1.32111) -- (-2.14361,1.30839) -- (-2.09082,1.36003) -- (-2.06694,1.29570)
  -- (-1.99031,1.28305) -- (-1.93771,1.33482) -- (-1.96175,1.39927) -- (-1.88536,1.38679) -- (-1.83329,1.43897)
  -- (-1.85762,1.50365) -- (-1.93379,1.51595) -- (-1.88205,1.56839) -- (-1.90657,1.63321) -- (-1.98255,1.64537)
  -- (-2.03429,1.59293) -- (-2.05863,1.65761) -- (-2.13479,1.66991) -- (-2.08305,1.72235) -- (-2.10758,1.78717)
  -- (-2.03170,1.77508) -- (-1.98076,1.82813) -- (-2.00571,1.89328) -- (-2.08126,1.90511) -- (-2.03078,1.95851)
  -- (-2.05597,2.02384) -- (-2.13128,2.03549) -- (-2.18176,1.98209) -- (-2.20671,2.04724) -- (-2.28226,2.05907)
  -- (-2.33320,2.00602) -- (-2.30858,1.94113) -- (-2.38456,1.95329) -- (-2.43630,1.90085) -- (-2.46063,1.96553)
  -- (-2.53680,1.97783) -- (-2.48506,2.03027) -- (-2.50958,2.09509) -- (-2.58556,2.10725) -- (-2.63730,2.05481)
  -- (-2.66164,2.11949) -- (-2.73780,2.13179) -- (-2.68606,2.18423) -- (-2.71059,2.24905) -- (-2.63471,2.23696)
  -- (-2.58377,2.29001) -- (-2.60872,2.35516) -- (-2.68427,2.36699) -- (-2.63379,2.42039) -- (-2.65898,2.48572)
  -- (-2.58381,2.47417) -- (-2.53435,2.52836) -- (-2.50876,2.46272) -- (-2.43385,2.45138) -- (-2.38495,2.50599)
  -- (-2.41097,2.57196) -- (-2.33664,2.56105) -- (-2.28893,2.61658) -- (-2.31557,2.68302) -- (-2.38943,2.69356)
  -- (-2.34236,2.74957) -- (-2.36931,2.81625) -- (-2.44286,2.82655) -- (-2.48994,2.77054) -- (-2.51657,2.83698)
  -- (-2.59044,2.84752) -- (-2.54336,2.90353) -- (-2.57031,2.97021) -- (-2.49693,2.96003) -- (-2.45115,3.01704)
  -- (-2.47876,3.08422) -- (-2.55165,3.09402) -- (-2.50654,3.15153) -- (-2.53449,3.21897) -- (-2.46210,3.20956)
  -- (-2.41832,3.26809) -- (-2.38988,3.20027) -- (-2.31782,3.19111) -- (-2.27471,3.25016) -- (-2.30365,3.31836)
  -- (-2.23225,3.30970) -- (-2.19044,3.36974) -- (-2.16101,3.30117) -- (-2.08994,3.29276) -- (-2.13175,3.23272)
  -- (-2.10264,3.16440) -- (-2.03125,3.15574) -- (-1.98943,3.21578) -- (-1.96001,3.14721) -- (-1.88893,3.13880)
  -- (-1.84775,3.19933) -- (-1.87764,3.26825) -- (-1.80717,3.26031) -- (-1.76718,3.32175) -- (-1.79765,3.39112)
  -- (-1.86768,3.39873) -- (-1.82825,3.46059) -- (-1.85898,3.53017) -- (-1.78935,3.52286) -- (-1.75094,3.58551)
  -- (-1.71983,3.51565) -- (-1.65043,3.50853) -- (-1.61248,3.57153) -- (-1.64392,3.64164) -- (-1.57495,3.63485)
  -- (-1.53780,3.69846) -- (-1.50608,3.62813) -- (-1.43730,3.62148) -- (-1.47445,3.55787) -- (-1.44292,3.48768)
  -- (-1.37395,3.48089) -- (-1.33680,3.54450) -- (-1.30508,3.47417) -- (-1.23630,3.46752) -- (-1.27345,3.40391)
  -- (-1.24191,3.33372) -- (-1.31098,3.34059) -- (-1.34893,3.27759) -- (-1.31782,3.20773) -- (-1.24843,3.20061)
  -- (-1.28684,3.13796) -- (-1.25597,3.06829) -- (-1.18634,3.06098) -- (-1.14793,3.12363) -- (-1.11682,3.05377)
  -- (-1.04742,3.04665) -- (-1.00947,3.10965) -- (-1.04091,3.17976) -- (-0.97194,3.17297) -- (-0.93479,3.23658)
  -- (-0.90307,3.16625) -- (-0.83429,3.15960) -- (-0.87144,3.09599) -- (-0.83991,3.02580) -- (-0.77094,3.01901)
  -- (-0.73379,3.08262) -- (-0.70207,3.01229) -- (-0.63329,3.00564) -- (-0.59647,3.06950) -- (-0.62842,3.14001)
  -- (-0.55994,3.13359) -- (-0.52365,3.19786) -- (-0.55583,3.26854) -- (-0.62415,3.27484) -- (-0.58805,3.33926)
  -- (-0.62031,3.41000) -- (-0.55210,3.40379) -- (-0.51623,3.46838) -- (-0.48390,3.39759) -- (-0.41573,3.39140)
  -- (-0.37990,3.45603) -- (-0.41224,3.52683) -- (-0.34406,3.52064) -- (-0.30816,3.58521) -- (-0.34042,3.65595)
  -- (-0.40866,3.66219) -- (-0.37264,3.72666) -- (-0.40482,3.79735) -- (-0.47314,3.80364) -- (-0.50916,3.73917)
  -- (-0.54142,3.80991) -- (-0.60966,3.81615) -- (-0.57364,3.88062) -- (-0.60582,3.95131) -- (-0.53745,3.94497)
  -- (-0.50105,4.00916) -- (-0.53300,4.07967) -- (-0.60155,4.08614) -- (-0.56488,4.15012) -- (-0.59669,4.22051)
  -- (-0.52790,4.21386) -- (-0.49057,4.27733) -- (-0.45903,4.20714) -- (-0.39007,4.20035) -- (-0.35233,4.26351)
  -- (-0.38355,4.33346) -- (-0.31415,4.32634) -- (-0.27550,4.38880) -- (-0.24464,4.31913) -- (-0.17500,4.31182)
  -- (-0.21365,4.24936) -- (-0.18255,4.17950) -- (-0.11315,4.17238) -- (-0.07450,4.23484) -- (-0.04364,4.16517)
  -- (0.02600,4.15786) -- (0.06517,4.21993) -- (0.03470,4.28930) -- (0.10487,4.28158) -- (0.14516,4.34280)
  -- (0.11527,4.41172) -- (0.04465,4.41978) -- (0.08553,4.48053) -- (0.05595,4.54922) -- (0.12703,4.54082)
  -- (0.16916,4.60061) -- (0.19826,4.53229) -- (0.26966,4.52363) -- (0.31245,4.58293) -- (0.28383,4.65088)
  -- (0.35589,4.64172) -- (0.40000,4.70000)%
}

\fill[epGray] (A)--(B)--(C)--cycle;

\path[fill=epGray]
  \epSnowflake -- cycle;

\path[fill=epGray]
  (B) --
  plot[epOscillation] ({epX(\s)},{epY(\s)})
  -- cycle;

\draw[
  black!80,
  line width=0.8pt,
  dash pattern=on 0pt off 3.2pt
]
  (A)--(B)--(C) (Ip)--(C);

\draw[
  green!45!black,
  line width=1.4pt,
  dash pattern=on 5pt off 2pt on 0pt off 2pt
]
  (A)--(Ip)
  node[midway,right=6pt,yshift=-8pt,font=\small] {$I$};

\draw[epRed,line width=1.18pt]
  \epSnowflake;

\draw[epRed,line width=1.18pt]
  plot[epOscillation] ({epX(\s)},{epY(\s)});

\draw[epTravel]
  (0.5369078,4.1685834)
  .. controls (0.45,3.98) and (0.20,4.02) .. (0.22,3.80)
  .. controls (0.24,3.58) and (0.52,3.50) .. (0.32,3.28)
  .. controls (0.24,3.19) and (0.15,3.13) .. (NonfatIn);

\draw[epTravel]
  (NonfatOut)
  .. controls (-0.37,3.09) and (-0.53,3.10) .. (-0.48,2.90)
  .. controls (-0.43,2.69) and (-0.75,2.58) .. (-0.62,2.31)
  .. controls (-0.56,2.12) and (-0.28,2.24) .. (-0.22,2.05)
  .. controls (-0.16,1.86) and (-0.42,1.76) .. (-0.22,1.62)
  .. controls (0.04,1.44) and (0.30,1.65) .. (CuspIn);

\draw[epTravel]
  (CuspOut)
  .. controls (1.14,2.66) and (0.84,2.70) .. (0.98,2.87)
  .. controls (1.08,2.99) and (1.25,3.03) .. (Gb);

\begin{scope}[epNonfatPlacement]
\draw[epRed,line width=1.18pt]
  plot[domain=0:0.46,samples=101,variable=\u]
    (\u,{0.8*\u^2});

\draw[epRed,line width=1.18pt]
  plot[domain=0:0.46,samples=101,variable=\u]
    (\u,{-0.8*\u^2});
\end{scope}

\begin{scope}[epCuspPlacement]
\def\epTwistedCusp{%
  plot[domain=0:0.72,samples=151,variable=\u]
    (\u,{epCuspCenter(\u)+0.36*\u^2})
  .. controls (0.88,0.225) and (0.88,-0.225)
  .. (0.72,-0.186624)
  -- plot[domain=0.72:0,samples=151,variable=\u]
    (\u,{epCuspCenter(\u)-0.36*\u^2})
  -- cycle
}

\path[fill=white]
  \epTwistedCusp;

\path[
  draw=epRed,
  line width=1.18pt,
  pattern=north east lines,
  pattern color=epRed
]
  \epTwistedCusp;
\end{scope}

\begin{scope}[shift={(0.5369078,4.1685834)},scale=1.8]
\draw[epBlue,line width=1.15pt]
  plot[domain=0:0.32,samples=101,variable=\u]
    (\u,{2*\u^2});

\draw[epBlue,line width=1.15pt]
  plot[domain=0:0.32,samples=101,variable=\u]
    (\u,{-2*\u^2});
\end{scope}

\fill[white] (Gb) circle[radius=3.6pt];
\fill[epBlue] (Gb) circle[radius=2.5pt];

\node[epBlue,font=\small,anchor=east]
  at (-0.62,2.42) {$\gamma$};

\node[epBlue,font=\small,anchor=west]
  at (1.05,4.4035834) {$\gamma(a)$};

\node[epBlue,font=\small,anchor=east]
  at (1.30,3.24) {$\gamma(b)$};

\foreach \P in {A,B,C}
  \fill[epRed] (\P) circle[radius=1.7pt];

\node[font=\large,anchor=south]
  at (0.4,4.83) {$a_p$};

\node[font=\large,anchor=east]
  at (-3.35,0.03) {$b_p$};

\node[font=\large,anchor=west]
  at (3.35,0) {$c_p$};

\node[font=\Large]
  at (-0.30,0.95) {$E_p$};

\node[font=\large,text=black!75]
  at (2.35,4.1) {$F_p$};

\end{tikzpicture}%
}
\caption{The solid curves connecting $a_p$, $b_p$, and $c_p$ represent $p\in\mathcal{C}_N(G)$. The two smaller sets along $\gamma$ represent other complementary components.
The gray region represents $E_p$, while the thick dashed path and
the dash-dotted segment represent $\gamma$ and $I$, respectively,
as in Lemma \ref{lem:sameside} and its proof.}
\label{fig:Ep}
\end{figure}

Suppose that $G \subset \hatc$ is a domain that contains $\infty$. We fix $p \in \mathcal{C}_N(G)$ and a non-degenerate solid triangle $\triangle_p$ with vertices $a_p,b_p,c_p \in p$. We denote $A_p := p \cup \triangle_p$, let $F_p$ be the connected component
	 of $\hatc\setminus A_p$ that contains $\infty$, and define (see Figure \ref{fig:Ep})  
	 \begin{equation} \label{eq:eepee}
     E_p := \hatc\setminus \overline F_p. 
\end{equation}

	 \begin{lemma}\label{lem:sameside}
	 	Suppose that $\gamma:[a,b] \to \hat G \setminus \{p\}$ is a path such that 
\begin{itemize}
    \item[(i)] both $\gamma(a)$ and $\gamma(b)$ intersect  
        $\partial \triangle_p \cap \partial E_p$ and 
    \item[(ii)] $\gamma(t) \subset E_p$ for all $a<t<b$.  
\end{itemize}
        Then there is an edge $e$ of $\triangle_p$ such that $e \cap \gamma(a)\neq \emptyset$ and $e \cap \gamma(b) \neq \emptyset$.  
\end{lemma}
\begin{proof}
There exists a unique component $B_p$ of $E_p \setminus p$ such that $\gamma(t) \in \pi_G(B_p)$ for all $a<t<b$. It follows that there exists a relative component $I$ of $\partial \triangle_p \setminus p$ such that 
\begin{equation} \label{eq:zepe1}
I \subset (\partial B_p)\setminus p \subset (\partial E_p) \setminus p \quad \text{and} \quad \gamma(a) \cap I \neq \emptyset. 
\end{equation}

Note that $I$ is a segment whose endpoints are in $p$. We apply the \emph{crosscut separation theorem} (\cite{newman1939elements}*{Theorem 11.7}) to $I$ and the simply connected domain $U_p:=\hatc \setminus p$ to show that $U_p \setminus I$ has exactly two components: a bounded component $E'_p$ and an unbounded component $F'_p$.

Since $\partial E'_p \subset \overline E_p$, we have $F_p \subset F'_p$, and hence
$E'_p \subset E_p\setminus p$. Together with $\partial E'_p\setminus p=I\subset\partial E_p$, it follows that $E'_p$ is a component of $E_p\setminus p$. Moreover, since  $B_p$ and $E'_p$ lie on the same side of $I$ and both are components of $E_p\setminus p$, we deduce that $B_p=E'_p$.
In particular, since 
$$
\emptyset \neq \gamma(b) \setminus E_p \subset \gamma(b) \setminus E'_p 
\quad \text{and} \quad p \notin |\gamma| 
$$ 
by assumption, there exists $a< c \leq b$ such that 
\begin{equation*} \label{eq:zepe2}
\emptyset \neq \gamma(c) \cap (\partial E_p' \setminus p)=\gamma(c) \cap I. 
\end{equation*}
On the other hand, since $I \subset \partial E_p$ by \eqref{eq:zepe1}, it follows that
$\gamma(c) \cap \partial E_p \neq \emptyset$. Combining with assumption (ii), we conclude that $c=b$. Thus the edge of $\triangle_p$ containing $I$ intersects
both $\gamma(a)$ and $\gamma(b)$. The proof is complete.
\end{proof}

Let $G$ be a finitely connected $\varepsilon$-UNC-domain, $0<\varepsilon<\nicefrac{1}{10}$, and fix $p \in \mathcal{C}_N(G)$. We choose a triangle $\triangle_p$ with vertices in $p$ satisfying \eqref{eq:non-collineardef}, and denote its incenter and inradius by $o_p$ and $r_p$. Let $\Gamma_p$ be the family of paths 
	 $$ 
	 	\alpha:[s_1,t_1] \to \hat{G} \setminus \{p\},  
	 $$
	 for which there are $s_1<s_2 \leq t_2 < t_1 $ with the following properties:  
	 \begin{itemize}
	 	\item[(i)] $\alpha(s_2) \cup \alpha(t_2) \subset \mathbb{D}(o_p, \varepsilon^5 \, r_p)$, 
	 	\item[(ii)] $\alpha(t) \cap \partial E_p \neq \emptyset$ for $t=s_1$ and $t=t_1$, 
        \item[(iii)] $\alpha(t) \subset E_p$ for all $s_1<t<t_1$;  
	 \end{itemize} 
	 here $E_p$ is the set in \eqref{eq:eepee}. Note that since $p \notin |\alpha|$, the definition of $E_p$ and the second condition imply that 
$$ 
\alpha(t) \cap \partial \Delta_p \neq \emptyset \text{ for }
t=s_1 \text{ and }t=t_1. 
$$
	 
\begin{proposition}\label{prop:blockingprop}
	Suppose that $\alpha \in \Gamma_p$ satisfies 
\begin{equation} \label{eq:duidi}
\diam (\alpha(t)) \le \varepsilon^5\,r_p \quad \text{for } t=s_1 \text{ and } t=t_1. 
\end{equation}
Then 
\begin{equation} \label{eq:dessed} 
\dist(\alpha(s_1),\alpha(t_1)) \leq \dist(\alpha(s_1),\alpha(s_{2}))+\dist(\alpha(t_1),\alpha(t_2))  - \varepsilon^5 \, r_p.  
\end{equation} 
\end{proposition}
		
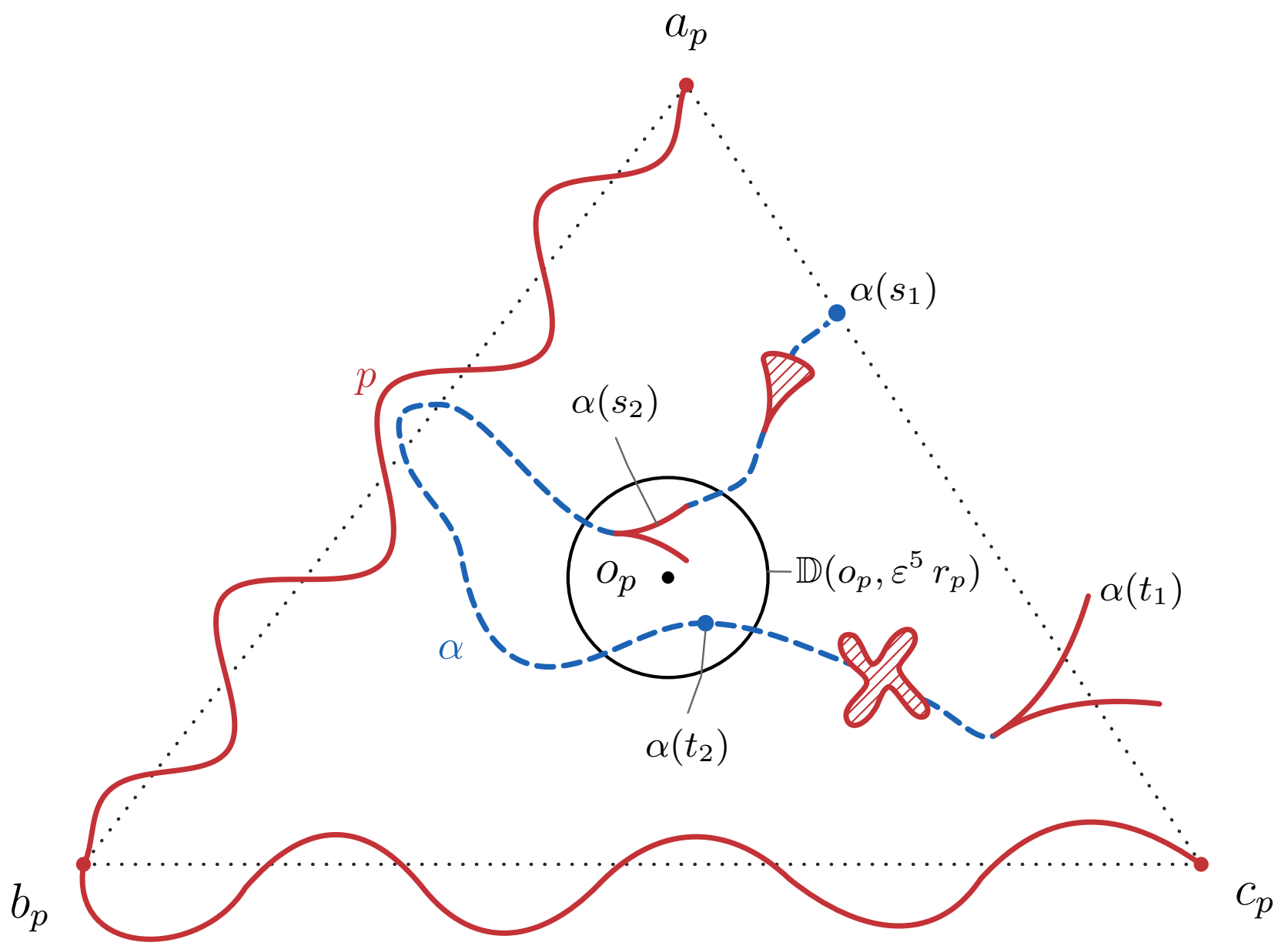
\begin{figure}[htbp]
\centering
\resizebox{0.88\linewidth}{!}{%
\begin{tikzpicture}[
  x=1.25cm,y=1.25cm,
  line cap=round,line join=round,
  every node/.style={font=\large},
  boundary/.style={
    draw=detourRed,
    solid,
    line width=1.30pt
  },
  travel/.style={
    draw=detourBlue,
    line width=1.35pt,
    dash pattern=on 4pt off 2.5pt
  },
  triangle/.style={
    draw=black!85,
    line width=0.80pt,
    dash pattern=on 0pt off 3.2pt
  },
  disk/.style={
    draw=black,
    solid,
    line width=0.80pt
  },
  leader/.style={
    draw=black!60,
    solid,
    line width=0.45pt
  }
]

\definecolor{detourRed}{HTML}{C33236}
\definecolor{detourBlue}{HTML}{1D63B5}

\coordinate (A) at (0.3,5.3);
\coordinate (B) at (-3.8,0);
\coordinate (C) at (3.8,0);
\coordinate (O) at (0.17468421,1.95040471);
\coordinate (Sone) at (1.32358491,3.75);
\coordinate (Ttwo) at (0.43,1.64);

\draw[triangle]
  (A)--(B)--(C)--cycle;

\draw[disk]
  (O) circle[radius=0.68];

\draw[boundary]
  plot[domain=0:1,samples=401,variable=\t]
  ({0.3-4.1*\t+5.3/sqrt(4.1^2+5.3^2)*
    (0.38+0.035*sin(540*\t+20))*
    sin(1440*\t+18*sin(360*\t))*sqrt(max(0,sin(180*\t)))},
   {5.3*(1-\t)-4.1/sqrt(4.1^2+5.3^2)*
    (0.38+0.035*sin(540*\t+20))*
    sin(1440*\t+18*sin(360*\t))*sqrt(max(0,sin(180*\t)))});

\draw[boundary]
  (B)
  .. controls (-3.90,-0.65) and (-3.00,-0.65) .. (-2.70,-0.16)
  .. controls (-2.30,0.30) and (-1.90,0.34) .. (-1.50,-0.15)
  .. controls (-1.15,-0.60) and (-0.70,-0.55) .. (-0.30,-0.14)
  .. controls (0.12,0.30) and (0.60,0.28) .. (1.00,-0.10)
  .. controls (1.50,-0.50) and (2.00,-0.55) .. (2.30,-0.10)
  .. controls (2.80,0.45) and (3.30,0.35) .. (C);

\draw[travel]
  (Sone)
  .. controls (1.21,3.60) and (1.05,3.57)
  .. (1.01127068,3.44803709);

\draw[travel]
  (0.83,2.95)
  .. controls (0.75,2.76) and (0.78,2.63) .. (0.59,2.55)
  .. controls (0.49,2.51) and (0.40,2.47) .. (0.30,2.43432);

\draw[travel]
  (-0.18,2.25)
  .. controls (-0.62,2.30) and (-1.051,3.132) .. (-1.381,3.127)
  -- plot[domain=0:1,samples=101,variable=\v]
    ({0.3-4.1*(0.41+0.06*\v)
       -5.3/sqrt(4.1^2+5.3^2)*0.16*sin(180*\v)},
     {5.3*(1-0.41-0.06*\v)
       +4.1/sqrt(4.1^2+5.3^2)*0.16*sin(180*\v)})
  .. controls (-1.561,2.537) and (-1.27,2.37) .. (-1.23,2.09)
  .. controls (-1.15,1.62) and (-0.94,1.29) .. (-0.54,1.35)
  .. controls (-0.20,1.40) and (0.07,1.65) .. (Ttwo)
  .. controls (0.82,1.63) and (1.14,1.46) .. (1.42,1.36);

\draw[travel]
  (1.908,1.12)
  .. controls (2.08,1.06) and (2.28,0.80)
  .. (2.40,0.88);

\begin{scope}[shift={(0.83,2.95)},rotate=70]
\path[
  boundary,
  pattern=north east lines,
  pattern color=detourRed
]
  plot[domain=0:0.44,samples=61,variable=\u]
    (\u,{0.85*\u^2})
  .. controls (0.56,0.25432) and (0.56,-0.25432)
  .. (0.44,-0.16456)
  -- plot[domain=0.44:0,samples=61,variable=\u]
    (\u,{-0.85*\u^2})
  -- cycle;
\end{scope}

\draw[boundary]
  plot[domain=0:0.48,samples=61,variable=\u]
    ({-0.18+\u},{2.25+0.8*\u^2});

\draw[boundary]
  plot[domain=0:0.48,samples=61,variable=\u]
    ({-0.18+\u},{2.25-0.8*\u^2});

\begin{scope}[
  shift={(2.40,0.88)},
  rotate=33.44,
  scale=2.5
]
\draw[boundary]
  plot[domain=0:0.42,samples=61,variable=\u]
    (\u,{\u^2});

\draw[boundary]
  plot[domain=0:0.42,samples=61,variable=\u]
    (\u,{-\u^2});
\end{scope}

\begin{scope}[shift={(1.62,1.28)},scale=0.8]
\path[
  boundary,
  pattern=north east lines,
  pattern color=detourRed
]
  (-0.22,0.34)
  .. controls (-0.35,0.40) and (-0.40,0.24) .. (-0.25,0.10)
  .. controls (-0.19,0.04) and (-0.10,-0.01) .. (-0.12,-0.07)
  .. controls (-0.15,-0.15) and (-0.32,-0.25) .. (-0.29,-0.36)
  .. controls (-0.26,-0.47) and (-0.12,-0.42) .. (-0.03,-0.25)
  .. controls (0.00,-0.18) and (0.02,-0.11) .. (0.07,-0.09)
  .. controls (0.12,-0.10) and (0.18,-0.22) .. (0.27,-0.31)
  .. controls (0.36,-0.40) and (0.46,-0.32) .. (0.36,-0.20)
  .. controls (0.28,-0.11) and (0.13,0.03) .. (0.16,0.09)
  .. controls (0.20,0.17) and (0.32,0.26) .. (0.29,0.35)
  .. controls (0.26,0.44) and (0.13,0.40) .. (0.04,0.23)
  .. controls (0.01,0.17) and (-0.02,0.11) .. (-0.07,0.10)
  .. controls (-0.13,0.15) and (-0.15,0.27) .. (-0.22,0.34)
  -- cycle;
\end{scope}

\fill[white] (Sone) circle[radius=3pt];
\fill[detourBlue] (Sone) circle[radius=2.1pt];
\fill[black] (O) circle[radius=1.5pt];
\fill[detourBlue] (Ttwo) circle[radius=2pt];

\node[anchor=east]
  at (0.08,1.94) {$o_p$};

\node[black,font=\footnotesize,anchor=west]
  at (1.30,3.90) {$\alpha(s_1)$};

\node[black,font=\footnotesize,anchor=west]
  at (3.00,1.86) {$\alpha(t_1)$};

\node[black,font=\footnotesize,anchor=south]
  (Ls) at (-0.18,2.90) {$\alpha(s_2)$};

\draw[leader]
  (Ls.south)--(-0.10,2.71)--(0.10,2.313);

\node[black,font=\footnotesize,anchor=north]
  (Lt) at (0.31,1.03) {$\alpha(t_2)$};

\draw[leader]
  (Lt.north)--(0.40,1.28)--(0.43,1.59);

\node[black,font=\footnotesize,anchor=west,inner sep=1pt]
  (DiskLabel) at (1.01,1.99)
  {$\mathbb{D}(o_p,\varepsilon^5\,r_p)$};

\draw[leader]
  (0.8534,1.992)--(DiskLabel.west);

\node[detourBlue,font=\small,anchor=west]
  at (-1.5,1.46) {$\alpha$};

\node[detourRed,font=\small,anchor=east]
  at (-1.70,3.28) {$p$};

\foreach \P in {A,B,C}
  \fill[detourRed] (\P) circle[radius=1.8pt];

\node[anchor=south]
  at (0.3,5.43) {$a_p$};

\node[anchor=north east]
  at (-3.92,-0.02) {$b_p$};

\node[anchor=north west]
  at (3.92,-0.02) {$c_p$};

\end{tikzpicture}%
}
\caption{The thick solid curves and regions represent
complementary components in $\mathcal{C}_N(G)$, while the thick
dashed path represents $\alpha\in\Gamma_p$. 
The disk $\mathbb{D}(o_p,\varepsilon^5\,r_p)$ and some
complementary components are enlarged for clarity.}
\label{fig:detour}
\end{figure}

The geometric meaning of \eqref{eq:dessed} is that the paths that ``avoid'' $p$ but travel near $\pi_G(o_p)$ must take a ``detour'' whose length is comparable to the diameter of $p$; see Figure \ref{fig:detour}. 

\begin{proof}
By Lemma \ref{lem:sameside}, there is an edge $e$ of $\triangle_p$ that both $\alpha(s_1)$ and $\alpha(t_1)$ intersect. We fix $a_p' \in \alpha(s_1) \cap e$ and $c_p' \in \alpha(t_1) \cap e$. Since $r_p \leq \ell(e)$ by Lemma \ref{lem:incenterradius}, Lemma \ref{lem:prepareforcosts'} and \eqref{eq:duidi} yield 
\begin{eqnarray}
\begin{aligned}
& \varepsilon^4\, r_p \leq \varepsilon^4\, \ell(e) \leq |a'_po_p|+|c'_po_p| - 
|a_p'c_p'|  \leq \dist(\alpha(s_1),o_p) \\& +\dist(o_p,\alpha(t_1))+2\varepsilon^5\,r_p-\dist(\alpha(s_1),\alpha(t_1)). 
\end{aligned} \label{eq:sepultura1}
\end{eqnarray}
Since $\alpha(s_2) \cup \alpha(t_2) \subset  \mathbb{D}(o_p,\varepsilon^{5} \, r_p)$, we also have 
\begin{eqnarray}
\begin{aligned}
& \dist(\alpha(s_1),o_p)+\dist(\alpha(t_1),o_p) \leq 
\dist(\alpha(s_1),\alpha(s_2)) \\ & +\dist(\alpha(t_1),\alpha(t_2))+ 
2\varepsilon^5\, r_p. 
\end{aligned}
\label{eq:sepultura2}
\end{eqnarray}
Since $0<\varepsilon <\nicefrac{1}{10}$, the claim follows by combining \eqref{eq:sepultura1} and \eqref{eq:sepultura2}. 
\end{proof}


\subsection{Construction of admissible function \texorpdfstring{$\rho$}{rho}} \label{sec:rhodefinition}

We now start the proof of Theorem \ref{thm:keyestimate}. Suppose that $G \subset \hatc$ is a finitely connected $\varepsilon$-UNC-domain for which the packing condition \eqref{eq:packingcondition} holds with constant $N$. Note that the UNC-condition is not affected if we decrease $\varepsilon$. Therefore, we lose no generality by assuming $\varepsilon< \nicefrac{1}{100}$. 

Our goal is to achieve an upper bound on the transboundary modulus of the paths that connect the projections of the boundary components of the annuli $\mathbb{A}(z_0,r)=\mathbb{D}(z_0,4r)\setminus \overline \D(z_0,r)$. Since the claim is not affected by scalings and translations, we lose no generality by assuming $z_0=0$ and $r=1$. 
Moreover, every $\gamma_1 \in \Gamma(z_0,1;G)$ has a subpath $\gamma$ that satisfies \eqref{eq:noutopoyta} below; note that we may have $\gamma:\{a\} \to \hat G$. Therefore, we conclude using the monotonicity of the transboundary modulus (see Proposition \ref{prop:propoftrans}) that Theorem \ref{thm:keyestimate} is true if we can prove the following theorem. 

\begin{theorem}\label{thm:keyestimate'} 
	There is an $M>0$ that depends only on $\varepsilon$ and $N$ such that $\Mod\Gamma \le M$, where $\Gamma$ is the family of paths $\gamma:[a,b] \to \hat G$ for which  
    \begin{equation}\label{eq:noutopoyta} 
    \gamma(a) \in \pi_{G}({\mathbb{S}}(0,4)),\, \gamma(b) \in \pi_{G}({\mathbb{S}}(0,1)), \text{ and }\gamma(t) \subset {\mathbb{A}}(0,1), \, a<t<b. 
    \end{equation}
\end{theorem} 

We start the proof of Theorem \ref{thm:keyestimate'} by constructing an admissible function $\rho: \hat G \to [0,\infty]$ for $\Gamma$ whose ``energy'' satisfies 
\begin{equation} \label{eq:energyestimate}
\int_{G} \rho^2 \, dA + \sum_{p \in \mathcal{C}(G)} \rho(p)^2 \leq M(\varepsilon,N). 
\end{equation}
The main challenge will be to establish the admissibility of $\rho$ for $\Gamma$.  

Since $G$ is finitely connected, part (4) of Proposition \ref{prop:propoftrans} shows that we lose no generality by assuming that 
all complementary components of $G$ are non-trivial, i.e., $\mathcal{C}(G)=\mathcal{C}_N(G)$. The definition of $\rho$ in $G$ is simple; we set 
\begin{eqnarray} \label{eq:rhodef1} \rho(z):=\left\{\begin{array}{ll}  
1, & z \in G \cap \mathbb{A}(0,1), \\ 
0, & z \in G \setminus \mathbb{A}(0,1). 
\end{array} \right.
\end{eqnarray}
We also set 
\begin{eqnarray} \label{eq:rhodef2} \rho(p):=\left\{\begin{array}{ll}  
0, & p \in \mathcal{C}(G), \, p \cap \mathbb{A}(0,1)=\emptyset, \\ 
1, & p \in \mathcal{C}_L,  \end{array} \right.
\end{eqnarray}
where 
$$
\mathcal{C}_L:=\{p \in \mathcal{C}(G): \, p \cap \mathbb{A}(0,1) \neq \emptyset, \, \,\diam(p) \geq \varepsilon^{10}\}. 
$$
The packing condition \eqref{eq:packingcondition} will guarantee that 
$\sum_{p \in \mathcal{C}_L}\rho(p)^2=\card(\mathcal{C}_L)$ is uniformly bounded. 

The definition of $\rho$ in the remaining set   
$$
\mathcal{C}_S :=\{p \in \mathcal{C}(G)\setminus \mathcal{C}_L: \, p \cap \mathbb{A}(0,1) \neq \emptyset\} 
$$
is more involved. The natural definition in terms of admissibility would be $\rho(p)=\diam(p)$. However, our assumptions do not imply uniform upper bounds on $\sum_{p \in \mathcal{C}_S} \diam(p)^2$; see \cite{esmayli2024conformal}*{Example 6.3} and the proof of Theorem \ref{thm:counterexample}. To overcome the problem, we identify a subcollection $\mathcal B \subset \mathcal C_S$ of components that are ``buried'' under some of the ``exposed'' components $p \in \mathcal E$ for which $\sum_{p \in \mathcal{E}} \diam(p)^2$ is uniformly bounded. We set $\rho=0$ on $\mathcal B$ and $\rho=C\diam(p)$ on $\mathcal E$, where $C=C(\varepsilon)>0$, and use Proposition \ref{prop:blockingprop} to prove that $\rho$ is admissible for $\Gamma$; compare to \cite{esmayli2024conformal}*{Section 4.2}. 

Given $p \in \mathcal{C}(G)$, fix a triangle $\triangle_p$ with vertices $a_p,b_p,c_p$, incenter $o_p$, and inradius $r_p$ as in Proposition \ref{prop:blockingprop}, so that \eqref{eq:non-collineardef} is satisfied.  
If $\mathcal C_S=\emptyset$, we set $\mathcal E=\mathcal B=\emptyset$.
Otherwise, we first choose $p_1\in\mathcal C_S$ such that
    $$r_{p_1}=\max_{p\in\mathcal C_S}r_p.$$
We denote $\triangle_1:=\triangle_{p_1}$, $r_1:=r_{p_1}$, $o_1:=o_{p_1}$, and let 
\begin{eqnarray*}
	\mathcal{E}_1&:=&\left\{p \in \mathcal{C}_S: r_p \geq \varepsilon^{10} r_1,\, |o_p-o_1| \leq \tfrac{r_1}{\varepsilon^{3}}\right\}, \text{ and } \\ 
	\mathcal{B}_1&:=&\left\{p \in \mathcal{C}_S: r_p < \varepsilon^{10} r_1, \,|o_p-o_1| \leq \varepsilon^6 r_1 \right\}.   
\end{eqnarray*}
Note that $p_1 \in \mathcal{E}_1$. The set $\mathcal{E}_1$ consists of the elements $p \in \mathcal{C}_S$ whose diameters are not too small compared to $p_1$ and which lie at a controlled distance not too far from $o_1$. The set $\mathcal{B}_1$ comprises those $p\in \mathcal{C}_S$ that are significantly smaller than $p_1$ and are close to $o_1$. 

Suppose now that $p_\ell \in \mathcal{C}_S$ and $\mathcal{E}_\ell, \mathcal{B}_\ell \subset \mathcal{C}_S$ are chosen for $1 \leq \ell \leq k-1$. We stop the process if 
$$\mathcal{C}_S^{k}:=\mathcal{C}_S \setminus \bigcup_{\ell=1}^{k-1} (\mathcal{E}_\ell \cup \mathcal{B}_\ell) =\emptyset. 
$$ 
Otherwise, we choose $p_{k} \in \mathcal{C}_S^k$ such that 
\begin{equation} \label{eq:rdef}
r_{p_{k}}=\max_{p \in \mathcal{C}_S^k} r_p.  
\end{equation}
We denote $\triangle_k:=\triangle_{p_k}$, $E_k:=E_{p_k}$,  $r_{k}:=r_{p_{k}}$ and $o_{k}:=o_{p_{k}}$, and let 
\begin{equation}\label{eq:edef}
	\mathcal{E}_{k} := \left\{p \in \mathcal{C}_S^k: \, r_p \geq \varepsilon^{10} r_{k},\, |o_p-o_{k}| \leq \tfrac{r_{k}}{\varepsilon^{3}}  \right\}, 
\end{equation} 
\text{ and }
\begin{equation}\label{eq:bdef}
	\mathcal{B}_{k} := \left\{p \in \mathcal{C}_S^k:\, r_p < \varepsilon^{10} r_{k}, \,|o_p-o_{k}|\leq \varepsilon^6 r_{k} \right\}. 
\end{equation}
Since $p_{k} \in \mathcal{E}_{k}$ for each $k\in \mathbb{N}$ and $G$ is finitely connected, the process stops after $L \leq \card(\mathcal{C}_S)$ steps and we have a partition of $\mathcal{C}_S$ into disjoint sets $\mathcal{E}_k$ and  $\mathcal{B}_k$, $k=1,\ldots, L$, so,
$$
	\mathcal{C}_S= \mathcal{E} \cup \mathcal{B}, \quad \text{where} \quad \mathcal{E}:= \bigcup_{k=1}^L \mathcal{E}_k, \quad \mathcal{B}:= \bigcup_{k=1}^L \mathcal{B}_k.
$$
Note that $\mathcal{E}_k$ comes with a distinguished $p_k\in \mathcal{E}_k$ and the associated $o_k$ and $r_k$. We now set 
\begin{eqnarray} \label{eq:rhodef3}
\rho(p):=\left\{\begin{array}{ll} \frac{\diam(p)}{\varepsilon^{15}}, & p \in \mathcal{E}, \\ 0, & p \in \mathcal{B}. \end{array}
\right. 
\end{eqnarray}
We have defined $\rho:\hat G \to [0,\infty]$. To prove Theorem \ref{thm:keyestimate}, it suffices to show that $\rho$ is admissible for $\Gamma$ and satisfies the energy estimate \eqref{eq:energyestimate}. 


\subsection{Energy of \texorpdfstring{$\rho$}{rho}} 
In this section, we show that the function $\rho$ defined in Section \ref{sec:rhodefinition} satisfies the energy estimate \eqref{eq:energyestimate}. We note that this is the only part of the proof of Theorem \ref{thm:keyestimate} that uses the packing condition \eqref{eq:packingcondition}. 

\begin{proposition}\label{prop:energy} We have 
$$
\int_{G} \rho^2 \, dA + \sum_{p \in \mathcal{C}(G)} \rho(p)^2 
\leq \varepsilon^{-75}N.  
$$  
\end{proposition}

We will need the following lemmas to prove Proposition \ref{prop:energy}. 
\begin{lemma}\label{lem:inter1} The disks 
	$\mathbb{D}(o_k,\varepsilon^{6} r_k/2)$, $1 \leq k \leq L$, 
	are pairwise disjoint. 
\end{lemma} 
\begin{proof}
Let $1 \leq k < m \leq L$. By \eqref{eq:rdef}, \eqref{eq:edef}, and \eqref{eq:bdef}, $r_m \leq r_k$ and 
$$
|o_k-o_m| > \varepsilon^6 r_k = \varepsilon^6(2r_k)/2\geq  \varepsilon^6(r_m+r_k)/2.  
$$ 
\end{proof}

\begin{lemma} \label{lem:voivoi} 
All $p \in \mathcal{C}(G)$ satisfy $E_p \subset \D(o_p,3r_p/\varepsilon)$, 
$$ 
\dist(o_p,p) \leq \diam(p), \quad \text{and} \quad r_p \leq \diam(p) \leq \frac{3r_p}{2\varepsilon}. 
$$
\end{lemma} 
\begin{proof} The second claim holds because the distance between $o_p$ and $a_p \in p$ is smaller than $\diam(\Delta_p) \leq \diam(p)$. Similarly, Lemma \ref{lem:incenterradius} shows that $r_p \leq \diam(\triangle_p) \leq \diam(p)$. Moreover, the non-collinearity condition \eqref{eq:non-collineardef} and the inradius formula \eqref{eq:inradiusformula} show that 
\begin{eqnarray*}
\varepsilon \diam(p)^2 &\leq& \area(\triangle_p)=
\frac{r_p(|a_p-b_p|+|a_p-c_p|+|b_p-c_p|)}{2} \\ 
&\leq& \frac{3 r_p \diam(\triangle_p)}{2} \leq \frac{3 r_p \diam(p)}{2}. 
\end{eqnarray*} 
The third claim follows. The first claim follows from the second and third claims and the triangle inequality.  
\end{proof}

\begin{proof}[Proof of Proposition \ref{prop:energy}] 
We first use \eqref{eq:rdef}, \eqref{eq:edef}, and the packing condition \eqref{eq:packingcondition} to bound the cardinality of $\mathcal{E}_k$. Fix $1 \leq k \leq L$. By \eqref{eq:rdef} and \eqref{eq:edef}, every $p \in \mathcal{E}_k$ satisfies   
\begin{equation}\label{eq:boo1}
	\varepsilon^{10} r_k \leq r_p \leq r_k \quad 
    \text{and} \quad |o_p-o_k| \leq \frac{r_k}{\varepsilon^{3}}. 
\end{equation} 
Moreover, by \eqref{eq:boo1} and Lemma \ref{lem:voivoi}, 
\begin{equation} \label{eq:edista} 
\dist(p,o_k) \leq |o_p-o_k|+\dist(o_p,p) \leq \frac{r_k}{\varepsilon^{3}} +\frac{3r_k}{2\varepsilon} \leq \frac{2r_k}{\varepsilon^{3}} 
\end{equation} 
for every $p \in \mathcal{E}_k$; here we also used the assumption $\varepsilon<\nicefrac{1}{100}$. 

By \eqref{eq:edista}, every $p \in \mathcal{E}_k$ intersects $D(k):=\mathbb{D}(o_k,2r_k/\varepsilon^3)$. Note also that 
\begin{equation} \label{eq:pele} 
\varepsilon^{10} r_k \leq r_p \leq \diam(p) \leq \frac{3r_k}{2\varepsilon}  \quad \text{for every } p \in \mathcal E_k
\end{equation} 
by \eqref{eq:edef} and \eqref{eq:boo1}. 
We may cover $D(k)$ with at most $\varepsilon^{-28}$ disks $D$ of radius $\varepsilon^{10}r_k$. Thus, combining \eqref{eq:pele} with the packing condition \eqref{eq:packingcondition} applied to each of the disks $D$ yields 
\begin{equation}\label{eq:boo2}
	\card \mathcal{E}_k \leq \varepsilon^{-28}N \quad \text{for every } 1 \leq k \leq L.  
\end{equation} 
We now estimate the energy of $\rho$. Since 
\begin{equation} \label{eq:crix}
\int_{G} \rho^2 \, dA \leq \area(\mathbb{D}(0,4)) = 16\pi, 
\end{equation} 
it suffices to estimate the sum of $\rho^2$ over $\mathcal{C}_L$ and $\mathcal{E}$. An argument similar to \eqref{eq:boo2} yields 
\begin{equation} \label{eq:crix2}
\sum_{p \in \mathcal{C}_L} \rho(p)^2 =	\card \mathcal{C}_L \leq \varepsilon^{-21}N. 
\end{equation}
Moreover, combining the definition of $\rho$ with \eqref{eq:pele} and \eqref{eq:boo2} yields 
\begin{equation} \label{eq:trap}
\sum_{p \in \mathcal{E}_k} \rho(p)^2 
	= \varepsilon^{-30}\sum_{p \in \mathcal{E}_k} \diam(p)^2 
	\leq 3\varepsilon^{-32}\,r^2_k\cdot \card \mathcal{E}_k \leq 
    3\varepsilon^{-60}N r^2_k 
\end{equation} 
for every $1 \leq k \leq L$. The disks 
$\mathbb{D}(o_k,\varepsilon^{6} r_k/2)$, $1 \leq k \leq L$, are pairwise disjoint by Lemma \ref{lem:inter1}. Since $p_k \in \mathcal{C}_S$, the definition of $\mathcal C_S$
and Lemma \ref{lem:voivoi} give
$$
r_k \leq \diam(p_k) <\varepsilon^{10} \quad \text{and} \quad 
\dist(o_k,p_k) \leq \diam(p_k) < \varepsilon^{10}. 
$$
Combining with the assumptions $p_k \cap \mathbb{D}(0,4) \neq \emptyset$ and $\varepsilon<\nicefrac{1}{100}$, we have 
$$ 
|o_k| \leq \dist(p_k,0)+\diam(p_k) + \dist(o_k,p_k) < 5, \quad 
 1 \leq k \leq L. 
$$

We conclude that the disks $\mathbb{D}(o_k,\varepsilon^6 r_k/2)$, $1 \leq k \leq L$, are pairwise disjoint subsets of $\D(0,6)$. In particular, we have 
\begin{equation} \label{eq:trap2}
	\frac{\pi\varepsilon^{12}}{4}\sum_{k=1}^L r_k^2=\sum_{k=1}^L \area(\mathbb{D}(o_k,\varepsilon^6 r_k/2)) \leq \area(\mathbb{D}(0,6))= 36\pi.  
\end{equation} 
Combining \eqref{eq:trap} and \eqref{eq:trap2} yields 
\begin{equation} \label{eq:crix3}
\sum_{p \in \mathcal{E}} \rho(p)^2=\sum_{k=1}^{L} \sum_{p \in \mathcal{E}_k} \rho(p)^2 \leq \varepsilon^{-74}N.
\end{equation} 
Since $\varepsilon<\nicefrac{1}{100}$, combining \eqref{eq:crix}, \eqref{eq:crix2}, and \eqref{eq:crix3} gives the desired estimate 
$$
\int_{G}\rho^2 \, dA + \sum_{p \in \mathcal{C}(G)} \rho(p)^2 
\leq 16\pi+\varepsilon^{-21}N+\varepsilon^{-74}N 
\leq \varepsilon^{-75}N. 
$$
\end{proof}

\subsection{Admissibility of \texorpdfstring{$\rho$}{rho} for 
\texorpdfstring{$\Gamma$}{Gamma}: grouping buried elements together}

Theorem \ref{thm:keyestimate} follows by Proposition \ref{prop:energy} and the admissibility of $\rho$ for $\Gamma$. 

\begin{proposition} \label{prop:admissibility}
The function $\rho$ is admissible for $\Gamma$. 
\end{proposition}

The rough idea of the proof of Proposition \ref{prop:admissibility} is to choose suitable subpaths of a given $\gamma \in \Gamma$, which completely avoid $\mathcal{B}$, so that we can apply Proposition \ref{prop:blockingprop} and the weights $\rho$ on the elements $p \in \mathcal{E}$ to compensate for the ``loss'' induced by the weight $\rho$ being zero on $\mathcal{B}$. 

The remainder of Section \ref{sec:key} is devoted to the proof of Proposition \ref{prop:admissibility}. The proof is divided into several subsections. We fix $\gamma:[a,b] \to \hat G$ in $\Gamma$, i.e.,   
\begin{equation} \label{eq:turki}
\gamma(a) \cap \mathbb{S}(0,4)\neq \emptyset, \, 
\gamma(b) \cap \mathbb{S}(0,1) \neq \emptyset, \text{ and } \gamma(t) \subset \mathbb A(0,1), \, a <t<b. 
\end{equation}
Recall the notation $\mathcal{C}(\gamma):=\mathcal{C}(G) \cap |\gamma|$ and 
$$ 
\hat{\ell}(\rho,\gamma):=\int_\gamma \rho \, ds 
+ \sum_{p \in \mathcal{C}(\gamma)} \rho(p). 
$$ 
In what follows, we also denote 
$$
\mathcal{E}(\gamma):=\mathcal{E} \cap |\gamma| \quad \text{and} 
\quad \mathcal{B}(\gamma):= \mathcal{B} \cap |\gamma|,  
$$ 
and define $\mathcal{E}_k(\gamma),\mathcal{B}_k(\gamma)$ similarly. 

Our goal is to prove $\hat{\ell}(\rho,\gamma) \geq 1$. Since $\rho=1$ in $\mathcal{C}_L$, we may assume that  
$$
\mathcal C(\gamma) \cap \mathcal{C}_L = \emptyset. 
$$ 
Moreover, since $\rho=1$ in  $G \cap \mathbb A(0,1)$ and $\rho(p)=\diam(p)/\varepsilon^{15}$ for $p \in \mathcal E$, the triangle inequality shows that $\hat \ell(\rho,\gamma) \geq 1$ if 
$\mathcal B(\gamma) =\emptyset$ (see Lemma \ref{lem:juures}). Therefore, we lose no generality by 
assuming that 
$$ 
\mathcal B(\gamma) \neq \emptyset. 
$$

We may also assume that for every $p \in \mathcal{C}(\gamma)$ there is a unique $t_p$ such that $\gamma(t_p)=p$: if there are $t_p<t'_p$ such that $\gamma(t_p)=p=\gamma(t'_p)$, then we can concatenate the paths $\gamma|_{[a,t_p]}$ and $\gamma|_{[t'_p,b]}$ to obtain a path $\tilde{\gamma} \in \Gamma$ such that 
$$
\hat{\ell}(\rho,\gamma) \geq \hat{\ell}(\rho,\tilde{\gamma}). 
$$






The following auxiliary results will be used to group together elements of $\mathcal{B}_k$ for different values of $k$. 

\begin{lemma} \label{lem:fishbone}
Let $1 \leq k <m \leq L$, and suppose that 
\begin{equation} \label{eq:ole}
\text{there is } p \in \mathcal{B}_k(\gamma) \text{ such that } 
\dist(o_m,p) \leq 2 \varepsilon^{-2}r_m. 
\end{equation} 
Then 
$$ 
|o_k-o_m| \leq 2\varepsilon^6 r_k \quad \text{and} \quad  
r_m \leq \varepsilon^{10}r_k. 
$$
\end{lemma}
 
\begin{proof} 
First, we have 
\begin{eqnarray*}
|o_k-o_m| &\leq& \dist(o_k,p)+\dist(o_m,p)+\diam(p) 
\leq \dist(o_p,p)+|o_p-o_k| \\ 
&+&\dist(o_m,p)+\diam(p) 
\leq 2 \diam(p)+|o_p-o_k|+\dist(o_m,p) \\ 
&\leq& 4\varepsilon^{-1}r_p+2\varepsilon^{-2}r_m+|o_p-o_k|
\end{eqnarray*} 
by the triangle inequality, Lemma \ref{lem:voivoi}, and our assumption. Moreover, since $p \in \mathcal{B}_k$, we have $r_p<\varepsilon^{10}r_k$ and $|o_p-o_k| \leq \varepsilon^6 r_k$ by \eqref{eq:bdef}. 
Since also $r_m \leq r_k$ (because $m >k$), we conclude that 
\begin{equation} \label{eq:kanski} 
|o_k-o_m| \leq (4\varepsilon^{9}+2\varepsilon^{-2}+\varepsilon^6)r_k \leq \varepsilon^{-3}r_k.  
\end{equation} 
Since $p_m\in\mathcal E_m$ and $m>k$, we have
$p_m\notin\mathcal E_k$. Thus, \eqref{eq:kanski} and \eqref{eq:edef}
imply that $r_m<\varepsilon^{10}r_k$. In particular, we can improve \eqref{eq:kanski} to 
$$ 
|o_k-o_m| \leq 4\varepsilon^{-1}r_p+2\varepsilon^{-2}r_m+|o_p-o_k| \leq (4\varepsilon^{9}+2\varepsilon^{8}+\varepsilon^6)r_k \leq 2\varepsilon^6 r_k. 
$$
\end{proof}

We are ready to group elements of the sets $\mathcal{B}_k(\gamma)$. We denote 
$$
J_\gamma:=\{1 \leq k \leq L: \, \mathcal B_k(\gamma) \neq \emptyset\}. 
$$

\begin{definition}  \label{def:roots}
\begin{enumerate} 
\item We call a finite sequence 
\begin{equation}\label{eq:bundle}
\mathbf K:=\{k_1,k_2,\ldots,k_l\} \subset J_\gamma, \quad 1 \leq k_1<k_2 < \cdots < k_l \leq L, 
\end{equation}
a \emph{$\gamma$-bundle} if $k=k_j$ and $m=k_{j+1}$, $1\leq j \leq l-1$, satisfy \eqref{eq:ole}.  
\item We call $m \in J_\gamma$ a \emph{$\gamma$-descendant} of $k \in J_\gamma$, $1 \leq k \leq m$, if there is a $\gamma$-bundle $\mathbf K$ as in \eqref{eq:bundle} such that $k=k_1$ and $m=k_l$. 
\item An index $K \in J_\gamma$ is a \emph{$\gamma$-root} if $K$ is not a $\gamma$-descendant of any $k \in J_\gamma$, $1 \leq k<K$, and a \emph{$\gamma$-root of $m$} if $K$ is a $\gamma$-root and $m$ is a $\gamma$-descendant of $K$. 
We denote the collection of $\gamma$-roots 
$$
\mathbf{R}:=\{K_1,K_2,\ldots,K_\ell\}, \quad 1 \leq K_1 < \cdots < K_\ell \leq L. 
$$
Given $K_\mu \in \mathbf{R}$, we denote $O_\mu:=o_{K_\mu}$, $R_\mu:=r_{K_\mu}$, and 
$$
\mathcal{D}_\mu :=\cup\{\mathcal{B}_m(\gamma): \, K_\mu \leq m \leq L, \,  K_\mu \text{ is a }\gamma \text{-root of } m \, \}. 
$$
\end{enumerate}
\end{definition}

The following properties of $\mathcal D_\mu$ will be important. First, we show that every $p \in \mathcal D_\mu$ lies near $O_\mu$. 

\begin{lemma} \label{lem:bundle} 
Every $m \in J_\gamma$ has exactly one $\gamma$-root, and thus the collection $\{\mathcal{D}_\mu \}$ is a partition of $\mathcal{B}(\gamma)$. All $1 \leq \mu \leq \ell$ satisfy 
\begin{equation} \label{eq:kevil}
\bigcup_{p \in \mathcal{D}_\mu} p \subset 
\D(O_\mu,4\varepsilon^6R_{\mu}). 
\end{equation}
\end{lemma}

\begin{proof} 
We first prove \eqref{eq:kevil}. Let $1 \leq \mu \leq \ell$, 
$\{K_\mu=k_1,k_2,\ldots, k_l\}$ a $\gamma$-bundle, and $p \in \mathcal{B}_{k_l}$. If $l=1$, then \eqref{eq:kevil} follows from \eqref{eq:bdef} and Lemma \ref{lem:voivoi}. If $l \geq 2$, we have 
\begin{eqnarray} \label{eq:kiksu} 
\begin{aligned}
& |o_p-O_\mu| \leq |o_p-o_{k_l}| + \sum_{n=1}^{l-1}|o_{k_{n+1}}-o_{k_{n}}| 
\leq \varepsilon^6 r_{k_l} +2\varepsilon^6 \sum_{n=1}^{l-1} r_{k_n} \\
& \leq \varepsilon^{6+10(l-1)} R_{\mu} + 2\varepsilon^6 R_\mu \sum_{n=1}^{l-1} \varepsilon^{10(n-1)} \leq 3 \varepsilon^6 R_\mu   
\end{aligned}
\end{eqnarray} 
by the definition of $\mathcal{B}_{k_l}$ (see \eqref{eq:bdef}) and Lemma \ref{lem:fishbone}. 
On the other hand,  
\begin{equation} \label{eq:kiksu2} 
\diam(p)+\dist(o_p,p) \leq  \frac{3r_p}{\varepsilon} \leq \varepsilon^6 R_\mu 
\end{equation}
by Lemmas \ref{lem:voivoi} and \ref{lem:fishbone}.  Combining \eqref{eq:kiksu} and \eqref{eq:kiksu2} with the triangle 
inequality yields \eqref{eq:kevil}. 

To start the proof of the first claim, note that every $n \in J_\gamma$ has at least one $\gamma$-root (see Definition \ref{def:roots}). Suppose towards a contradiction that $k$ and $m$, $k < m <n$, are $\gamma$-roots of $n$. Then there are $\gamma$-bundles 
$$
\mathbf K=\{k=k_1, \ldots,k_\sigma= n\} \quad \text{and} \quad \{m=m_1, \ldots,m_\nu= n\}. 
$$ 
We have $k_j \neq m$ for all $j=2,\ldots,\sigma$, since otherwise $m$ is not a root. Let $k \leq k_{j(0)} <n$ be the largest element of $\mathbf K$ such that $k_{j(0)} < m$. By the definition of a $\gamma$-bundle, there is  
$p \in \mathcal{B}_{k_{j(0)}}(\gamma)$ such that 
$$
\dist(p,o_{k_{j(0)+1}}) \leq 2\varepsilon^{-2}r_{k_{j(0)+1}}. 
$$
Applying Lemma \ref{lem:fishbone} and estimating as in \eqref{eq:kiksu} yields (note that $r_{k_{j(0)+1}} \leq r_m$ since $m < k_{j(0)+1}$)
$$
|o_m-o_{k_{j(0)+1}}| \leq \sum_{j=j(0)+1}^{\sigma-1} |o_{k_j}-o_{k_{j+1}}|
+\sum_{j=1}^{\nu-1} |o_{m_j}-o_{m_{j+1}}| \leq 6\varepsilon^6 r_m. 
$$
Thus, by \eqref{eq:edef} and since $p_{k_{j(0)+1}} \notin \mathcal{E}_m$, we have $r_{k_{j(0)+1}} <\varepsilon^{10} r_m$. Combining with the above estimates yields 
\begin{eqnarray*} \begin{aligned}
& \dist(p,o_m) \leq \dist(p,o_{k_{j(0)+1}})+ |o_m-o_{k_{j(0)+1}}|\\ 
& \leq 
2\varepsilon^{-2}r_{k_{j(0)+1}} + 6\varepsilon^6 r_m \leq 2\varepsilon^{-2} r_m.   
\end{aligned}
\end{eqnarray*} 
We conclude using \eqref{eq:ole} that $\{k,\ldots,k_{j(0)},m\}$ is a $\gamma$-bundle and, therefore, $m$ is not a $\gamma$-root. This is a contradiction. The proof is complete. 
\end{proof} 

We establish a criterion similar to \eqref{eq:ole}, for $\mathcal D_\mu$ and $O_\nu$, $1 \leq \mu < \nu \leq \ell$. 

\begin{lemma} \label{lem:fishboned}
If $1 \leq \mu < \nu \leq \ell$, then 
$$
\dist(p,O_\nu) > \varepsilon^{-2}R_\nu \quad \text{for every } p \in \mathcal D_\mu. 
$$
\end{lemma}

\begin{proof} 
Suppose towards a contradiction that there is $p \in \mathcal B_m(\gamma) \subset \mathcal D_\mu$ such that 
\begin{equation} \label{eq:kaaleppi} 
\dist(p,O_\nu) \leq \varepsilon^{-2}R_\nu. 
\end{equation}
If $m < K_\nu$, then \eqref{eq:kaaleppi} and \eqref{eq:ole} show that 
$K_\nu$ is a $\gamma$-descendant of $m$. This is a contradiction, since $K_\nu$ is a $\gamma$-root. The case $m=K_\nu$ is not possible, since otherwise $K_\nu$ would be a $\gamma$-descendant of $K_\mu$, and therefore not a root. 

Finally, suppose that $K_\nu <m$. Since $m$ is a $\gamma$-descendant of $K_\mu$, there is a $\gamma$-bundle 
$\mathbf K=\{K_\mu=k_1,\ldots,k_\sigma=m\}$. Let $K_\mu \leq k_{j(0)}<m$ be the largest element of $\mathbf K$ such that $k_{j(0)}<K_\nu$. We have 
\begin{eqnarray} \label{eq:cheesecake} \begin{aligned} 
|o_{k_{j(0)+1}}-O_\nu|&\leq |o_{k_{j(0)+1}}-o_m|+|o_m-O_\nu| \\ 
& \leq  \sum_{j=j(0)+1}^{\sigma-1} |o_{k_j}-o_{k_{j+1}}| + \dist(o_m,p) \\ 
&+\diam(p)+\dist(p,O_\nu) 
\leq 5 \varepsilon^{6}r_{k_{j(0)+1}}+ \varepsilon^{-2}R_\nu <\varepsilon^{-3}R_\nu  
\end{aligned}
\end{eqnarray} 
by Lemmas \ref{lem:voivoi} and \ref{lem:fishbone}, and \eqref{eq:kaaleppi}. 
Since $K_\nu \leq k_{j(0)+1}$, the definition \eqref{eq:edef} of $\mathcal E_{K_\nu}$ 
and \eqref{eq:cheesecake} show that $r_{k_{j(0)+1}}< \varepsilon^{10}R_\nu$. 

Since $k_{j(0)+1}$ is a $\gamma$-descendant of $k_{j(0)}$, there is $q \in \mathcal B_{k_{j(0)}}(\gamma)$ such that 
\begin{equation} \label{eq:neverending} 
\dist(q,o_{k_{j(0)+1}}) \leq 2 \varepsilon^{-2}r_{k_{j(0)+1}}. 
\end{equation}
Combining \eqref{eq:cheesecake} and \eqref{eq:neverending} gives 
$$ 
\dist(q,O_\nu) \leq 3 \varepsilon^{-2}r_{k_{j(0)+1}}+\varepsilon^{-2}R_\nu \leq 2 \varepsilon^{-2} R_\nu. 
$$
We conclude that $K_\nu$ is a $\gamma$-descendant of $k_{j(0)}$ (recall \eqref{eq:ole}). This is a contradiction. The proof is complete. 

\end{proof}

\subsection{Admissibility of \texorpdfstring{$\rho$}{rho} for 
\texorpdfstring{$\Gamma$}{Gamma}: Modified endpoints}

We prepare for the choices of subpaths of $\gamma$ that avoid $\mathcal B$. The following lemma will be useful for choosing their endpoints. Recall the $\gamma$-roots $K_\mu \in \mathbf R$ and sets $\mathcal D_\mu$ in Definition \ref{def:roots}. 

\begin{lemma} \label{lem:jruusing}
Given $x \in \hat G$, let $1  \leq \mu < \nu \leq \ell$ be 
such that 
\begin{equation} \label{eq:kaivanto0}
\dist(x,O_\mu) > \varepsilon^{-3} R_\mu/2 \quad \text{and} \quad 
\dist(x,O_\nu) \leq \varepsilon^{-3} R_\nu/2. 
\end{equation}
Then 
\begin{equation} \label{eq:kaivanto}
\dist(p,O_\mu) > \varepsilon^{-2}R_\mu \quad \text{for all } p \in \mathcal{D}_\nu. 
\end{equation}
\end{lemma}

\begin{proof} 
Suppose towards a contradiction that there is $p \in \mathcal D_\nu$ such that 
\begin{equation} \label{eq:kaivanto2} 
\dist(p,O_\mu) \leq \varepsilon^{-2}R_\mu. 
\end{equation}
Now \eqref{eq:kaivanto0}, Lemma \ref{lem:bundle}, and \eqref{eq:kaivanto2} yield 
\begin{eqnarray} \label{eq:krontah} \begin{aligned}
\dist(x,O_\mu) &\leq \dist(x,O_\nu)+\dist(O_\nu,p)+ \diam(p) \\ &+\dist(p,O_\mu) \leq \varepsilon^{-3}R_\nu/2+4\varepsilon^{6}R_\nu
+\varepsilon^{-2}R_\mu.  
\end{aligned}
\end{eqnarray}

Moreover, since $R_\nu \leq R_\mu$, Lemma \ref{lem:bundle} and \eqref{eq:kaivanto2} yield 
\begin{eqnarray*} 
|O_\nu-O_\mu| &\leq& \dist(p,O_\mu)+\dist(p,O_\nu)+
\diam(p) \\ 
&\leq& \varepsilon^{-2}R_\mu+4 \varepsilon^6 R_\nu \leq \varepsilon^{-3}R_\mu, 
\end{eqnarray*}
and thus $R_\nu< \varepsilon^{10} R_\mu$ by \eqref{eq:edef}. Now \eqref{eq:krontah} gives 
$$ 
\dist(x,O_\mu) \leq \varepsilon^{-3}R_\mu/2. 
$$
This contradicts \eqref{eq:kaivanto0}. The proof is complete. 
\end{proof}

We now start to choose the subpaths of $\gamma$ by first modifying the endpoints $a$ and $b$. Let $\mathbf{A}$ be the collection of those $1 \leq \mu \leq \ell$ for which 
$$
\dist(\gamma(a),O_\mu) \leq \varepsilon^{-3}R_\mu/2. 
$$
We define $a \leq \tilde a < b$ as follows: If $\mathbf{A}=\emptyset$, then we set $\tilde a:= a$. Otherwise, let $\mu(a)$ be the smallest element of $\mathbf{A}$ and set 
$$
\tilde a := \max\{a \leq t < b: \, \gamma(t) \in \mathcal D_{\mu(a)}\}. 
$$

Similarly, let $\mathbf{B}$ be the collection of those $1 \leq \mu \leq \ell$ for which $\gamma(t) \in \mathcal D_\mu$ for some $\tilde a < t <b$ and 
$$
\dist(\gamma(b),O_\mu) \leq \varepsilon^{-3}R_\mu/2. 
$$
We define $\tilde a < \tilde b \leq b$ as follows: If $\mathbf{B}=\emptyset$, then we set $\tilde b:= b$. Otherwise, let $\mu(b)$ be the smallest element of $\mathbf{B}$ and set 
$$ 
\tilde b := \min\{\tilde a < t \leq b: \, \gamma(t) \in \mathcal D_{\mu(b)}\}. 
$$

\begin{lemma} \label{lem:scamming} The following properties are true: 
\begin{enumerate} 
\item $\dist(\gamma(\tilde a),\gamma(\tilde b)) \geq 1$. 
\item The sets $(\tilde a,\tilde b) \cap \gamma^{-1}(\mathcal D_\mu)$  
form a partition of $(\tilde a,\tilde b) \cap \gamma^{-1}(\mathcal B)$. 
\item If $1 \leq \mu \leq \ell$ is such that $(\tilde a, \tilde b) \cap \gamma^{-1}(\mathcal D_\mu) \neq \emptyset$, then  
\begin{equation}  \label{eq:isi}
\dist(\gamma(\tilde a) \cup \gamma(\tilde b),O_\mu) >\varepsilon^{-2}R_\mu.  
\end{equation}
\end{enumerate} 
\end{lemma}

\begin{proof} 

By \eqref{eq:turki} and the assumption $|\gamma|\cap\mathcal C_L=\emptyset$, we have
\begin{equation}\label{eq:|z||w|}
    |z|\ge 4-\varepsilon^{10}
\quad\text{for all } z\in\gamma(a),
\qquad
|w|\le 1+\varepsilon^{10}
\quad\text{for all } w\in\gamma(b).
\end{equation}
Combining the definitions of $\widetilde a$ and $\widetilde b$ with \eqref{eq:kevil} gives
$$\dist(z,\gamma(a))\le\dist (O_{\mu(a)},\gamma(a))+|z-O_{\mu(a)}|<\varepsilon^7\quad\text{for all }z\in\gamma(\widetilde a),$$
and 
$$\dist(w,\gamma(b))<\varepsilon^7\quad\text{for all } w\in\gamma(\widetilde b).$$
Consequently, by \eqref{eq:|z||w|} and the above estimates, we conclude that
$$\begin{aligned}
\dist\bigl(\gamma(\widetilde a),\gamma(\widetilde b)\bigr)
&\ge \dist\bigl(\gamma(a),\gamma(b)\bigr)-2\varepsilon^7\\
&\ge 3-2\varepsilon^{10}-2\varepsilon^7
>1.
\end{aligned}$$
This proves the first claim. The second claim follows
from Lemma \ref{lem:bundle}.

To prove the last claim, we show only that 
\begin{equation}  \label{eq:isi2}
\dist(\gamma(\tilde a),O_\mu) >\varepsilon^{-2}R_\mu; 
\end{equation} 
the corresponding estimate for $\gamma(\tilde b)$ is proved similarly. 
We may assume that $\mathbf A \neq \emptyset$, since otherwise \eqref{eq:isi2} is trivial. If $\mu > \mu(a)$, then \eqref{eq:isi2} holds by Lemma \ref{lem:fishboned}. 

If $\mu=\mu(a)$ then $(\tilde a, \tilde b) \cap \gamma^{-1}(\mathcal D_\mu) = \emptyset$ by the definition of $\tilde a$, so the last claim is trivial. If $\mu< \mu(a)$, then the assumptions of Lemma \ref{lem:jruusing} hold for $\mu$, $\nu:=\mu(a)$, and $x:=\gamma(a)$, since $\mu(a)$ is the smallest element of $\mathbf A$. We therefore conclude from the lemma that 
$$ 
\dist(p,O_\mu) > \varepsilon^{-2}R_\mu \quad \text{for all } p \in \mathcal D_{\mu(a)}. 
$$
Since $\gamma(\tilde a) \in \mathcal D_{\mu(a)}$ by the definition of $\tilde a$, \eqref{eq:isi2} follows. 
\end{proof}


\subsection{Admissibility of \texorpdfstring{$\rho$}{rho} for 
\texorpdfstring{$\Gamma$}{Gamma}: choosing subpaths}

We continue to prove the admissibility of $\rho$ for $\Gamma$ (Proposition \ref{prop:admissibility}). We next construct subpaths of $\gamma$. They will avoid $\mathcal{B}$, allowing us to apply the following consequence of the triangle inequality; recall that $\rho=1$ on $(\mathbb{D}(0,4)\setminus \overline \D(0,1)) \cap G$  and the notation $\mathcal{E}(\alpha)=\{p \in \mathcal E \, \cap |\alpha|\}$. 

\begin{lemma} \label{lem:juures} 
Let $\tilde a \leq c < d \leq \tilde b$ be such that $\gamma(t) \notin \mathcal{B}$ 
for every $c <t <d$. Then 
$$
\dist(\gamma(c),\gamma(d)) \leq \int_{\gamma|_{(c,d)}} \rho \, ds 
+  \sum_{p \in \mathcal{E}(\gamma|_{(c,d)})} \diam(p). 
$$
\end{lemma} 

We are ready to construct the desired subpaths of $\gamma$. 

\begin{proposition} \label{prop:kukinto} For every $0 \leq \mu \leq \ell$ there is a collection $\mathcal I_\mu$ of open and pairwise disjoint segments $I \subset (\tilde a,\tilde b)$ such that $\mathcal I_0=\{(\tilde a,\tilde b)\}$ and the following properties hold: 

\begin{enumerate}
\item \label{it:kukinto1} If $I \in \mathcal I_\mu$, $1 \leq \mu \leq \ell$, then $I \subset (\tilde a, \tilde b) \setminus \cup_{\nu=1}^\mu \gamma^{-1}(\mathcal D_\nu)$. 

\item \label{it:kukinto2} If $I=(c,d) \in \mathcal I_{\mu}$, $0 \leq \mu \leq \ell-1$, and $I \cap \gamma^{-1}(\mathcal D_{\nu}) \neq \emptyset$ for some $\mu+1 \leq \nu \leq \ell$, then 
\begin{equation} \label{eq:rain}  
\dist(\gamma(c) \cup \gamma(d),O_{\nu}) > \varepsilon^{-2}R_\nu.    
\end{equation}

\item \label{it:kukinto3} If $I=(c,d) \in \mathcal I_{\mu-1}$, $1 \leq \mu \leq \ell$, and if 
$I_j=(c_j,d_j)$, $j=1,\ldots,N_I$, are the elements of $\mathcal I_{\mu}$ for which $I_j \cap I \neq \emptyset$, then 
\begin{equation} \label{eq:kraha}
I_j \subset I \quad \text{for every } j=1,\ldots,N_I,   
\end{equation}
and 
\begin{eqnarray} \label{eq:uio} \begin{aligned} 
& \dist(\gamma(c),\gamma(d)) \leq \sum_{j=1}^{N_I} \dist(\gamma(c_j),\gamma(d_j))
+\psi(\mu,I), 
\end{aligned}
\end{eqnarray}
where 
\begin{eqnarray*}
\psi(\mu,I)= \sum_{p \in \mathcal{E} \cap \gamma(I \setminus \cup_{j=1}^{N_I} I_j)} \frac{\diam(p)}{\varepsilon^{14}}. 
\end{eqnarray*} 

\item \label{it:kukinto4} If $1 \leq \mu \leq \ell$, then 
\begin{equation} \label{eq:uio2} 
1 \leq \sum_{(c,d) \in \mathcal I_\mu} \dist(\gamma(c),\gamma(d)) 
+ \varphi(\mu), 
\end{equation} 
where 
\begin{eqnarray*}
\varphi(\mu)= \sum_{p \in \mathcal{E} \cap \gamma((\tilde a,\tilde b) \setminus \cup \{I \in \mathcal I_\mu\})} \frac{\diam(p)}{\varepsilon^{14}}. 
\end{eqnarray*} 
\end{enumerate}
\end{proposition}

Before proving Proposition \ref{prop:kukinto}, we use it to establish Proposition \ref{prop:admissibility}.

\begin{proof}[Proof of Proposition \ref{prop:admissibility}]
We need to show that $\hat \ell(\rho,\gamma)\geq 1$. 
By Lemma \ref{lem:scamming} and parts (\ref{it:kukinto1}) and 
(\ref{it:kukinto4}) of Proposition \ref{prop:kukinto} (applied with $\mu=\ell$),  
we have 
\begin{eqnarray} \label{eq:cleaners0}
& &I \cap \gamma^{-1}(\mathcal B) = \emptyset \quad \text{for every } I \in \mathcal I_\ell \quad \text{and} \\
\label{eq:cleaners}
& & 1 \leq \sum_{(c,d) \in \mathcal I_\ell} \dist(\gamma(c),\gamma(d)) + 
\sum_{p \in \mathcal E(\gamma)} \frac{\diam(p)}{\varepsilon^{14}}.  
\end{eqnarray} 
By \eqref{eq:cleaners0}, we may apply Lemma \ref{lem:juures} to the intervals of $\mathcal I_\ell$. Combining with \eqref{eq:cleaners} and the definition of $\rho$ (see \eqref{eq:rhodef1}, \eqref{eq:rhodef2} and \eqref{eq:rhodef3}) gives 
$$ 
1 \leq \int_\gamma \rho \, ds+ 
\sum_{p \in \mathcal E(\gamma)}\diam(p)\Big(1+ \varepsilon^{-14}\Big) 
\leq \hat \ell(\rho,\gamma), 
$$ 
as desired. 
\end{proof}

Theorem \ref{thm:keyestimate} follows by combining Propositions \ref{prop:energy} and \ref{prop:admissibility}. Therefore, the proof of Theorem \ref{thm:keyestimate} is complete after we have proved Proposition \ref{prop:kukinto}. 


\section{Proof of Proposition \ref{prop:kukinto}} 

We construct $\mathcal I_\mu$ iteratively. First, the only non-trivial claim in conditions \ref{it:kukinto1}-\ref{it:kukinto4} that concerns  $\mathcal I_0=\{(\tilde a,\tilde b)\}$ is \eqref{eq:rain}, which follows by \eqref{eq:isi}. Suppose next that $\mathcal I_{\mu-1}$, $1 \leq \mu \leq \ell$, are defined and satisfy conditions (\ref{it:kukinto1})-(\ref{it:kukinto4}). 

We define $\mathcal I_{\mu}$. Let $I \in \mathcal I_{\mu-1}$. If $I \cap \gamma^{-1}(\mathcal D_\mu)=\emptyset$, then we include $I$ in $\mathcal I_\mu$. Assume then that $I=(c,d) \in \mathcal{I}_{\mu-1}$ is such that $I\cap \gamma^{-1}(\mathcal D_\mu) \neq \emptyset$. Let 
\begin{eqnarray} \label{eq:ironmaiden} \begin{aligned}
b_\mu=b_\mu(I) :=& \min\{c < t < d: \, \gamma(t) \in \mathcal D_\mu\}  \quad \text{and} \\   
a_\mu=a_\mu(I) :=& \max\{c < t < d: \, \gamma(t) \in \mathcal D_\mu\}.  
\end{aligned}
\end{eqnarray}
The construction of the intervals $I' \in \mathcal I_\mu$ for which $I \cap I' \neq \emptyset$ is divided into cases. Let $T_\mu$ be the collection of those $c<t<d$ for which $\gamma(t)\in \mathcal E_k$ for some 
$1 \leq k \leq K_\mu$ and 
\begin{equation} \label{eq:tarve}
\dist(\gamma(t),O_\mu) \leq \varepsilon^{-2}R_\mu. 
\end{equation}

\medskip
\noindent
\textbf{Case 1: $T_\mu \neq \emptyset$. } We include the subintervals $I_1,I_2$ of $I$ in $\mathcal I_\mu$ as follows: Let $\beta_\mu$ and $\alpha_\mu$, $c<\beta_\mu \leq \alpha_\mu<d$, be the smallest, respectively, largest, elements of $T_\mu$. 
Moreover, let $\boldsymbol \beta$ be the collection of those $\mu+1 \leq \nu \leq \ell$ for which there are $c<t <\beta_\mu$ such that 
\begin{equation} \label{eq:tilka} 
\gamma(t) \in \mathcal D_\nu \quad \text{and} \quad \dist(O_\nu,\gamma(\beta_\mu)) \leq \varepsilon^{-3}R_\nu/2. 
\end{equation}
If $\boldsymbol{\beta}=\emptyset$, we define 
\begin{equation} \label{eq:illi}
I_1:=(c,\min\{\beta_\mu,b_\mu\}). 
\end{equation}
Otherwise, let $\nu(\beta)$ be the smallest element of $\boldsymbol{\beta}$ and define 
\begin{equation} \label{eq:illi2}
I_1:=(c,\min\{\beta^-_\mu,b_\mu\}), 
\end{equation}
where 
\begin{equation} \label{eq:killi}
\beta^-_\mu:=\min\{c<t<\beta_\mu: \, \gamma(t) \in \mathcal D_{\nu(\beta)} \}. 
\end{equation}

Similarly, let $\boldsymbol \alpha$ be the collection of those $\mu+1 \leq \nu \leq \ell$ for which there are $\alpha_\mu<t <d$ such that 
$$ 
\gamma(t) \in \mathcal D_\nu \quad \text{and} \quad \dist(O_\nu,\gamma(\alpha_\mu)) \leq \varepsilon^{-3}R_\nu/2. 
$$
If $\boldsymbol{\alpha}=\emptyset$, we define 
$$
I_2:=(\max\{\alpha_\mu,a_\mu\},d). 
$$
Otherwise, let $\nu(\alpha)$ be the smallest element of $\boldsymbol{\alpha}$ and define 
$$
I_2:=(\max\{\alpha^+_\mu,a_\mu\},d), 
$$
where 
$$ 
\alpha^+_\mu:=\max\{\alpha_\mu<t<d: \, \gamma(t) \in \mathcal D_{\nu(\alpha)} \}. 
$$
We have concluded the construction in case 1, i.e., when $T_\mu \neq \emptyset$. 
\medskip

Now assume that $T_\mu=\emptyset$. We will consider cases 2 and 3 using the set $H_\mu(I)$ defined below in \eqref{eq:HmuI}. Note that $p_{K_\mu} \notin \gamma(I)$, since otherwise the failure of \eqref{eq:tarve} would contradict Lemma \ref{lem:voivoi}. Let $E_\mu$ be the set $E_{p}$ defined in \eqref{eq:eepee} for $p:=p_{K_\mu}$. Moreover, let  
\begin{eqnarray} \label{eq:pilli} \begin{aligned}
s_1&:=\max\{c < s < b_{\mu}:\, \gamma(s) \cap \partial E_{\mu} \neq \emptyset  \}  \quad \text{and } \\ 
t_1&:=\min\{a_{\mu} < t < d:\, \gamma(t) \cap \partial E_{\mu} \neq \emptyset  \}. 
\end{aligned}
\end{eqnarray} 
To guarantee that $s_1,t_1$ are well-defined, we first note that 
$\gamma(b_{\mu}) \cup \gamma(a_{\mu}) \subset E_{\mu}$ by Lemma \ref{lem:bundle}, Lemma \ref{lem:incenterradius}, and the definition \eqref{eq:eepee} of $E_{\mu}$. On the other hand,  
$$
(\gamma(c) \cup \gamma(d)) \cap \overline E_{\mu}=\emptyset  
$$
by condition (\ref{it:kukinto2}) applied to $\mu-1$ and Lemma \ref{lem:voivoi}, so $s_1$ and $t_1$ are indeed well defined.  Let
\begin{equation}\label{eq:HmuI}
    H_\mu(I):=\{b_\mu<t<a_\mu:
\gamma(t)\cap\partial E_\mu\ne\emptyset\}.
\end{equation}

\medskip
\noindent 
\textbf{Case 2: $T_\mu=\emptyset$ and $H_\mu(I) \neq \emptyset$. } 
We include the subintervals $I_1,I_2,I_3,I_4$ of $I$ in $\mathcal I_\mu$ as follows: fix $\tau\in H_\mu(I)$ and define
$$\xi:=\max\{t\in[b_\mu,\tau):\gamma(t)\in\mathcal D_\mu\},
\qquad
\eta:=\min\{t\in(\tau,a_\mu]:\gamma(t)\in\mathcal D_\mu\}.$$
Note that $\gamma(\tau)\notin\mathcal D_\mu$ by Lemma~\ref{lem:bundle}, and therefore 
$$
J\cap\gamma^{-1}(\mathcal D_\mu)=\emptyset,  \quad \text{where } J:=(\xi,\eta). 
$$

Let $\mathbf S_\tau$ be the collection of those $\mu+1 \leq \nu \leq \ell$ for which there are $\xi<t<\eta$ such that 
\begin{equation} \label{eq:rolling1} 
\gamma(t) \in \mathcal D_\nu \quad \text{and} \quad \dist(O_\nu,\gamma(\tau)) \leq \varepsilon^{-3} R_\nu/2.  
\end{equation} 
If $\mathbf S_\tau=\emptyset$, let $\sigma_-=\sigma_+=\tau$.
Otherwise, let $\nu_0(s_\tau)$ be the smallest element of $S_\tau$ and 
\begin{eqnarray}\label{eq:defofsigma} \begin{aligned}
  &  \sigma_-:=\min\{t\in J:\gamma(t)\in\mathcal D_{\nu_0(s_\tau)}\}, \\
& \sigma_+:=\max\{t\in J:\gamma(t)\in\mathcal D_{\nu_0(s_\tau)}\}. 
\end{aligned}
\end{eqnarray}
We include the following subintervals of $I=(c,d)$ in $\mathcal I_\mu$:
\begin{equation}\label{eq:gap-four-intervals}
I_1:=(c,b_\mu),\quad
I_2:=(\xi,\sigma_-),\quad
I_3:=(\sigma_+,\eta),\quad
I_4:=(a_\mu,d); 
\end{equation}
note that 
$$c<b_\mu\le\xi<\sigma_-\le\sigma_+<\eta\le a_\mu<d. $$
 We have concluded the construction in case 2.

\medskip
\noindent 
\textbf{Case 3: $T_\mu=\emptyset$ and $H_\mu(I) = \emptyset$. } We include the subintervals $I_1,I_2,I_3,I_4$ of $I$ in $\mathcal I_\mu$ as follows: let $\mathbf S$ be the collection of those $\mu+1 \leq \nu \leq \ell$ for which there are $c<t<b_\mu$ such that 
\begin{equation} \label{eq:rolling} 
\gamma(t) \in \mathcal D_\nu \quad \text{and} \quad \dist(O_\nu,\gamma(s_1)) \leq \varepsilon^{-3} R_\nu/2.  
\end{equation} 
We define 
$$ 
I_1:=(c,s_1^-) \quad \text{and} \quad I_2:=(s_1^+,b_\mu), 
$$ 
where $s_1^-,s_1^+$ are defined as follows: if $\mathbf S = \emptyset$, we set $s_1^-=s_1^+=s_1$. Otherwise, let $\nu(s)$ be the smallest element of $\mathbf S$ and define 
\begin{eqnarray} \label{eq:bowling} \begin{aligned}
s_1^-&:=\min\{c<t<b_\mu: \, \gamma(t)\in \mathcal D_{\nu(s)}\} \quad \text{and }  \\ 
s_1^+&:=\max\{c<t<b_\mu: \, \gamma(t)\in \mathcal D_{\nu(s)}\}. 
\end{aligned}
\end{eqnarray}

Secondly, let $\mathbf T$ be the collection of those $\mu+1 \leq \nu \leq \ell$ for which there are $a_\mu<t<d$ such that 
$$
\gamma(t) \in \mathcal D_\nu \quad \text{and} \quad \dist(O_\nu,\gamma(t_1)) \leq \varepsilon^{-3} R_\nu/2.  
$$ 
We define 
$$ 
I_3:=(a_\mu,t^-_1) \quad \text{and} \quad I_4:=(t^+_1,d), 
$$ 
where $t_1^-,t_1^+$ are defined as follows: if $\mathbf T = \emptyset$, we set $t_1^-=t_1^+=t_1$. Otherwise, let $\nu(t)$ be the smallest element of $\mathbf T$ and define 
\begin{eqnarray*}
t_1^-&:=&\min\{a_\mu<t<d: \, \gamma(t)\in \mathcal D_{\nu(t)}\} \quad \text{and }  \\ 
t_1^+&:=&\max\{a_\mu<t<d: \, \gamma(t)\in \mathcal D_{\nu(t)}\}. 
\end{eqnarray*}

The collection $\mathcal I_\mu$ consists of the intervals included above. We now verify conditions (\ref{it:kukinto1})-(\ref{it:kukinto4}) for $\mathcal I_\mu$.

\medskip
\noindent
\textbf{Proof of condition (\ref{it:kukinto1})}:
It follows from the construction that every $I_j\in \mathcal{I}_\mu$ is contained in its ``parent'' interval $I\in\mathcal I_{\mu-1}$ and 
$$
I_j\cap\gamma^{-1}(\mathcal D_\nu)
\subset I\cap\gamma^{-1}(\mathcal D_\nu)
=\emptyset\quad \text{for every} \,\, 1\leq\nu\leq\mu-1.$$

It remains to consider $\mathcal D_\mu$. If $I\cap\gamma^{-1}(\mathcal D_\mu)=\emptyset$, then $I_j=I$ and condition (\ref{it:kukinto1}) is satisfied. If one of cases 1 or 3 occurs, then 
$I_j \subset I\setminus[b_\mu,a_\mu]$, so condition (\ref{it:kukinto1}) is satisfied by the definitions \eqref{eq:ironmaiden} of $b_\mu$ and $a_\mu$. If case 2 occurs, we have $I_1,I_4\subset I\setminus[b_\mu,a_\mu]$ and $I_2,I_3\subset J$. Since $J\cap\gamma^{-1}(\mathcal D_\mu)=\emptyset$, we conclude that condition (\ref{it:kukinto1}) holds in all cases.

\medskip
\noindent
\textbf{Proof of condition (\ref{it:kukinto2})}: Since we assume that the condition holds for $\mu-1$, it suffices to prove the claim for the segments $I=I_j \in \mathcal I_\mu$ constructed in the three cases above. That is, we assume that $I_j \cap \gamma^{-1}(\mathcal D_\nu) \neq \emptyset$ for some $\mu+1 \leq \nu \leq \ell$ and claim that \eqref{eq:rain} is true for the endpoints of $I_j$. 

\medskip
\noindent
\textbf{Case 1: $T_\mu \neq \emptyset$. }  
Since the arguments for $I_1$ and $I_2$ are almost identical, it suffices to prove \eqref{eq:rain} for all possible endpoints of $I_1$. First, 
\begin{equation} \label{eq:uino1}
\dist(\gamma(c),O_\nu) > \varepsilon^{-2}R_\nu, 
\end{equation}
since we assume condition (\ref{it:kukinto2}) for $\mu-1$. Suppose that $I_1=(c,b_\mu)$. Since $\gamma(b_\mu) \in \mathcal D_\mu$,  Lemma \ref{lem:fishboned} shows that  
\begin{equation} \label{eq:uino2}
\dist(\gamma(b_\mu),O_\nu) > \varepsilon^{-2}R_\nu. 
\end{equation}
Combining \eqref{eq:uino1} and \eqref{eq:uino2} proves \eqref{eq:rain} in this case. 

If the set $\boldsymbol{\beta}$ defined before \eqref{eq:illi} is empty 
and $I_1=(c,\beta_\mu)$, then the definition of $\boldsymbol{\beta}$ shows that 
\begin{equation} \label{eq:uino3}
\dist(\gamma(\beta_\mu),O_\nu) > \varepsilon^{-3}R_\nu/2 \geq \varepsilon^{-2}R_\nu.  
\end{equation}
Therefore, \eqref{eq:rain} holds by \eqref{eq:uino1} and \eqref{eq:uino3}. 

Finally, suppose that $\boldsymbol{\beta} \neq \emptyset$ and $I_1=(c,\beta^-_\mu)$. Recall that $\gamma(\beta^-_\mu) \in \mathcal D_{\nu(\beta)}$, where $\nu(\beta)$ is the smallest element of $\boldsymbol{\beta}$. We claim that 
\begin{equation} \label{eq:uino4}
\dist(\gamma(\beta^-_\mu),O_\nu) > \varepsilon^{-2}R_\nu. 
\end{equation}
If $\nu(\beta) < \nu$, then \eqref{eq:uino4} follows as \eqref{eq:uino2}, with $\mu$ replaced by $\nu(\beta)$. 

The definition \eqref{eq:killi} of $\beta^-_\mu$ shows that the case $\nu(\beta)=\nu$ is not possible. Next, assume that $\nu(\beta) > \nu$. Now we can apply Lemma \ref{lem:jruusing} with $\nu$ instead of $\mu$ and $\nu(\beta)$ instead of $\nu$, to conclude that \eqref{eq:uino4} is true. Combining \eqref{eq:uino1} and \eqref{eq:uino4} gives \eqref{eq:rain}. We have proved condition (\ref{it:kukinto2}) in case 1. 

\medskip
\noindent
\textbf{Case 2: $T_\mu=\emptyset$ and $H_\mu(I) \neq \emptyset$. } 
Recall that $I_1,I_2,I_3,I_4\subset I$ are defined in \eqref{eq:gap-four-intervals}. The endpoints $c$ and $d$ satisfy \eqref{eq:uino1}. Moreover, since  
$$
\gamma(b_\mu),\gamma(a_\mu),\gamma(\xi),\gamma(\eta)\in\mathcal D_\mu, 
$$ 
\eqref{eq:uino2} holds for $b_\mu,a_\mu,\xi,\eta$ by Lemma \ref{lem:fishboned}. Therefore, it suffices to prove \eqref{eq:rain} for the endpoints $\sigma_-$ and $\sigma_+$ of $I_2$ and $I_3$, respectively. Without loss of generality, we consider the subinterval $I_2=(\xi,\sigma_-)\in \mathcal{I}_\mu$.

Suppose that $I_2 \cap \gamma^{-1}(\mathcal D_{\nu}) \neq \emptyset$ for some $\mu+1 \leq \nu \leq \ell$. If the set $\mathbf S_\tau$ is empty, then its definition \eqref{eq:rolling1} implies    
$$\dist(\gamma(\sigma_-),O_\nu)=\dist(\gamma(\tau),O_\nu)>\varepsilon^{-3}R_\nu/2
\ge\varepsilon^{-2}R_\nu, 
$$
and so \eqref{eq:rain} holds for $I_2$.  

Now suppose that $\mathbf S_\tau\neq\emptyset$. If the smallest element $\nu_0(s_\tau)$ of $\mathbf S_\tau$ satisfies $\nu_0(s_\tau)<\nu$, applying
Lemma~\ref{lem:fishboned} to $\gamma(\sigma_-)\in \mathcal D_{\nu_0(s_\tau)}$ gives \eqref{eq:rain}.
The definition \eqref{eq:defofsigma} of $\sigma_-$ shows that the case $\nu=\nu_0(s_\tau)$ is impossible. Finally, assume that $\nu_0(s_\tau)>\nu$. Then $\nu\notin\mathbf S_\tau$, and \eqref{eq:rolling1} gives 
$$\dist(\gamma(\tau),O_\nu)>\varepsilon^{-3}R_\nu/2,\qquad
\dist(\gamma(\tau),O_{\nu_0(s_\tau)})\le\varepsilon^{-3}R_{\nu_0(s_\tau)}/2.$$
Therefore, Lemma \ref{lem:jruusing} shows that 
$$\dist(p,O_\nu) > \varepsilon^{-2}R_\nu \quad \text{for all } p \in \mathcal{D}_{\nu_0(s_\tau)}.$$
Since $\gamma(\sigma_-)\in \mathcal{D}_{\nu_0(s_\tau)}$, \eqref{eq:rain} holds for $I_2$ also in this case.

\medskip
\noindent
\textbf{Case 3: $T_\mu=\emptyset$ and $H_\mu(I) =\emptyset$. } Since the arguments for the pairs $I_1,I_2$ and $I_3,I_4$ are almost identical, it suffices to prove \eqref{eq:rain} for all possible endpoints of $I_1$ and $I_2$. Note that \eqref{eq:uino1} and \eqref{eq:uino2} remain true. 

If the set $\mathbf S$ defined before \eqref{eq:rolling} is empty, then $I_1=(c,s_1)$, $I_2=(s_1,b_\mu)$, and any $\nu$ for which $(c,b_\mu) \cap \gamma^{-1}(\mathcal D_\nu) \neq \emptyset$ satisfies 
\begin{equation} \label{eq:raimo1}
\dist(\gamma(s_1),O_\nu) > \varepsilon^{-3}R_\nu/2 \geq \varepsilon^{-2}R_\nu. 
\end{equation}
Combining \eqref{eq:uino1}, \eqref{eq:uino2}, and \eqref{eq:raimo1} shows that \eqref{eq:rain} holds for $I_1$ and $I_2$. 

Next, suppose that $\mathbf S \neq \emptyset$. Then 
$I_1:=(c,s_1^-) \quad \text{and} \quad I_2:=(s_1^+,b_\mu)$, where  
$\gamma(s_1^-),\gamma(s_1^+) \in \mathcal D_{\nu(s)}$ and $\nu(s)$ is the smallest element of $\mathbf S$. Since the argument for $I_2$ and $s_1^+$ is almost identical, it suffices to prove \eqref{eq:rain} for $I_1$. We claim that 
\begin{equation} \label{eq:raimo2}
\dist(\gamma(s_1^-),O_\nu) > \varepsilon^{-2}R_\nu.  
\end{equation}
If $\nu(s) < \nu$, then \eqref{eq:raimo2} follows as \eqref{eq:uino2}, with $\mu$ replaced by $\nu(s)$. 

The definition \eqref{eq:bowling} of $s_1^-$ shows that the case $\nu(s)=\nu$ is not possible. Finally, assume that $\nu(s) > \nu$. Now we can apply Lemma \ref{lem:jruusing} with $\nu$ instead of $\mu$ and $\nu(s)$ instead of $\nu$, to conclude that \eqref{eq:raimo2} is true. Combining \eqref{eq:uino1} and \eqref{eq:raimo2} gives \eqref{eq:rain}. The proof of condition (\ref{it:kukinto2}) is complete.

\medskip
\noindent
\textbf{Proof of condition (\ref{it:kukinto3})}: Let $1 \leq \mu \leq \ell$ and $(c,d) \in \mathcal I_{\mu-1}$. The first claim \eqref{eq:kraha} follows from the construction. It suffices to prove the second claim \eqref{eq:uio} in the three cases above. Recall that the set $T_\mu$ is defined using \eqref{eq:tarve}. 

\medskip
\noindent
\textbf{Case 1: $T_\mu \neq \emptyset$. } We have $I_1=(c,b_* )$ and $I_2=(a^*,d)$, where 
$$
b_* \in \{b_\mu,\beta_\mu,\beta_\mu^-\} \quad \text{and} \quad 
a^* \in \{a_\mu,\alpha_\mu,\alpha_\mu^+\}. 
$$
Recall that 
$$
\gamma(\beta_\mu) \in \mathcal E_{k(\beta)} \quad \text{and} \quad 
\gamma(\alpha_\mu) \in \mathcal E_{k(\alpha)}, \quad 
\text{where } 1 \leq k(\beta),k(\alpha) \leq K_\mu. 
$$
By the triangle inequality, 
\begin{equation} \label{eq:samesame}
\dist(\gamma(c),\gamma(d)) \leq \dist(\gamma(c),\gamma(b_*))+
\dist(\gamma(a^*),\gamma(d))+D,  
\end{equation}
where 
\begin{equation} \label{eq:Dee} 
D:=\diam(\gamma(b_*))+\dist(\gamma(b_*),\gamma(a^*))+\diam(\gamma(a^*)). 
\end{equation} 

We claim that 
\begin{equation} \label{eq:jurha}
\diam(\gamma(b_*)) \leq \frac{3\diam(\gamma(\beta_\mu))}{2\varepsilon}. 
\end{equation} 
The claim is obvious if $b_*=\mathcal \beta_\mu$. If 
$b_* \in \{b_\mu,\beta_\mu^-\}$, then $p:=\gamma(b_*) \in \mathcal D_\nu$ for some $1 \leq \mu \leq \nu$. In particular, $r_p<\varepsilon^{10} R_\mu$. 
On the other hand, $q:=\gamma(\beta_\mu) \in \mathcal E_k$ for some $1 \leq k \leq K_\mu$, so $R_\mu \leq r_k$ and $r_q \geq \varepsilon^{10}r_k \geq \varepsilon^{10}R_\mu$. We conclude that $r_p< r_q$. Combining with Lemma \ref{lem:voivoi} proves \eqref{eq:jurha}. 

A similar estimate holds between $a^*$ and $\alpha_\mu$, and we conclude that 
\begin{equation} \label{eq:movin} 
\diam(\gamma(b_*))+\diam(\gamma(a^*)) \leq 
\frac{3(\diam(\gamma(\beta_\mu))+\diam(\gamma(\alpha_\mu))}{\varepsilon}; 
\end{equation}
note that $\beta_\mu$ and $\alpha_\mu$ may coincide. 

Next, we estimate $\dist(\gamma(b_*),\gamma(a^*))$. Since $\gamma(b_\mu),\gamma(a_\mu) \in \mathcal D_\mu$, Lemma \ref{lem:bundle} and the triangle inequality give 
\begin{equation} \label{eq:kojo}
\dist(\gamma(b_*),\gamma(a^*)) \leq \dist(\gamma(b_*),\gamma(b_\mu))
+\dist(\gamma(a_\mu),\gamma(a^*)) + 8\varepsilon^6 R_\mu. 
\end{equation} 

The rest of the estimate depends on the choices of $b_*$ and $a^*$. 
We bound the term $\dist(\gamma(b_*),\gamma(b_\mu))$, which in the case $b_*=b_\mu$ is trivial. We recall from \eqref{eq:tarve} that 
\begin{equation} \label{eq:frank}
\dist(\gamma(\beta_\mu),O_\mu) \leq \varepsilon^{-2} R_\mu. 
\end{equation} 
On the other hand, since $\gamma(b_\mu) \in \mathcal D_\mu$, Lemma \ref{lem:bundle} yields 
\begin{equation} \label{eq:frank2} 
\dist(O_\mu,\gamma(b_\mu)) \leq 4\varepsilon^6 R_\mu.
\end{equation}
Combining \eqref{eq:frank} and \eqref{eq:frank2} gives 
\begin{equation} \label{eq:imuli}
\dist(\gamma(\beta_\mu),\gamma(b_\mu)) \leq 2 \varepsilon^{-2}R_\mu. 
\end{equation} 

Moreover, if the point $\beta_\mu^- \in \mathcal D_{\nu(\beta)}$ in \eqref{eq:killi} is defined, then \eqref{eq:tilka} holds for $\nu=\nu(\beta)$, i.e., 
\begin{equation} \label{eq:haras}
\dist(O_{\nu(\beta)},\gamma(\beta_\mu)) \leq \varepsilon^{-3}R_{\nu(\beta)}/2. 
\end{equation}
Since $\nu(\beta) \geq \mu+1$, we have $R_{\nu(\beta)} \leq R_\mu$. Therefore, combining \eqref{eq:haras} and Lemma \ref{lem:bundle} gives  
\begin{equation} \label{eq:haras2} 
\dist(\gamma(\beta_\mu^-),\gamma(\beta_\mu)) \leq \varepsilon^{-3}R_{\mu}. 
\end{equation}

By combining \eqref{eq:imuli} and \eqref{eq:haras2} with the triangle inequality, we have 
\begin{equation} \label{eq:stonu1} 
\dist(\gamma(b_*),\gamma(b_\mu)) \leq 2 \varepsilon^{-3}R_\mu+\diam(\gamma(\beta_\mu)).  
\end{equation}
An identical argument gives 
\begin{equation} \label{eq:stonu2} 
\dist(\gamma(a^*),\gamma(a_\mu)) \leq 2 \varepsilon^{-3}R_\mu+\diam(\gamma(\alpha_\mu)). 
\end{equation}

We denote $M:=\max\{\diam(\gamma(\beta_\mu)),\diam(\gamma(\alpha_\mu))\}$. Combining \eqref{eq:movin}, \eqref{eq:kojo}, \eqref{eq:stonu1}, and \eqref{eq:stonu2} shows that the number $D$ in \eqref{eq:Dee} satisfies 
\begin{equation} \label{eq:sting} 
D \leq 6 \varepsilon^{-1}M+8\varepsilon^6R_\mu+4\varepsilon^{-3}R_\mu+2M 
\leq 7\varepsilon^{-1}M+5\varepsilon^{-3}R_\mu. 
\end{equation} 

To complete the estimate on $D$, we claim that 
\begin{equation} \label{eq:sting2} 
R_\mu \leq \varepsilon^{-10} M. 
\end{equation}
Indeed, recall that $q:=\gamma(\beta_\mu)$ is in $\mathcal E_{k}$ for some $1 \leq k \leq K_\mu$, so 
$$ 
r_q \geq \varepsilon^{10} r_k \geq \varepsilon^{10} R_\mu
$$
by the definition of $\mathcal E_k$ (see \eqref{eq:rdef} and \eqref{eq:edef}), and $r_q \leq \diam(q) \leq M$ by Lemma 
\ref{lem:voivoi}, so \eqref{eq:sting2} is true. Combining \eqref{eq:sting} and \eqref{eq:sting2} gives 
$$ 
D \leq 8\varepsilon^{-13} M \leq \psi(\mu,I), 
$$
which together with \eqref{eq:samesame} completes the proof \eqref{eq:uio} and condition (\ref{it:kukinto3}) in case 1.  

\medskip
\noindent
\textbf{Case 2: $T_\mu=\emptyset$ and $H_\mu(I) \neq \emptyset$. } 
Recall that the number $\tau$ defined after \eqref{eq:HmuI} satisfies $\gamma(\tau)\cap\partial E_\mu\ne\emptyset$. In particular, Lemma \ref{lem:voivoi} gives
\begin{equation} \label{eq:jodi} 
\dist(\gamma(\tau),O_\mu)\leq 
3R_\mu/\varepsilon.  
\end{equation}
We claim that
\begin{equation}\label{eq:gap-boundary-small}
\diam(\gamma(\tau))\leq \varepsilon^8R_\mu.
\end{equation}
If $\gamma(\tau)$ is a singleton, then \eqref{eq:gap-boundary-small} is obvious. 

Otherwise, we have $\gamma(\tau)\in \mathcal{C}_S$ by \eqref{eq:turki} and the assumption $|\gamma|\cap \mathcal{C}_L=\emptyset$. Since $T_\mu=\emptyset$, we have $\gamma(\tau)\notin \bigcup_{1\le k\le K_\mu}\mathcal E_k$. Moreover, since condition \eqref{it:kukinto1} is true for $I=(c,d)\in \mathcal{I}_{\mu-1}$, we also have $\gamma(\tau)\notin \bigcup_{1\le k< K_\mu}\mathcal B_k$. It follows that 
\begin{equation} \label{eq:jedi}
\gamma(\tau)\in\mathcal{C}_S^{K_\mu}\setminus\mathcal E_{K_\mu}, \quad \text{where} \quad  \mathcal{C}_S^{k}:=\mathcal{C}_S \setminus \bigcup_{\ell=1}^{k-1} (\mathcal{E}_\ell \cup \mathcal{B}_\ell). 
\end{equation}

Next, we have $r_{\gamma(\tau)} \leq R_\mu$ by \eqref{eq:jedi} and the definition \eqref{eq:rdef} of $R_\mu$. Therefore, Lemma \ref{lem:voivoi} shows that 
\begin{equation} \label{eq:jodi2}
\dist(o_{\gamma(\tau)},\gamma(\tau)) \leq \diam(\gamma(\tau)) \quad \text{and} \quad  \diam(\gamma(\tau)) \leq \frac{3R_\mu}{2 \varepsilon}.   
\end{equation}  
Combining \eqref{eq:jodi} and \eqref{eq:jodi2} gives 
\begin{equation} \label{eq:nakka} 
|o_{\gamma(\tau)}-O_\mu|\le2\diam(\gamma(\tau))+\dist(\gamma(\tau),O_\mu)
\leq 6 \varepsilon^{-1}R_\mu<\varepsilon^{-3}R_\mu. 
\end{equation} 
Combining \eqref{eq:nakka} and the definition \eqref{eq:edef} of $\mathcal E_{K_\mu}$ shows that $r_{\gamma(\tau)}<\varepsilon^{10}R_\mu$. Since $2 \varepsilon \diam(\gamma(\tau)) \leq 3 r_{\gamma(\tau)}$ by Lemma \ref{lem:voivoi}, the claim  \eqref{eq:gap-boundary-small} follows.

Next, we claim that 
\begin{equation}\label{eq:gap-endpoint-near}
\diam(\gamma(\sigma_-))+\dist(\gamma(\sigma_-),\gamma(\tau))\le\varepsilon^6R_\mu.
\end{equation}
If $\mathbf S_\tau=\emptyset$, then $\gamma(\sigma_-)=\gamma(\tau)$, so \eqref{eq:gap-endpoint-near} follows from \eqref{eq:gap-boundary-small}. Otherwise, we have $R_{\nu_0(s_\tau)}\le R_\mu$ (recall \eqref{eq:defofsigma}). Therefore, Lemma \ref{lem:voivoi}, \eqref{eq:jodi}, and \eqref{eq:gap-boundary-small} show that 
\begin{eqnarray} \label{eq:nuiju} \begin{aligned} 
|O_{\nu_0(s_\tau)}-O_\mu|
&\le\dist(O_{\nu_0(s_\tau)},\gamma(\tau))+\diam(\gamma(\tau))\\ 
&+\dist(\gamma(\tau),O_\mu) \le\varepsilon^{-3}R_\mu. 
\end{aligned}
\end{eqnarray}

Since $p_{K_{\nu_0(s_\tau)}}\in\mathcal C_S^{K_{\mu}}\setminus\mathcal E_{K_\mu}$, it follows from \eqref{eq:nuiju} and the definition \eqref{eq:edef} of $\mathcal E_{K_\mu}$ that $R_{\nu_0(s_\tau)}<\varepsilon^{10}R_\mu$.  Since 
we also have $\gamma(\sigma_-)\in\mathcal D_{\nu_0(s_\tau)}$, \eqref{eq:rolling1} and Lemma~\ref{lem:bundle} show that 
$$\diam(\gamma(\sigma_-))+\dist(\gamma(\sigma_-),\gamma(\tau))\le \bigl(8\varepsilon^6+4\varepsilon^6+\varepsilon^{-3}/2\bigr)R_{\nu_0(s_\tau)}<\varepsilon^6R_\mu.  
$$
We have proved \eqref{eq:gap-endpoint-near}. An identical argument shows that 
$$\diam(\gamma(\sigma_+))+\dist(\gamma(\sigma_+),\gamma(\tau))\le\varepsilon^6R_\mu.$$

We fix  $z\in \gamma(\tau)\cap\partial E_\mu$, and note that $|z-O_\mu|\ge R_\mu$ by the definition \eqref{eq:eepee} of $E_\mu$. 
Given $y\in \gamma(\sigma_-)$, \eqref{eq:gap-boundary-small} and
\eqref{eq:gap-endpoint-near} imply
\begin{align*}
|y-O_\mu|
&\ge |z-O_\mu| -\Bigl(
    \diam(\gamma(\sigma_-))
    +\dist(\gamma(\sigma_-),\gamma(\tau))
    +\diam(\gamma(\tau))
  \Bigr) \\
&\ge R_\mu-\varepsilon^8R_\mu-\varepsilon^6R_\mu.
\end{align*}
On the other hand, Lemma~\ref{lem:bundle} gives
$\gamma(\xi),\gamma(\eta)\subset\mathbb D(O_\mu,4\varepsilon^6R_\mu)$.
Since $1-5\varepsilon^6-\varepsilon^8\ge1/2$, we obtain
\begin{equation}\label{eq:gap-middle-distances}
\dist(\gamma(\xi),\gamma(\sigma_-))\ge R_\mu/2, \quad 
\text{similarly} \quad \dist(\gamma(\sigma_+),\gamma(\eta))\ge R_\mu/2. 
\end{equation}

Finally, by Lemma~\ref{lem:bundle}, we have $\gamma(b_\mu)\cup\gamma(a_\mu)
\subset\mathbb D(O_\mu,4\varepsilon^6R_\mu)$. Therefore, the triangle inequality gives 
\begin{equation} \label{eq:baila} 
\dist(\gamma(c),\gamma(d))\le \dist(\gamma(c),\gamma(b_\mu))
 +\dist(\gamma(a_\mu),\gamma(d))+8\varepsilon^6R_\mu. 
\end{equation}

Combining \eqref{eq:gap-middle-distances} and \eqref{eq:baila}, 
and recalling that the distances appearing there involve the endpoints of $I_j$ (see \eqref{eq:gap-four-intervals}), we conclude that
$$
\dist(\gamma(c),\gamma(d))
\le\sum_{j=1}^{4}\dist(\gamma(c_j),\gamma(d_j)). 
$$
We have proved \eqref{eq:uio}, and therefore condition \ref{it:kukinto3}, in case 2.

\medskip
\noindent
\textbf{Case 3: $T_\mu=\emptyset$ and $H_\mu(I) =\emptyset$. }
We claim that the assumptions of Proposition~\ref{prop:blockingprop} hold for the numbers $s_1,t_1$ defined in \eqref{eq:pilli}, $s_2:=b_\mu$, $t_2:=a_\mu$, and the path $\alpha:=\gamma|_{[s_1,t_1]}$.

Recall that $p_{K_\mu}\notin\gamma([s_1,t_1])$ (see the argument before \eqref{eq:pilli}). Moreover, Lemma \ref{lem:bundle} gives  
$$ 
\gamma(b_\mu),\gamma(a_\mu) \subset \mathbb D(O_\mu,4\varepsilon^6 R_\mu) \subset D(O_\mu,\varepsilon^5 R_\mu). 
$$
Since we assume $H_\mu(I) =\emptyset$, the definitions of $s_1,t_1$ in \eqref{eq:pilli} and $H_\mu(I)$ in \eqref{eq:HmuI} show that 
$$
\gamma(t) \subset E_\mu
\qquad\text{for all }s_1<t<t_1.$$
Thus, conditions~(i)--(iii) before Proposition~\ref{prop:blockingprop} are valid.

The remaining assumption is \eqref{eq:duidi}, i.e., 
\begin{equation} \label{eq:duidi'}
\diam (\gamma(t)) \le \varepsilon^5\,R_\mu \quad \text{for } t=s_1 \text{ and } t=t_1. 
\end{equation}
Recall that 
\begin{equation}\label{eq:creation}
E_\mu \subset \mathbb D(O_\mu,\varepsilon^{-2}R_\mu) 
\end{equation} 
by Lemma \ref{lem:voivoi}. Since the set $T_\mu$ is empty by assumption, its definition \eqref{eq:tarve} and \eqref{eq:creation} 
show that 
\begin{equation} \label{eq:slammi1} 
\gamma(s_1),\gamma(t_1) \notin \cup_{k=1}^{K_\mu} \mathcal E_k. 
\end{equation}
On the other hand, condition (\ref{it:kukinto1}) and the definition of $\mathcal D_\mu$ show that 
\begin{equation} \label{eq:slammi2} 
\gamma(s_1),\gamma(t_1) \notin 
 \cup_{k=1}^{K_\mu-1} \mathcal B_k. 
\end{equation} 
Combining \eqref{eq:creation}, \eqref{eq:slammi1}, and \eqref{eq:slammi2} with the definition \eqref{eq:edef} of $\mathcal E_\mu$ and Lemma \ref{lem:voivoi} shows that 
\begin{equation} \label{eq:pitka}
\diam(\gamma(s_1))+\diam(\gamma(t_1)) \leq 3\varepsilon^{9} R_\mu \leq \varepsilon^8 R_\mu. 
\end{equation}

We have proved \eqref{eq:duidi'} and thus verified the assumptions of Proposition \ref{prop:blockingprop}. We conclude that 
\begin{equation} \label{eq:dessed'} 
\dist(\gamma(s_1),\gamma(t_1)) \leq \dist(\gamma(s_1),\gamma(b_\mu))+\dist(\gamma(t_1),\gamma(a_\mu))  - \varepsilon^5 R_\mu.  
\end{equation} 

The desired estimate \eqref{eq:uio} follows if we can show that 
\begin{equation} \label{eq:nedd} 
\dist(\gamma(c),\gamma(d)) \leq \sum_{j=1}^4 \dist(\gamma(c_j),\gamma(d_j)), 
\end{equation}
where we denote the endpoints of $I_j$ by $c_j$ and $d_j$. 
We have 
\begin{equation} \label{eq:noitu1}
c_1=c, \, d_2=b_\mu, \, c_3=a_\mu, \text{ and } d_4=d. 
\end{equation}
Moreover, 
\begin{equation} \label{eq:noitu2}
d_1=s_1 \text{ or }s_1^-, \, c_2=s_1 \text{ or } s_1^+, \,  d_3=t_1 \text{ or } t_1^-, \text{ and }c_4=t_1 \text{ or } t_1^+. 
\end{equation}

By \eqref{eq:rolling} and the definition \eqref{eq:bowling} of $s_1^-$, we have 
$\gamma(s_1^-) \in \mathcal D_{\nu(s)}$ and 
\begin{equation} \label{eq:mirhi} 
\dist(O_{\nu(s)},\gamma(s_1)) \leq \varepsilon^{-3} R_{\nu(s)}/2. 
\end{equation}

Since $\nu(s) \geq \mu+1$, we have $R_{\nu(s)} \leq R_\mu$. 
Therefore, since $\dist(\gamma(s_1),O_\mu) \leq \varepsilon^{-2}R_\mu$ by Lemma \ref{lem:voivoi} and the definition of $s_1$, \eqref{eq:mirhi} and \eqref{eq:pitka} give 
\begin{eqnarray} \label{eq:juoh} \begin{aligned}
|O_{\nu(s)}-O_\mu| &\leq \dist(O_{\nu(s)},\gamma(s_1))+\diam(\gamma(s_1)) \\ &+\dist(\gamma(s_1),O_\mu) \leq \varepsilon^{-3}R_\mu. 
\end{aligned} \end{eqnarray}
Thus, the definition \eqref{eq:edef} of $\mathcal E_k$ shows that 
$R_{\nu(s)}< \varepsilon^{10}R_\mu$. Combining this with \eqref{eq:mirhi}, Lemma \ref{lem:voivoi} and Lemma \ref{lem:bundle}, we have 
\begin{equation} \label{eq:tirhi} 
\diam(\gamma(s^-_1))+ \dist(\gamma(s_1^-),\gamma(s_1)) \leq \varepsilon^6 R_\mu. 
\end{equation}

We can repeat the argument used to prove \eqref{eq:tirhi} for $s_1^+$, $t_1^-$, and $t_1^+$, to conclude that 
\begin{eqnarray} \label{eq:tirhi2} \begin{aligned}
 \diam(\gamma(s^+_1))&+\diam(\gamma(t_1^-))+\diam(\gamma(t_1^+)) 
+ \dist(\gamma(s_1^+),\gamma(s_1))\\ &+ \dist(\gamma(t_1^-),\gamma(t_1))
+\dist(\gamma(t_1^+),\gamma(t_1)) \leq 3 \varepsilon^{6}R_\mu. 
\end{aligned}
\end{eqnarray}

We now combine the estimates to prove \eqref{eq:nedd}. By 
\eqref{eq:pitka}, \eqref{eq:dessed'}, and the triangle inequality, we have 
\begin{eqnarray} \label{eq:iusto} \begin{aligned}
\dist(\gamma(c),\gamma(d)) & \leq 
\dist(\gamma(c),\gamma(s_1)) +\diam(\gamma(s_1)) \\
&+\dist(\gamma(s_1),\gamma(t_1))+\diam(\gamma(t_1))+\dist(\gamma(t_1),\gamma(d))\\ 
&\leq \dist(\gamma(c),\gamma(s_1)) + \dist(\gamma(s_1),\gamma(b_\mu))+\dist(\gamma(t_1),\gamma(a_\mu)) \\ &+\dist(\gamma(t_1),\gamma(d)) +\varepsilon^8R_\mu -\varepsilon^5 R_\mu. 
\end{aligned} 
\end{eqnarray}
Moreover, \eqref{eq:noitu1}, \eqref{eq:noitu2}, \eqref{eq:tirhi}, and the triangle inequality give 
\begin{eqnarray*} 
\dist(\gamma(c),\gamma(s_1)) &\leq& \dist(\gamma(c_1),\gamma(d_1))+\diam(\gamma(s_1^-))+
\dist(\gamma(s_1^-),\gamma(s_1))\\  
&\leq& \dist(\gamma(c_1),\gamma(d_1))+\varepsilon^6R_\mu. 
\end{eqnarray*}
We repeat the method and now use the estimates in \eqref{eq:tirhi2} to the term 
$$ 
\dist(\gamma(s_1),\gamma(b_\mu))+\dist(\gamma(t_1),\gamma(a_\mu))+\dist(\gamma(t_1),\gamma(d)), 
$$ 
so that estimates combined with \eqref{eq:iusto} give 
\begin{eqnarray*} 
\dist(\gamma(c),\gamma(d)) &\leq& \sum_{j=1}^4 \dist(\gamma(c_j),\gamma(d_j)) +4 \varepsilon^6 R_\mu + \varepsilon^8 R_\mu-\varepsilon^{5}R_\mu \\ &\leq& \sum_{j=1}^4\dist(\gamma(c_j),\gamma(d_j)). 
\end{eqnarray*}
We have proved the desired estimate \eqref{eq:nedd}. The proof of condition \eqref{it:kukinto3} is complete. 

\medskip
\noindent
\textbf{Proof of condition (\ref{it:kukinto4}}): The desired estimate \eqref{eq:uio2} for $\mu$ follows by combining \eqref{eq:kraha} and  \eqref{eq:uio} with the assumption \eqref{eq:uio2} for $\mu-1$ if $2 \leq \mu \leq \ell$ and Lemma \ref{lem:scamming} if $\mu=1$.


\section{Proof of Theorem \ref{thm:counterexample}}\label{sec:counterexample}

We prove Theorem \ref{thm:counterexample} by constructing a domain 
$\Omega \subset \hatc$ for which $\mathcal{C}_N(\Omega)$ consists of a countable collection of horizontal segments and $\mathcal{C}_P(\Omega)=\{\infty\}$. We then obtain the desired domain $\Omega'$ by applying a Möbius transformation that
maps $\Omega$ to $\Omega'$ and $\infty$ to $0$. 

\begin{theorem}\label{thm:counterexample2}
There is a domain $\Omega \subset \hatc$ such that $\mathcal C(\Omega)$ is countable, $0 \in \Omega$, and the following properties hold: 
\begin{itemize} 
\item[(i)] $\{\infty\}$ is the only element of $\mathcal{C}_P(\Omega)$ and $\mathcal C_N(\Omega)$ is locally finite in $\mathbb C$.  
\item[(ii)] If $f:\Omega \to D$ is a conformal map onto a circle domain $D$, then $\hat f(\{\infty\}) \in \mathcal C_N(D)$. 
\item[(iii)] $\Omega$ and $\Omega':=\varphi(\Omega)$, $\varphi(z)=z^{-1}$, 
satisfy the packing condition \eqref{eq:packingcondition}. 
\end{itemize}
\end{theorem}

Before proving Theorem \ref{thm:counterexample2}, we show that the domain $\Omega'$ in part (iii) satisfies the conditions in Theorem  \ref{thm:counterexample}. Since $\Omega$ is countably connected, so is $\Omega'$. Moreover, $\Omega'$ satisfies the packing condition \eqref{eq:packingcondition} by Theorem \ref{thm:counterexample2}. 

The only element of $\mathcal{C}_P(\Omega')$ is $\{\varphi(\infty)\}=\{0\}$. If $f:\Omega'\to D$ is a conformal map onto a circle domain $D$, then $g:=f \circ \varphi:\Omega \to D$ is also conformal, and part (ii) shows that $\hat{g}(\{\infty\})=\hat f(\{0\}) \in \mathcal C_N(D)$. We have proved that Theorem \ref{thm:counterexample} follows from Theorem \ref{thm:counterexample2}.  

We start the proof of Theorem \ref{thm:counterexample2} by constructing $\Omega$. In Section \ref{sec:firstconstruction}, we define ``slit squares'', i.e., the unit square from which a finite collection of horizontal slits with varying  parameters is removed. Hakobyan and Li use a similar construction \cite{HaLi23}*{Section 7}. 

In Section \ref{sec:constructionOs}, we construct the non-trivial complementary components of $\Omega$ by carrying out the slit square construction of Section \ref{sec:firstconstruction} for every $k \in \mathbb{N}$ with suitably chosen parameters depending on $k$, and translating the resulting slits to the right by $k$ units. In particular, the non-trivial complementary components of $\Omega$ lie in the strip $(0,\infty)\times (0,1)$.

\subsection{Slit squares}\label{sec:firstconstruction}
We recall the collections of dyadic squares. 

\begin{definition}
Let $\mathbb{U}=(0,1)\times(0,1)$. 
We call $Q=Q^n_{\sigma,\tau}\subset\overline{\mathbb{U}}$ a \emph{dyadic square of generation $n$} if there exist $\sigma,\tau \in \{0,1,2,\ldots,2^n-1\}$ such that 
\begin{align*}
Q^n_{\sigma,\tau}=\left[\frac{\sigma}{2^n},\frac{\sigma+1}{2^n}\right] \times \left[\frac{\tau}{2^n},\frac{\tau+1}{2^n}\right].
\end{align*}
\end{definition}

\noindent We denote by 
$$
\mathcal{D}_n=\{Q^n_{\sigma,\tau}: \, 0 \leq \sigma,\tau \leq 2^{n}-1\}, \quad n=0,1,2,\ldots, 
$$ 
the collection of all dyadic squares of generation $n$, and let $\mathcal{D}:=\cup_{n=0}^{\infty}\mathcal{D}_n$ be the collection of dyadic squares in ${\mathbb{U}}$. Notice that  
\begin{itemize} 
\item[(i)] The side length of $Q \in \mathcal{D}_n$ is $\ell(Q):=1/2^n$. 
\item[(ii)] The center of $Q^n_{\sigma,\tau} \in \mathcal{D}_n$ 
is 
$$
    (x^n_{\sigma},y^n_{\tau}):=
\left(\frac{2\sigma+1}{2^{n+1}},\frac{2\tau+1}{2^{n+1}}\right). 
$$
\end{itemize}
We fix $0<\varepsilon < 1/100$, integers $j,J$ with $\log_2(1/\varepsilon) \leq j \leq J$, and real numbers $r_n$ ($j \leq n \leq J$) satisfying 
\begin{equation} \label{naama}
1-\frac{\varepsilon}{100^n}<r_n<1, \quad j \leq n \leq J. 
\end{equation}
Let 
\begin{equation}\label{slitsdef}
  s^n_{\sigma,\tau} := \left[x^n_{\sigma}-\frac{r_n}{2^{n+1}},  x^n_{\sigma}+\frac{r_n}{2^{n+1}}\right]\times \{y^n_{\tau}\}, \quad p^n_{\sigma,\tau}:=f_n(s^n_{\sigma,\tau});   
\end{equation}
here 
\begin{equation} \label{fn}
f_n(x+iy):=\alpha_n(x-2^{-n-1})+\varepsilon+iy,  \quad \text{where} \quad \alpha_n:=\frac{1-2\varepsilon}{1-2^{-n}}<1. 
\end{equation}

The map $f_n$ scales and horizontally translates the segments $s^n_{\sigma,\tau}$. Its purpose is to modify them so that the resulting segments $p^n_{\sigma,\tau}$ satisfy Lemma \ref{centerlemma} below. Note that we have $s^n_{\sigma,\tau} \subset Q^n_{\sigma,\tau}$ and $p^n_{\sigma,\tau} \subset\mathbb{U}$, but $p^n_{\sigma,\tau}$ is not always a subset of $Q^n_{\sigma,\tau}$. We denote 
\begin{equation} \label{Sn}
\mathcal{P}_n=\{p^n_{\sigma,\tau}: \, 0 \leq \sigma,\tau \leq 2^{n}-1\}, \quad j \leq n \leq J, 
\end{equation} 
and define the slit square 
\begin{equation} \label{kinka} 
S := {\mathbb{U}}\backslash C, \quad \text{where} \quad  C:= \bigcup_ 
{p\in \mathcal C}p \quad \text{and} \quad \mathcal C:=\cup_{n=j}^J \mathcal{P}_n.  
\end{equation}

The collections of ``leftmost'', resp., ``rightmost'', slits in $\mathcal{P}_n$, namely, 
\begin{eqnarray*}
\begin{aligned}
\mathcal{P}_{n,L}
&:=\{p^n_{0,\tau}: \, 0 \leq \tau \leq 2^n-1\} \quad \text{and } \\
\mathcal{P}_{n,R}
&:=\{p^n_{2^n-1,\tau}: \, 0 \leq \tau \leq 2^n-1\}, 
\end{aligned}
\end{eqnarray*}
have the following property which follows directly from the construction. 

\begin{lemma}\label{centerlemma} 
The center of each $p^n_{0,\tau} \in \mathcal{P}_{n,L}$ is $(\varepsilon,y^n_\tau)$, and the center of each 
$p^n_{2^n-1,\tau} \in \mathcal{P}_{n,R}$ is 
$(1-\varepsilon,y^n_\tau)$. 
\end{lemma}

\subsection{Construction of \texorpdfstring{$\Omega$}{Omega}}\label{sec:constructionOs}
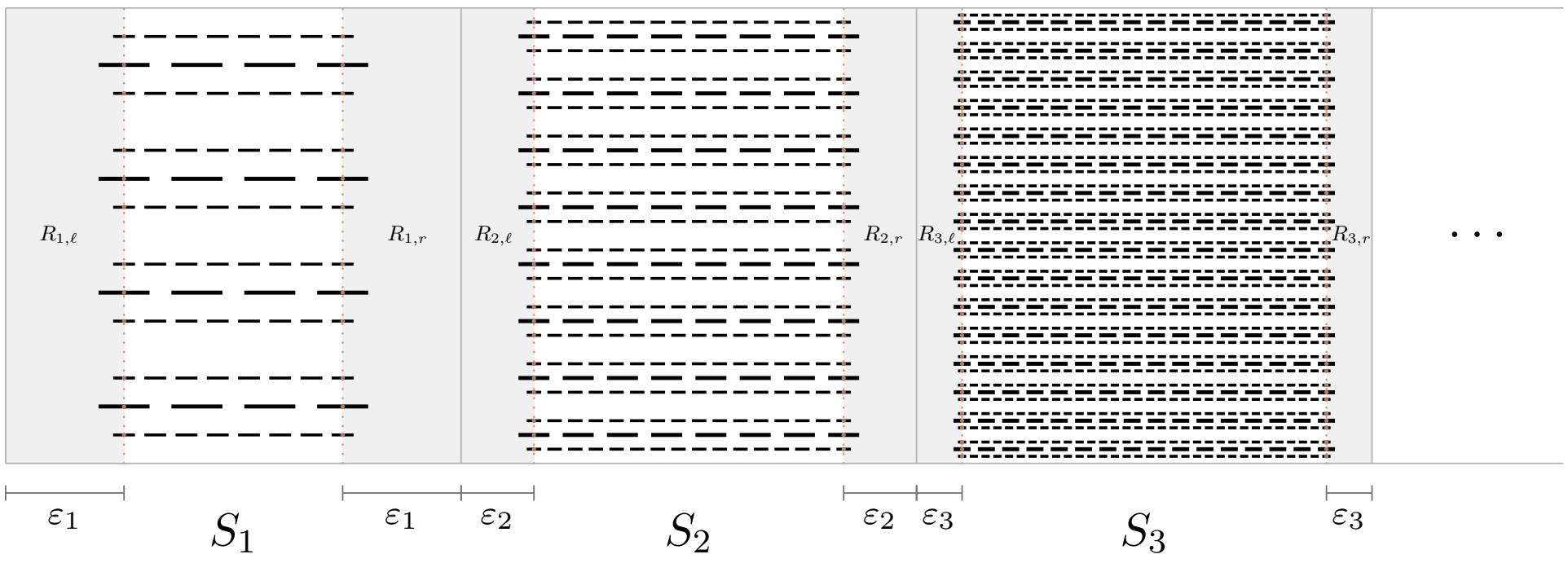
\begin{figure}[tbp]
\centering
\resizebox{\linewidth}{!}{%
\begin{tikzpicture}[
  x=3.8cm,y=3.8cm,
  line cap=butt,line join=round,
  font=\small,
  slitBoundary/.style={
    draw=black!30,line width=.35pt,solid
  },
  slitGuide/.style={
    draw=slitAccent!75,line width=.48pt,
    line cap=round,
    dash pattern=on 0pt off 2.1pt
  },
  slitDimension/.style={
    draw=black!55,line width=.35pt
  },
  slitCoarse/.style={
    draw=black,line width=.95pt,solid,line cap=butt
  },
  slitFine/.style={
    draw=black,line width=.70pt,solid,line cap=butt
  }
]

\definecolor{slitAccent}{HTML}{B47C62}

\foreach \slitK/\slitEps/\slitFirst in
  {1/.26/2,2/.16/3,3/.10/4}{
  \begin{scope}[shift={({\slitK-1},0)}]

    \fill[black!6]
      (0,0) rectangle (\slitEps,1);
    \fill[black!6]
      ({1-\slitEps},0) rectangle (1,1);

    \draw[slitGuide]
      (\slitEps,0) -- (\slitEps,1);
    \draw[slitGuide]
      ({1-\slitEps},0) -- ({1-\slitEps},1);

    \foreach \slitLevel in {0,1}{
      \pgfmathtruncatemacro{\slitN}
        {\slitFirst+\slitLevel}
      \pgfmathtruncatemacro{\slitCount}
        {2^\slitN}
      \pgfmathtruncatemacro{\slitLast}
        {\slitCount-1}

      \pgfmathsetmacro{\slitAlpha}
        {(1-2*\slitEps)/(1-1/\slitCount)}
      \pgfmathsetmacro{\slitHalf}
        {\slitAlpha*.70/(2*\slitCount)}

      \ifnum\slitLevel=0
        \tikzset{slitCurrent/.style={slitCoarse}}
      \else
        \tikzset{slitCurrent/.style={slitFine}}
      \fi

      \foreach \slitTau in {0,...,\slitLast}{
        \pgfmathsetmacro{\slitY}
          {(2*\slitTau+1)/(2*\slitCount)}

        \foreach \slitSigma in {0,...,\slitLast}{
          \pgfmathsetmacro{\slitX}
            {\slitEps+\slitAlpha*\slitSigma/\slitCount}

          \draw[slitCurrent]
            ({\slitX-\slitHalf},\slitY) --
            ({\slitX+\slitHalf},\slitY);
        }

        \fill[slitAccent]
          (\slitEps,\slitY) circle[radius=.50pt];
        \fill[slitAccent]
          ({1-\slitEps},\slitY) circle[radius=.50pt];
      }
    }

    \node[
      font=\fontsize{9}{11}\selectfont,
      inner sep=0pt,
      outer sep=0pt
    ] at ({.45*\slitEps},.5)
      {\scalebox{.56}{$R_{\slitK,\ell}$}};

    \node[
      font=\fontsize{9}{11}\selectfont,
      inner sep=0pt,
      outer sep=0pt
    ] at ({1-.45*\slitEps},.5)
      {\scalebox{.56}{$R_{\slitK,r}$}};

    \foreach \slitA/\slitB in
      {0/\slitEps,{1-\slitEps}/1}{
      \draw[slitDimension]
        (\slitA,-.065) -- (\slitB,-.065);
      \draw[slitDimension]
        (\slitA,-.048) -- (\slitA,-.082);
      \draw[slitDimension]
        (\slitB,-.048) -- (\slitB,-.082);

      \node[
        font=\scriptsize,
        anchor=north,
        inner sep=1.5pt
      ] at ({(\slitA+\slitB)/2},-.086)
        {$\varepsilon_{\slitK}$};
    }

    \node[
      font=\normalsize,
      anchor=north,
      inner sep=1.5pt
    ] at (.5,-.095)
      {$S_{\slitK}$};

  \end{scope}
}

\draw[slitBoundary] (0,0) -- (3.42,0);
\draw[slitBoundary] (0,1) -- (3.42,1);

\foreach \slitBoundaryX in {0,1,2,3}{
  \draw[slitBoundary]
    (\slitBoundaryX,0) -- (\slitBoundaryX,1);
}

\node[font=\large,inner sep=0pt]
  at (3.24,.5) {$\cdots$};

\end{tikzpicture}%
}
\caption{The elements of $\mathcal{C}_k$ for small values of $k$. The centers of the leftmost and rightmost segments in $S_k$ have $x$-coordinates $k+\varepsilon_k$ and $k+1-\varepsilon_k$, respectively. }
\label{slits}
\end{figure}

Given $k=1,2,\ldots$, we carry out the construction of the slit square $S=:S'_k$ in Section \ref{sec:firstconstruction} with the following parameters: 
\begin{eqnarray}\label{jJ}
& &    j=j(k):=(k+2)^4, \quad J=J(k):=j(k)+32(k+2)^2, \\
\label{wkdef}
& & \varepsilon=\varepsilon_k:=\frac{1}{2^{k+6}}, \quad \text{and} \\ 
\label{gaplength} 
& & r_n=r_n(k):=1-\frac{2^{n-J(k)-3}}{e^{(J(k)+2)\cdot{16}^{J(k)}}}, \quad j(k)\le n \le J(k). 
\end{eqnarray}
We then let $S_k:=S'_k+k$. 
We now construct the domain $\Omega \subset \hatc$ by specifying its complementary components. First, $\{\infty\}$ is the only element of $\mathcal{C}_P(\Omega)$. 
Secondly, recall from \eqref{kinka} that we denote the collection of the translated slits in $S_k$ by $\mathcal{C}=:\mathcal{C}_k$; see Figure \ref{slits} for an illustration. The collection $\mathcal C_N(\Omega)$ is defined as  
\begin{equation}\label{COmegas}
    \mathcal{C}_N(\Omega):=\bigcup_{k=1}^\infty \mathcal{C}_k.    
\end{equation} 

\subsection{Main estimate in \texorpdfstring{$\Omega$}{omega}} 

The most substantial part of the proof of Theorem \ref{thm:counterexample2} is part (ii). The proof is based on the following transboundary modulus estimate. We denote 
\begin{equation} \label{L}
T:=(1,\infty)\times (0,1), \quad L_0:=(1,\infty)\times \{0\}, 
\quad L_1:=(1,\infty)\times \{1\},  
\end{equation}
and let $\Gamma$ be the family of paths $\gamma:[a,b] \to \hat{\Omega}$ that join 
$L_0$ and $L_1$ in $\pi_{\Omega}(T)$. More precisely, we assume that $\gamma(a) \in L_0$, $\gamma(b) \in L_1$, and $\gamma(t) \in \pi_\Omega(T)$ for every $a<t<b$. 

\begin{theorem} \label{dualthm}
We have $\Mod \Gamma < \infty$. 
\end{theorem}

We now deduce part (ii) of Theorem \ref{thm:counterexample2} assuming Theorem \ref{dualthm}. Let $f:\Omega \to D$ be a conformal map onto a circle domain $D$, and suppose towards a contradiction that $w_0:=\hat f(\{\infty\})$ is a point component. After postcomposing with a Möbius transformation, we may assume that $w_0=0$. 

The closures of the curves $f(L_0)$ and $f(L_1)$ meet at $0$. Moreover, there is $r_0>0$ such that for every $0<r<r_0$ the length parametrized circle $\mathbb S(0,r)$ contains a subpath $\gamma_r$ that joins $f(L_0)$ and $f(L_1)$ in the domain $U \subset \hatc$ bounded by the union of $f(\partial T \cap \Omega)$ and $\{0\}$. It follows that 
\begin{equation} \label{eq:monotino}
\Lambda_U:=\{\pi_D\circ \gamma_r: \, 0<r<r_0\} \subset \hat f(\Gamma).  
\end{equation} 

We claim that 
\begin{equation} \label{eq:miksivitussa}
\modu\Lambda_U=\infty. 
\end{equation}

Given $0<t<r_0$, let $\mathcal G_t:=\{\gamma_r: 0<r<t\}$. Since $\gamma_r$ is a subpath of $\mathbb S(0,r)$, by 
\cite{Raj25b}*{Proposition 3.2} we have $\modu_C{\mathcal G_t}=\infty$ for every $0<t<r_0$; 
here $\modu_C$ is the \emph{classical modulus} of a family of paths in $\hatc$, which coincides with the transboundary modulus (Definition \ref{lem:transmoddef}) for $G=\hatc$. 

Suppose towards a contradiction that there exists an admissible function $\rho$ for $\Lambda_U$ with
$$\int_{D} \rho^2 \, dA + \sum_{p \in \mathcal{C}(D)} \rho(p)^2<\infty. $$
Given $0<t<r_0$, let 
$$\rho_t(p):=
\begin{cases}
\rho(p)\chi_{\mathbb D(0,t)}(p), & p\in D, \\
\rho(p),&
p\ne\{0\}\ \text{and}\ p\cap\mathbb D(0,t)\ne\emptyset,\\
0,&\text{otherwise}.
\end{cases}$$
Since $\gamma_r\subset\mathbb D(0,t)\setminus\{0\}$ whenever
$0<r<t$, $\rho_t$ is admissible for
$\Lambda_U(t):=\{\pi_D\circ\gamma_r:0<r<t\} \subset \Lambda_U$. By the dominated convergence theorem and the absolute continuity of the integral, we have
$$E(\rho_t):=\int_D\rho_t^2\,dA+
\sum_{p\in\mathcal C(D)}\rho_t(p)^2
\longrightarrow0
\qquad\text{as }t\to0.$$

On the other hand, by \cite{Nta26}*{Proposition 5.6} the transboundary modulus is bounded from below by the classical modulus as follows: 
$$1=\min\{1,\modu_C\mathcal G_t\}
\leq C\,
\modu
\Lambda_U(t), 
$$
where $C$ is independent of $t$. Since $\modu \Lambda_U(t) \leq E(\rho_t) \to 0$ as $t \to 0$, we arrive at a contradiction. We have proved \eqref{eq:miksivitussa}. 

It follows from \eqref{eq:monotino}, \eqref{eq:miksivitussa}, and the monotonicity (Proposition \ref{prop:propoftrans}) and conformal invariance (Lemma \ref{lem:confinvariance}) of the transboundary modulus that $\modu \Gamma = \infty$. This contradicts Theorem \ref{dualthm}. We have proved part (ii) of Theorem \ref{thm:counterexample2}, assuming Theorem \ref{dualthm}.

The proof of Theorem \ref{dualthm} is based on the following estimate on the slit squares $S_k$. Given $k \in \mathbb{N}$, we denote $\hat{S}_k:=\pi_{\Omega}(\overline{\mathbb{U}}+k)$. If $A \subset \overline T$, where $T$ is the strip in \eqref{L}, we denote 
\begin{equation}  \label{eq:heya} \begin{aligned}
\operatorname{pr}_y(A)&:=\{0 \leq y \leq 1: \, \text{there exists } x \text{ with }x+iy \in A\}, \\ \eta(A)&:=\mathcal{H}^1(\operatorname{pr}_y(A)).  
\end{aligned}
\end{equation}
We also define the following two subrectangles of $\overline{\mathbb{U}}+k$:
\begin{equation}\label{Lk'}
    R_{k,\ell}:=[k,k+\varepsilon_k]\times [0,1], \quad R_{k,r}:=[k+1-\varepsilon_k,k+1]\times [0,1];   
\end{equation} 
see Figure \ref{slits} for an illustration. 
Finally, we will use the following notation: 
$$
\text{if } \hat{A} \subset \hat{\Omega}, \text{ then } A:=\pi^{-1}_{\Omega}(\hat{A}). 
$$

\begin{proposition} \label{admissiblelemma} For every $k=1,2,\ldots$ 
there is a Borel function $\rho_k:\hat{\Omega} \to [0,\infty]$ satisfying 
$\rho_k(p)=0$ for $p \notin \hat{S}_k$ and 
\begin{equation}\label{kenergy}
    \int_{\Omega} \rho_k^2 \, dA + \sum_{p \in \mathcal{C}(\Omega)} \rho_k^2(p) \leq (k+2)^{-2}, 
\end{equation}
such that if a Borel set $\hat{A} \subset \hat{S}_k$ satisfies at least one of the conditions 
\begin{itemize} 
\item[(i)] $A \subset R_{k,\ell}\cup R_{k,r}$,  
\item[(ii)] $A$ is a continuum that intersects a vertical edge of $\mathbb{U}+k$, or 
\item[(iii)] $A$ is a continuum that connects the horizontal edges of $\mathbb{U}+k$,
\end{itemize}
then (recall notation $\eta(A)$ in \eqref{eq:heya}) 
\begin{equation} \label{aasti} 
\int_{\hat A \cap \Omega} \rho_k \, d\mathcal{H}^1 + \sum_{{p \in \hat{A} \, \cap \,\mathcal{C}(\Omega)}} \rho_k(p) \geq \,\frac{\eta(A)}{32}. 
\end{equation}
\end{proposition} 
We postpone the proof of Proposition \ref{admissiblelemma}, and first show how it yields Theorem \ref{dualthm}. We fix $\gamma \in \Gamma$ and recall that $\gamma$ connects the rays $L_0$ and $L_1$ in the projection $\pi_\Omega(T)$ of the strip $T$ (see \eqref{L}).  

\begin{lemma}\label{kirsku} 
Given $k=1,2,\ldots$, the restriction $\hat{E}_k$ of 
$|\gamma|$ to $\hat{S}_k$ is the union of pairwise disjoint Borel sets 
$\hat{A}_k(j)$, $j=0,1,2, \ldots$, each of which satisfies at least one of the conditions in Proposition \ref{admissiblelemma}. 
\end{lemma}

\begin{proof} 
We first notice that since $\gamma:[a,b] \to \hat{\Omega}$ is continuous, 
the collection of the connected components of $\hat{E}_k$ which intersect the interior of $\hat{S}_k$ in $\hat \Omega$ is at most countable. By adding empty sets if needed, we may denote their collection by 
$\hat{A}_k(j)$, $j=1,2, \ldots$, and let 
$$
\hat{A}_k(0):= \hat{E}_k \setminus \Big(\cup_{j=1}^\infty \hat{A}_k(j) \Big). 
$$ 

Notice also that if $|\gamma| \subset \hat{S}_k$, then there is only one non-empty $\hat{A}_k(j)$, which equals $|\gamma|$ and therefore connects the horizontal edges of $\mathbb{U}+k$. In particular, $A_k(j)$ satisfies condition (iii) (note that the preimage of a continuum under $\pi_\Omega$ is a continuum). Otherwise, every non-empty $A_k(j)$, $j=1,2,\ldots$, intersects a vertical segment of $\mathbb{U}+k$, and condition (ii) is true. 

The set $A_k(0)$ lies in  
the union of the vertical edges of $\mathbb{U}+k$. In particular, $A_k(0)$ satisfies condition (i). 
\end{proof}

\begin{corollary} \label{kirsku2} Let $\rho_k$ be the function in Proposition \ref{admissiblelemma}. Then 
$$
\int_{E_k\cap \Omega} \rho_k \, d\mathcal{H}^1 + \sum_{{p \in \hat{E}_k \, \cap \,\mathcal{C}(\Omega)}} \rho_k(p) \geq \,\frac{\eta(E_k)}{32}. 
$$
\end{corollary}

\begin{proof}
The estimate follows by applying \eqref{aasti} to the sets $\hat{A}=\hat{A}_k(j)$, summing over $j=0,1,2,\ldots$, and applying the subadditivity of $\mathcal{H}^1$. 
\end{proof}

We now define $\rho:=32\sum_{k=1}^\infty \rho_k$. Corollary \ref{kirsku2} 
yields 
$$
\hat \ell(\rho,\gamma)\geq 32 
\sum_{k=1}^\infty \Big(\int_{E_k\cap \Omega} \rho_k \, d\mathcal{H}^1 + \sum_{{p \in \hat{E}_k \, \cap \,\mathcal{C}(\Omega)}} \rho_k(p)\Big) \geq 
\sum_{k=1}^\infty \eta(E_k) \geq 1.  
$$
Since the estimate holds for all $\gamma \in \Gamma$, we conclude that $\rho$ is admissible for $\Gamma$. On the other hand, the supports of the functions $\rho_k$ intersect only in the edges of the squares $\overline{\mathbb{U}}+k$. Thus, Proposition \ref{admissiblelemma} shows that 
$$
\int_{\Omega} \rho^2 \, dA + \sum_{p \in \mathcal{C}(\Omega)} \rho(p)^2\le (32)^2\sum_{k=1}^\infty(k+2)^{-2}<\infty. 
$$
We have established the desired estimate $\Mod \Gamma < \infty$ assuming Proposition \ref{admissiblelemma}, which remains to be proved. 

\subsection{Proof of Proposition \ref{admissiblelemma}} \label{admissiblesection} 
We fix $k =1,2,\ldots$. To simplify notation, we suppress the dependence on $k$ and write, for example, $j=j(k)$ and $J=J(k)$.
Recall the collection $\mathcal P_n$ of slits of generation $j \leq n \leq J$ defined in \eqref{Sn}. Abusing notation, we replace the elements $p$ of $\mathcal{P}_n$ with the translated slits $p+k$, so that $\mathcal{P}_n$ 
consists of complementary components of $S_k$. 

By \eqref{slitsdef}, \eqref{fn}, and our choice of $r_n$ in \eqref{gaplength}, the distance between $p^n_{\sigma,\tau}$ and $p^n_{\sigma+1,\tau}$ is smaller than 
\begin{equation}\label{wn}
    \frac{1-r_n}{2^n}=\frac{2^{-J-3}}{e^{(J+2)\cdot{16}^{J}}}=:\omega_J 
\end{equation}
for every $\sigma=0,1,\ldots,2^{n}-2$ and $\tau =0,1,\ldots, 2^{n}-1$.  

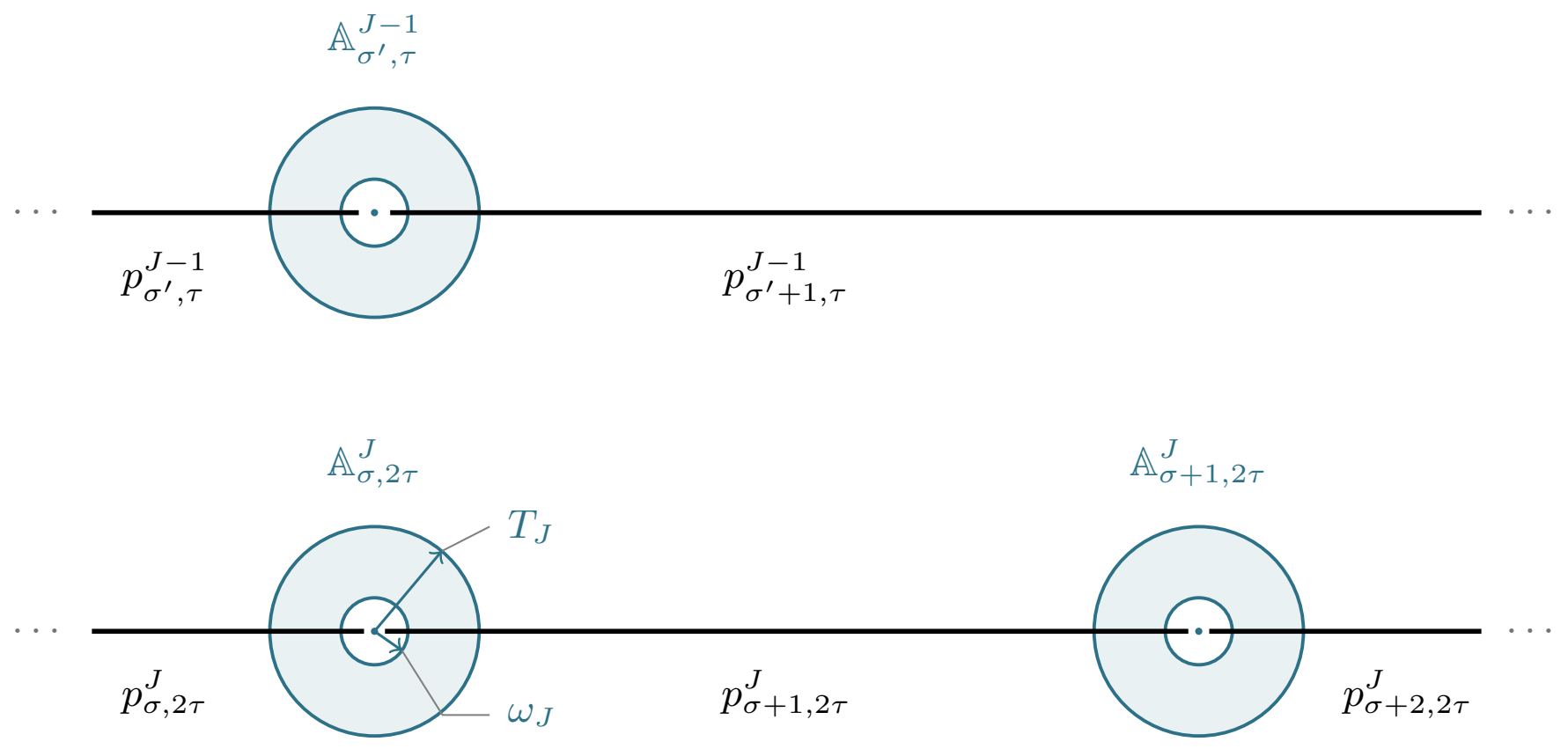
\begin{figure}[tbp]
\centering
\resizebox{0.92\linewidth}{!}{%
\begin{tikzpicture}[
  x=.82cm,y=.82cm,
  >=to,
  line cap=butt,line join=round,
  every node/.style={
    font=\fontsize{9}{11}\selectfont,
    inner sep=1.5pt
  },
  daCircle/.style={
    draw=daBlue,line width=.75pt,solid
  },
  daSlit/.style={
    draw=black,line width=1.1pt,solid
  },
  daRadius/.style={
    draw=daBlue,line width=.6pt,->
  },
  daLeader/.style={
    draw=black!50,line width=.4pt,solid
  },
  daMath/.style={
    fill=white,inner sep=1.5pt
  }
]

\definecolor{daBlue}{HTML}{2C7187}

\def\daT{1}
\def\daW{.32}
\def\daD{7.875}
\def\daH{4}
\def\daGLower{.20}
\def\daGUpper{.30}

\foreach \daX/\daY in {0/0,\daD/0,0/\daH}{
  \begin{scope}[shift={(\daX,\daY)}]
    \path[fill=daBlue!10,even odd rule]
      (0,0) circle[radius=\daT]
      (0,0) circle[radius=\daW];

    \draw[daCircle]
      (0,0) circle[radius=\daT];
    \draw[daCircle]
      (0,0) circle[radius=\daW];

    \fill[daBlue]
      (0,0) circle[radius=.8pt];
  \end{scope}
}

\draw[daSlit]
  (-2.7,0) -- ({-\daGLower/2},0);
\draw[daSlit]
  ({\daGLower/2},0) -- ({\daD-\daGLower/2},0);
\draw[daSlit]
  ({\daD+\daGLower/2},0) -- ({\daD+2.7},0);

\node[anchor=north] at (-2.0,-.30)
  {$p^J_{\sigma,2\tau}$};
\node[anchor=north] at ({\daD/2},-.30)
  {$p^J_{\sigma+1,2\tau}$};
\node[anchor=north] at ({\daD+2.0},-.30)
  {$p^J_{\sigma+2,2\tau}$};

\draw[daSlit]
  (-2.7,\daH) -- ({-\daGUpper/2},\daH);
\draw[daSlit]
  ({\daGUpper/2},\daH) -- ({\daD+2.7},\daH);

\node[anchor=north] at (-2.0,{\daH-.30})
  {$p^{J-1}_{\sigma',\tau}$};
\node[anchor=north] at ({\daD/2},{\daH-.30})
  {$p^{J-1}_{\sigma'+1,\tau}$};

\foreach \daY in {0,\daH}{
  \node[text=black!55,inner sep=0pt]
    at (-3.2,\daY) {$\cdots$};
  \node[text=black!55,inner sep=0pt]
    at ({\daD+3.2},\daY) {$\cdots$};
}

\draw[daRadius]
  (0,0) -- (50:\daT);
\draw[daLeader]
  (50:\daT) -- (1.1,1.0);
\node[text=daBlue,anchor=west] at (1.2,1.0)
  {$T_J$};

\draw[daRadius]
  (0,0) -- (-35:\daW);
\draw[daLeader]
  (-35:\daW) -- (.65,-.8) -- (1.1,-.8);
\node[text=daBlue,anchor=west] at (1.2,-.8)
  {$\omega_J$};

\node[daMath,text=daBlue,anchor=south] at (0,1.3)
  {$\mathbb A^J_{\sigma,2\tau}$};
\node[daMath,text=daBlue,anchor=south] at (\daD,1.3)
  {$\mathbb A^J_{\sigma+1,2\tau}$};
\node[daMath,text=daBlue,anchor=south] at (0,{\daH+1.3})
  {$\mathbb A^{J-1}_{\sigma',\tau}$};

\end{tikzpicture}%
}
\caption{Three annuli associated with two adjacent
slit rows of generations $J$ and $J-1$ in $S_k$. The annuli are pairwise disjoint and have inner radius $\omega_J$
and outer radius $T_J$.}
\label{fig:disjoint-slit-annuli}
\end{figure}

We denote by $G^n_{\sigma,\tau}$, $0 \leq \sigma \leq 2^{n}-2$, the relatively open interval with center $O^n_{\sigma,\tau}$ between $p^n_{\sigma,\tau}$ and $p^n_{\sigma+1,\tau}$, and denote 
\begin{equation}\label{Tndef}
    T_J:=2^{-J-3}< \frac{\ell(p^n_{\sigma,\tau})}{2} \quad \text{for every } p^n_{\sigma,\tau}\in \mathcal{P}_n.
\end{equation}
Here, the last inequality follows by construction (see \eqref{naama} and \eqref{slitsdef}).  Then \eqref{wn} yields  
\begin{equation}\label{logT/g}
    \log\frac{T_J}{\omega_J}={(J+2)\cdot{16}^{J}}. 
\end{equation}
We define the annuli with centers $O^n_{\sigma,\tau}$ (see Figure \ref{fig:disjoint-slit-annuli}), 
\begin{equation}\label{an}
    \mathbb{A}^n_{\sigma,\tau}:=\mathbb{D}(O^n_{\sigma,\tau},T_J)
\setminus
\overline{\mathbb{D}}(O^n_{\sigma,\tau},\omega_J), 
\end{equation} 
and denote 
\begin{equation}\label{annuli}
       \mathcal{A}:= \bigl\{\mathbb{A}^n_{\sigma,\tau}: j\le n\le J,\, 0 \leq \sigma \leq 2^n-2 ,\, 0\le \tau\le 2^n-1 \bigr\}. 
    \end{equation}

\begin{lemma} \label{disjointlemma}
The annuli $\mathbb{A}^n_{\sigma,\tau}$ in $\mathcal{A}$ are pairwise disjoint and satisfy 
\begin{equation}\label{AnotinLR}
    \mathbb{A}^n_{\sigma,\tau} \subset (\mathbb U+k) \setminus (R_{k,\ell}\cup R_{k,r}). 
\end{equation} 
\end{lemma} 
\begin{proof} 
We first prove disjointness. Let $\mathbb{A}^{n_1}_{\sigma_1,\tau_1},\mathbb{A}^{n_2}_{\sigma_2,\tau_2}\in\mathcal A$ be distinct annuli. It suffices to prove that
\begin{equation}\label{disjointA}
    \dist(O^{n_1}_{\sigma_1,\tau_1}, O^{n_2}_{\sigma_2,\tau_2})\ge 2T_J. 
\end{equation}

If $O^{n_1}_{\sigma_1,\tau_1}$ and $O^{n_2}_{\sigma_2,\tau_2}$ lie on the same horizontal line, then \eqref{disjointA} follows from \eqref{Tndef}. Otherwise, the distance of their heights is at least $2^{-\max\{n_1,n_2\}-1}$ by construction. Since $n_1,n_2 \leq J$, we therefore have 
$$
\dist(O^{n_1}_{\sigma_1,\tau_1}, O^{n_2}_{\sigma_2,\tau_2})\ge 2^{-\max\{n_1,n_2\}-1} \geq {2^{-J-1}}.
$$
Combining with \eqref{Tndef}, we obtain \eqref{disjointA}. 

It remains to prove \eqref{AnotinLR}. By Lemma \ref{centerlemma}, the center of each $p^n_{0,\tau} \in \mathcal{P}_{n,L}$ is $(k+\varepsilon,y^n_\tau)$. 
Combining with \eqref{Tndef}, we obtain $\mathbb{D}(O^n_{0,\tau},T_J)\cap R_{k,\ell}=\emptyset$, which implies 
$$
\mathbb{A}_{0,\tau}^{n}\subset (\mathbb U+k) \setminus R_{k,\ell}.
$$
Similarly, we have $\mathbb{A}_{2^n-2,\tau}^{n}\subset (\mathbb U+k) \setminus R_{k,r}$. Thus, \eqref{AnotinLR} holds.   
\end{proof}

We are ready to define the function $\rho=\rho_k$ in Proposition \ref{admissiblelemma}. 

\begin{definition}\label{rhokdef}
    We define $\rho=\rho_k: \hat{\Omega} \to [0,\infty]$ as follows:
    \begin{enumerate}
\item $\rho(z)=1$ for $z\in (R_{k,\ell}\cup R_{k,r})\cap \Omega$,
\item $\rho(z)=\dfrac{1}{|z-O^n_{\sigma,\tau}|\,\log\,\bigl(\tfrac{T_J}{\omega_J}\bigr)}$ for $z\in  \mathbb{A}^n_{\sigma,\tau}\cap \Omega$, where $\mathbb{A}^n_{\sigma,\tau}\in \mathcal{A}$,
\item $\rho(p)=\dfrac{2^{-n+2}}{J-j+1}$ for $p\in \mathcal{P}_n$, where $j\le n\le J$, 
\item $\rho(z)=0$ elsewhere.
\end{enumerate}
\end{definition}

\begin{lemma} \label{energylemma}
The function $\rho=\rho_k$ satisfies \eqref{kenergy}; we have 
$$
\int_{\Omega} \rho^2 \, dA + \sum_{p \in \mathcal{C}(\Omega)} \rho^2(p) \leq (k+2)^{-2}. 
 $$
\end{lemma}

\begin{proof}
By the definition of the rectangles $R_{k,\ell},R_{k,r}\subset\overline{\mathbb{U}}+k$ in \eqref{Lk'}, Lemma \ref{disjointlemma}, and the choice of $\varepsilon=\varepsilon_k$ in \eqref{wkdef}, we have  
\begin{equation*}
    \int_{(R_{k,\ell}\cup R_{k,r})\cap \Omega} \rho^2 \, dA  \le \area(R_{k,\ell}\cup R_{k,r})=2 \varepsilon= 2^{-k-5}. 
\end{equation*}

Next, by \eqref{logT/g} and \eqref{an}, every $\mathbb{A}^n_{\sigma,\tau} \in \mathcal{A}$ satisfies 
\begin{equation*}
    \int_{\mathbb{A}^n_{\sigma,\tau}\cap \Omega} \rho^2 \, dA\le \dfrac{2\pi}{\log(\tfrac{T_J}{\omega_J})}=\dfrac{2\pi}{{(J+2)\cdot{16}^{J}}}.
\end{equation*}
Combining with \eqref{annuli} and noting that $\card\mathcal{A}\le J\cdot4^{J}$, we obtain
\begin{equation*}
   \sum_{\mathbb{A}^n_{\sigma,\tau}\in \mathcal{A}} \, \int_{\mathbb{A}^n_{\sigma,\tau}\cap S_k} \rho^2 \, dA\le\dfrac{2\pi}{{4^{J}}}.
\end{equation*}

Recall that $\card\mathcal{P}_n=4^{n}$. Hence, the definition of $\rho$ on $\mathcal{P}_n$ yields 
$$\sum_{p \in \mathcal{C}(\Omega)} \rho^2(p)=\sum_{n=j}^{J}\sum_{p\in \mathcal{P}_n}\rho^2(p)\le \dfrac{16}{J-j+1}.$$

Combining the above estimates with \eqref{jJ}, we conclude that
$$ \int_{\Omega} \rho^2 \, dA + \sum_{p \in \mathcal{C}(\Omega)} \rho^2(p) \le 2^{-k-5}+ \dfrac{2\pi}{{4^{(k+2)^4}}}+\dfrac{1}{2(k+2)^{2}}\le \dfrac{1}{(k+2)^{2}}. $$ 
The proof is complete. 
\end{proof}
Recall the notation $A:=\pi^{-1}_\Omega(\hat A)$ for $\hat A \subset \hat{\Omega}$. In view of Lemma \ref{energylemma}, Proposition \ref{admissiblelemma} follows if we can prove \eqref{aasti}, i.e., 
\begin{equation} \label{aasti2} 
\int_{\hat A \cap \Omega} \rho \, d\mathcal{H}^1 + \sum_{{p \in \hat{A} \, \cap \,\mathcal{C}(\Omega)}} \rho(p) \geq \,\frac{\eta(A)}{32}, 
\end{equation}
for any set $\hat{A} \subset \hat{S}_k$ that satisfies at least one of the conditions in the proposition, namely, 
\begin{itemize} 
\item[(i)] $A \subset R_{k,\ell}\cup R_{k,r}$,  
\item[(ii)] $A$ is a continuum that intersects a vertical edge of $\mathbb{U}+k$, or 
\item[(iii)] $A$ is a continuum that connects the horizontal edges of $\mathbb{U}+k$. 
\end{itemize}

The rest of this section is devoted to the proof of \eqref{aasti2}. The next result follows directly from the definition of $\rho$. 

\begin{lemma} \label{smola}
If $\hat{A} \subset \hat{S}_k$ satisfies condition (i), then \eqref{aasti2} holds.  
\end{lemma}

It remains to prove \eqref{aasti2} under conditions (ii) and (iii).  
\begin{lemma} \label{smolo}
Suppose that $\hat{A} \subset \hat{S}_k$ satisfies $\eta(A)\leq 2^{-k-7}$ and condition (ii). Then \eqref{aasti2} holds. 
\end{lemma}

\begin{proof}
We may assume without loss of generality that $A$ intersects the left edge $I_\ell$. By Lemma \ref{centerlemma} and the definitions \eqref{wkdef} and \eqref{gaplength} of $\varepsilon=\varepsilon_k$ and $r_n=r_n(k)$, respectively, we have  
$$ 
\dist(I_\ell,p^n_{\sigma,\tau}) \geq \varepsilon-2^{-n-1}r_n \geq 2^{-k-7}\geq \eta(A) \quad 
\text{for every } p^n_{\sigma,\tau} \in \mathcal{C}_k  
$$ 
(see \eqref{COmegas}). Since $A$ is a continuum, it follows that   
$$
\int_{\hat A \cap \Omega} \rho \, d\mathcal{H}^1 \geq 
\mathcal{H}^1(A \cap R_{k,\ell} \cap \Omega) \geq  \eta(A), 
$$
which implies \eqref{aasti2}. 
\end{proof} 

Note that if $\eta(A) \leq 2^{-k-7}$ then $\hat A$ cannot satisfy condition (iii). 
We now consider continua $\hat A$ for which 
$\eta(A) \geq 2^{-k-7}$. Recall that $G_{\sigma,\tau}^n$ is the relatively open interval between $p_{\sigma,\tau}^n$ and $p_{\sigma+1,\tau}^n$. Suppose that there are $j\le n\le J$, $\sigma=0,1,\ldots,2^{n}-2$, and $\tau =0,1,\ldots, 2^{n}-1$, so that 
\begin{equation} \label{case1} 
A\cap G_{\sigma,\tau}^n\neq\emptyset \quad \text{and} \quad 
p^n_{\sigma',\tau} \notin \hat A \quad \text{for every }0 \leq \sigma' \leq 2^n-1. 
\end{equation}
\begin{lemma} \label{smallgap}
If $\hat{A}$ satisfies \eqref{case1} and $\eta(A)\geq 2^{-k-7}$, then \eqref{aasti2} holds. 
\end{lemma}
\begin{proof} 
First, note that if $G_{\sigma,\tau}^n$ and $p^{n'}_{\sigma',\tau'} \in \mathcal{P}_{n'}$ do not lie on the same horizontal line, then the distance between their heights is at least $2^{-\max\{n,n'\}-1}\geq 2^{-J-1}>T_J$. Since $T_J$ is the outer radius of 
$\mathbb{A}^n_{\sigma,\tau}$, the assumptions \eqref{case1} and $\eta(A)\geq 2^{-k-7} > 2T_J$ show that $A$ contains a subcontinuum in $\Omega$ that intersects both boundary circles $\mathbb{S}(O_{\sigma,\tau}^n,T_J)$ and $\mathbb{S}(O_{\sigma,\tau}^n,\omega_J)$ of $\mathbb{A}^n_{\sigma,\tau}$. Therefore, by Definition \ref{rhokdef} we have 
    \begin{equation*}
      \int_{\hat A \cap \Omega} \rho \, d\mathcal{H}^1\ge \int_{\hat A \cap \Omega}\dfrac{\chi_{\mathbb{A}_{\sigma,\tau}^n}(z)}{|z-O_{\sigma,\tau}^n|\,\log\,\bigl(\tfrac{T_J}{\omega_J}\bigr)}\,d\mathcal{H}^1(z)\ge 1  \geq \,\eta(A), 
    \end{equation*} 
    which implies \eqref{aasti2}. 
\end{proof}

Next, recall that each segment $p^n_{\sigma,\tau}\in \mathcal{P}_n$ lies on the horizontal line $$y=y_\tau^n=\dfrac{2\tau+1}{2^{n+1}},\qquad  0\le \tau \le 2^n-1. $$ 
We denote the set of these numbers by $Y_n$, and notice that every closed interval $I\subset [0,1]$ of length $0<|I|\le 1$ satisfies 
\begin{equation} \label{contintersect} 
\card (Y_n\cap I)\ \ge\ \lfloor 2^n |I|\rfloor>2^n|I|-1. 
\end{equation}
We set $I_{\hat A,n}:=\left\{\tau:y_{\tau}^n\in \operatorname{pr}_y(A)\right\}\subset\{0,1,2,\ldots,2^n-1\}$. 
Then \eqref{contintersect} yields 
   \begin{equation}\label{Yncappr'}
\card(I_{\hat A,n})=
\card\bigl(Y_n\cap \operatorname{pr}_y(A)\bigr) > 2^n\eta(A)-1.
\end{equation} 
We also denote $I_{\hat A,n,c}:=\{\tau\in I_{\hat A,n}: \,p^n_{\sigma,\tau} \in \hat A \,\,\text{for some}\,\, 0 \leq \sigma \leq 2^n-1\}$. 
\begin{lemma} \label{case2lemma} 
If $\hat{A}$ satisfies condition (ii) or (iii),  $\eta(A)\geq 2^{-k-7}$, and if every $j \le n \le J$ satisfies 
\begin{equation} \label{case2}
\card(I_{\hat A,n,c})\ge \tfrac{1}{2}\card(I_{\hat A,n}), 
\end{equation} 
then \eqref{aasti2} holds. 
\end{lemma} 
\begin{proof} 
By Definition \ref{rhokdef} of $\rho$, the definition \eqref{jJ} of $j$ and $J$, and \eqref{Yncappr'}, we obtain
    \begin{align}
        \sum_{p\in \hat A\cap\mathcal{C}(\Omega)}\rho(p)\ge \frac{1}{2}\sum_{n=j}^{J} \frac{\left(2^n \eta(A)-1\right)\cdot2^{-n+2}}{J-j+1}\notag
        \ge \eta(A),\notag 
    \end{align}
which implies \eqref{aasti2}.  
\end{proof}
We now consider the set $I_{\hat A,n,d}:=I_{\hat A,n} \setminus I_{\hat A,n,c}$. 

\begin{lemma} \label{squarmish}
Suppose that $\hat{A}$ satisfies condition (ii) or (iii),
$\eta(A) \geq 2^{-k-7}$, and that condition \eqref{case1}
does not hold. Let $j \leq n \leq J$, and suppose that
$\tau \in I_{\hat A,n,d}$ is such that
$\tau+1 \in I_{\hat A,n}$. Then
\begin{equation} \label{cigam}
\mathcal{H}^1
\bigl(
A \cap (R_{k,\ell} \cup R_{k,r})
\cap \Omega \cap E_{n,\tau}
\bigr)
\geq 2^{-n-3},
\end{equation}
where
$$E_{n,\tau}
:=
\bigl\{
(x,y) \in \overline{\mathbb{U}}+k:
y_\tau^n < y < y_\tau^n+2^{-n-1}
\bigr\}.$$
\end{lemma}

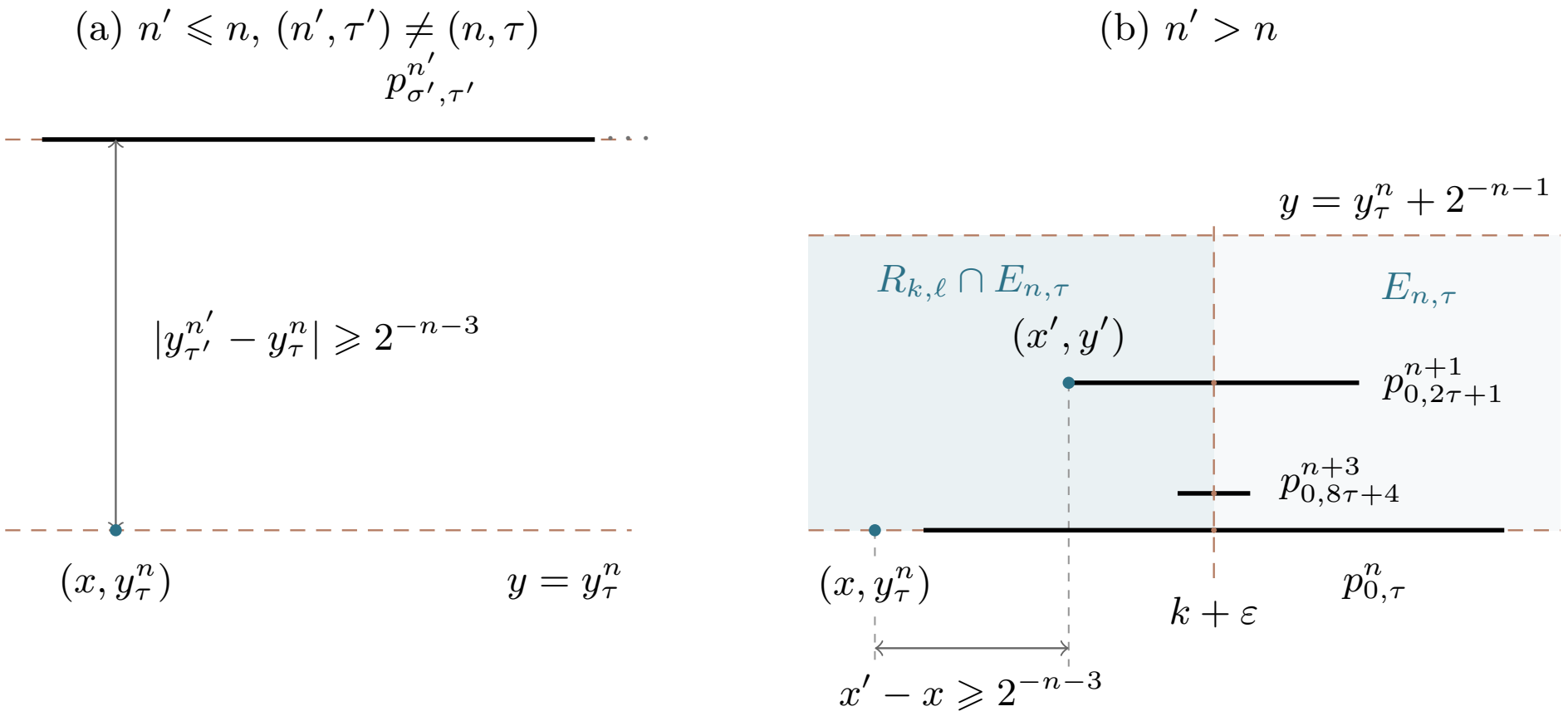
\begin{figure}[tbp]
\centering
\resizebox{0.98\linewidth}{!}{%
\begin{tikzpicture}[
  x=.60cm,y=.60cm,
  >=to,
  line cap=butt,line join=round,
  every node/.style={
    font=\fontsize{9}{11}\selectfont,
    inner sep=1.5pt
  },
  sdSlit/.style={
    draw=black,line width=1.05pt,solid
  },
  sdGuide/.style={
    draw=sdWarm,line width=.45pt,
    dash pattern=on 3.5pt off 2.5pt
  },
  sdProject/.style={
    draw=black!40,line width=.4pt,
    dash pattern=on 2pt off 2pt
  },
  sdDim/.style={
    draw=black!60,line width=.45pt,<->
  },
  sdText/.style={
    fill=white,inner sep=2pt
  },
  sdPoint/.style={
    fill=sdBlue,draw=none
  }
]

\definecolor{sdBlue}{HTML}{2C7187}
\definecolor{sdWarm}{HTML}{B47C62}

\begin{scope}
  \node at (4.6,6.8)
    {(a) $n'\leq n$, $(n',\tau')\neq(n,\tau)$};

  \draw[sdGuide] (.5,0) -- (9,0);
  \draw[sdGuide] (.5,5.3) -- (9,5.3);

  \draw[sdSlit] (1.0,5.3) -- (8.5,5.3);

  \node[anchor=south] at (6.3,5.65)
    {$p^{n'}_{\sigma',\tau'}$};

  \node[text=black!55,inner sep=0pt] at (9,5.3)
    {$\cdots$};

  \draw[sdDim] (2,0) -- (2,5.3);

  \node[align=left,anchor=west] at (2.4,2.65)
    {$|y^{n'}_{\tau'}-y^n_\tau|\geq 2^{-n-3}$};

  \fill[sdPoint] (2,0) circle[radius=1.4pt];

  \node[anchor=north] at (2,-.35)
    {$(x,y^n_\tau)$};

  \node[anchor=north east] at (9,-.35)
    {$y=y^n_\tau$};
\end{scope}

\begin{scope}[shift={(11.4,0)}]
  \def\sdC{5.5}
  \def\sdHN{3.94}
  \def\sdHF{1.97}
  \def\sdHM{.4925}
  \def\sdPX{.9}

  \pgfmathsetmacro{\sdXPrime}{\sdC-\sdHF}

  \node at (5.15,6.8)
    {(b) $n'>n$};

  \path[fill=sdBlue!4]
    (0,0) rectangle (10.2,4);
  \path[fill=sdBlue!10]
    (0,0) rectangle (\sdC,4);

  \draw[sdGuide] (0,0) -- (10.2,0);
  \draw[sdGuide] (0,4) -- (10.2,4);
  \draw[sdGuide] (\sdC,-.65) -- (\sdC,4.25);

  \node[anchor=north] at (\sdC,-.8)
    {$k+\varepsilon$};

  \node[text=sdBlue] at (2.25,3.35)
    {$R_{k,\ell}\cap E_{n,\tau}$};

  \node[text=sdBlue] at (8.3,3.25)
    {$E_{n,\tau}$};

  \node[anchor=south east] at (10.2,4.12)
    {$y=y^n_\tau+2^{-n-1}$};

  \draw[sdProject]
    (\sdPX,0) -- (\sdPX,-1.85);
  \draw[sdProject]
    (\sdXPrime,2) -- (\sdXPrime,-1.85);

  \draw[sdSlit]
    ({\sdC-\sdHN},0) -- ({\sdC+\sdHN},0);

  \node[anchor=north] at (7.7,-.35)
    {$p^n_{0,\tau}$};

  \draw[sdSlit]
    ({\sdC-\sdHF},2) -- ({\sdC+\sdHF},2);

  \node[anchor=west] at (7.7,2)
    {$p^{n+1}_{0,2\tau+1}$};

  \draw[sdSlit]
    ({\sdC-\sdHM},.5) -- ({\sdC+\sdHM},.5);

  \node[anchor=west] at (6.3,.65)
    {$p^{n+3}_{0,8\tau+4}$};

  \foreach \sdY in {0,.5,2}{
    \fill[sdWarm]
      (\sdC,\sdY) circle[radius=.65pt];
  }

  \fill[sdPoint]
    (\sdXPrime,2) circle[radius=1.4pt];

  \node[anchor=south] at (\sdXPrime,2.25)
    {$(x',y')$};

  \draw[sdDim]
    (\sdPX,-1.6) -- (\sdXPrime,-1.6);

  \node[anchor=north]
    at ({(\sdPX+\sdXPrime)/2},-1.82)
    {$x'-x\geq 2^{-n-3}$};

  \fill[sdPoint]
    (\sdPX,0) circle[radius=1.4pt];

  \node[sdText,anchor=north] at (\sdPX,-.35)
    {$(x,y^n_\tau)$};
\end{scope}

\end{tikzpicture}%
}
\caption{The distance between $(x,y_\tau^n)$ and any
$p^{n'}_{\sigma',\tau'}$ for which $n' \neq n$, or
$n=n'$ and $\tau' \neq \tau$, is at least $2^{-n-3}$. }
\label{fig:gapdistance}
\end{figure}

\begin{proof}
Recall that $A=\pi_\Omega^{-1}(\hat A)$ is a continuum.
By assumption, there is a subcontinuum $A_\tau$ of $A$
that joins the horizontal edges of $E_{n,\tau}$ in
$\overline{E}_{n,\tau}$. Let
\begin{equation}\label{defofF}
    F:=A_\tau\cap\{(x,y):y=y_\tau^n\}.
\end{equation}
Since $\tau \in I_{\hat A,n,d}$ but \eqref{case1}
does not hold, every point $(x,y_\tau^n)\in F$ satisfies
\begin{equation} \label{eq:gardenia}
x \in
\Bigl[k,k+\varepsilon-\frac{\alpha_n r_n}{2^{n+1}}\Bigr)
\cup
\Bigl(k+1-\varepsilon+\frac{\alpha_n r_n}{2^{n+1}},k+1\Bigr].
\end{equation}
Here $r_n$ and $\alpha_n$ are the constants in
\eqref{naama} and \eqref{fn}, respectively.

Without loss of generality, we assume
that $x$ lies in the segment on the left in
\eqref{eq:gardenia}. We claim that the distance between
$(x,y_\tau^n)$ and any $p^{n'}_{\sigma',\tau'}$ for which
$n'\neq n$, or $n'=n$ and $\tau'\neq\tau$, is at least
$2^{-n-3}$.

First, if $n'<n$, or $n'=n$ and $\tau'\neq\tau$, then
the $y$-coordinate $y_{\tau'}^{n'}$ of
$p^{n'}_{\sigma',\tau'}$ satisfies
$$|y_\tau^n-y_{\tau'}^{n'}|
\geq 2^{-n-1}
\geq 2^{-n-3},$$
so the claim holds in this case.

Secondly, recall that the length and the $x$-coordinate
of the center of every $p^{n'}_{\sigma',\tau'}$ are
$\alpha_{n'}r_{n'}/2^{n'}$ and
$k+\varepsilon+\alpha_{n'}\sigma'/2^{n'}$, respectively.
In particular, every
$(x',y')\in p^{n'}_{\sigma',\tau'}$ satisfies
$x'\geq
k+\varepsilon-\frac{\alpha_{n'}r_{n'}}{2^{n'+1}}.$
Therefore, if $n'>n$, we conclude using
\eqref{eq:gardenia} that
$$x'-x
\geq
\frac{\alpha_n r_n}{2^{n+1}}-\frac{\alpha_{n'}r_{n'}}{2^{n'+1}}
\geq
\frac{\alpha_n r_n}{2^{n+1}}-\frac{\alpha_{n+1}r_{n+1}}{2^{n+2}}
\geq 2^{-n-3}.$$
Thus the claim holds in all cases (see Figure \ref{fig:gapdistance}).

By \eqref{eq:gardenia}, we also have 
$$\dist\bigl(
F,
(\overline{\mathbb U}+k)
\setminus(R_{k,\ell}\cup R_{k,r})
\bigr)
\geq
\frac{\alpha_n r_n}{2^{n+1}}
\geq 2^{-n-3}.$$
We now consider the set
\begin{equation}\label{defofBtau}
    B_\tau:=\{z\in A_\tau:0<\dist(z,F)<2^{-n-3}\}\subset A_\tau.
\end{equation}
Since $F$ contains all points of $A_\tau$ on $y=y_\tau^n$, every $z=(x_0,y_0)\in B_\tau$ satisfies
$$0<y_0-y_\tau^n
\leq \dist(z,F)
<2^{-n-3}
<2^{-n-1}.$$
In particular, $B_\tau\subset E_{n,\tau}$ and $B_\tau$ does not
intersect any slit in $\overline{\mathbb{U}}+k$.
Consequently, we have
$$B_\tau\subset
A\cap(R_{k,\ell}\cup R_{k,r})
\cap\Omega\cap E_{n,\tau}.$$
Combining this with \eqref{defofF} and \eqref{defofBtau}, we conclude that
$$\mathcal H^1
\bigl(
A\cap(R_{k,\ell}\cup R_{k,r})
\cap\Omega\cap E_{n,\tau}
\bigr)\ge \mathcal H^1(B_\tau)\ge 2^{-n-3},$$
which implies \eqref{cigam}.
\end{proof}

\begin{lemma} \label{lastcase}
Suppose that $\hat{A} \subset \hat{S}_k$ satisfies
condition (ii) or (iii), that $\eta(A) \geq 2^{-k-7}$,
and that neither \eqref{case1} nor \eqref{case2} holds.
Then \eqref{aasti2} holds.
\end{lemma}

\begin{proof} 
Since \eqref{case2} fails, there is $j \leq n_0 \leq J$ such that 
\begin{equation} \label{countdetour}
\card(I_{\hat A,n_0,d})\ge \tfrac{1}{2}\card(I_{\hat A,n_0})>2^{n_0-1}\eta(A)-1; 
\end{equation}
here the last inequality follows from \eqref{Yncappr'}. On the other hand, all but at most two of the values $\tau \in I_{\hat A,n_0,d}$ 
satisfy the assumptions of Lemma \ref{squarmish}. Since the sets 
$E_{n_0,\tau}$ are pairwise disjoint, combining \eqref{cigam} with \eqref{countdetour} and the definition of $\rho$ gives 
\begin{equation} \label{kuti}
\int_{\hat{A}\cap \Omega} \rho \, d\mathcal{H}^1 \geq \frac{2^{n_0-1}\eta(A)-3}{ 2^{n_0+3}}\geq 2^{-4}\eta(A)-2^{-n_0-1} \geq \frac{\eta(A)}{32}; 
\end{equation}
here the last inequality follows from the assumption $\eta(A) \geq 2^{-k-7}$, the definition of $j$ in \eqref{jJ}, and $n_0 \geq j$. Estimate \eqref{aasti2} follows from \eqref{kuti}. 
\end{proof}

The desired estimate \eqref{aasti2} now follows by combining Lemmas \ref{smola}, \ref{smolo}, \ref{smallgap}, \ref{case2lemma}, and \ref{lastcase}. The proof of Proposition \ref{admissiblelemma} is complete.

\subsection{Proof of Theorem \ref{thm:counterexample}: the packing condition }\label{packcondition}

We now prove the last part of Theorem \ref{thm:counterexample2}, i.e., the packing condition \eqref{eq:packingcondition} for $\Omega$ and $\Omega':=\varphi(\Omega)$; here $\Omega$ is the domain in Section \ref{sec:constructionOs} and $\varphi(z):=z^{-1}$. We first claim that 
if $z_0 \in \mathbb{C}$ and $r>0$, then 
     \begin{equation}\label{packingslits}
	\card \{p \in \mathcal{C}(\Omega): \,\diam(p) \geq r, \,  p \cap \mathbb{D}(z_0,r) \neq \emptyset \} \leq 90. 
\end{equation}

Since $\diam(p)<1$ for every $p \in \mathcal{C}(\Omega)$, we lose no generality by assuming $r<1$. Moreover, since $\mathbb{D}(z_0,r)$ intersects at most three of the squares $\mathbb{U}_k:=\mathbb{U}+k$, it suffices to show that every $k \in \mathbb{N}$ satisfies 
\begin{equation}\label{cardPk}
    \card \mathcal{P}_{k}(z_0,r)\le 30, 
\end{equation}
where (recall that $\mathcal{C}_k$ is the collection of the elements of $\mathcal{C}(\Omega)$ in $\mathbb{U}_k$) 
$$\mathcal{P}_{k}(z_0,r):=\{p \in \mathcal{C}_k: \,\diam(p) \geq r, \,  p \cap \mathbb{D}(z_0,r)\neq \emptyset\}. $$

We fix $k \in \mathbb{N}$, and suppress the dependence on $k$ in the notation as in Section \ref{admissiblesection}. Note that 
$$
2^{-n-2} <\diam(p^n_{\sigma,\tau})< 2^{-n} \quad \text{for every } p^n_{\sigma,\tau} \in \mathcal{C}_k 
$$ 
by \eqref{naama}-\eqref{fn}. 
In particular, if $p^n_{\sigma,\tau} \in \mathcal{P}_k(z_0,r)$, then 
\begin{equation}\label{n0def}
n < n_0:=\Bigl\lceil \log_2\frac1r\Bigr\rceil. 
\end{equation}

Given an $n<n_0$ as well as the dyadic squares $Q^n_{\sigma,\tau}$ and the affine map $f_n$ defined in Section \ref{sec:firstconstruction}, the rectangles $f_n(Q^n_{\sigma,\tau})$ have pairwise disjoint interiors and their heights are at least $r$, so the disk $\mathbb{D}(z_0,r)$ meets the interiors of at most three rows and at most three slits of length at least $r$ in each row. It follows that (recall that $\mathcal{P}_n$ is the collection of the segments $p^n_{\sigma,\tau}$) 
\begin{equation} \label{sqaares}  
\card (\mathcal{P}_k(z_0,r) \cap \mathcal{P}_n)\le 9.  
\end{equation}

Next, we prove that there are at most two elements in $\mathcal{P}_k(z_0,r)$ 
whose diameters are significantly larger than $r$. Combining with \eqref{sqaares} then shows that $\card \mathcal{P}_k(z_0,r) \leq 3 \cdot 9+2<30$, so \eqref{cardPk} and \eqref{packingslits} hold.  

\begin{lemma} \label{largeradius} 
There are at most two elements $p^n_{\sigma,\tau}$ in $\mathcal{P}_k(z_0,r)$ such that $n \leq n_0-4$. 
\end{lemma}

\begin{proof} 
Let $n_1$ be the smallest number for which some $p^{n_1}_{\sigma_1,\tau_1}$ is in $\mathcal{P}_k(z_0,r)$, and assume that $n_1 \leq n_0-4$. Then 
$\diam(p^{n_1}_{\sigma_1,\tau_1})>2r$, which shows that there is at most one 
$\sigma_1' \neq \sigma_1$ for which also $p^{n_1}_{\sigma_1',\tau_1} \in \mathcal{P}_k(z_0,r)$. 

It remains to show that if $p^{n_2}_{\sigma_2,\tau_2} \in \mathcal{P}_k(z_0,r)$, $n_2 \leq n_0-4$, then $n_2=n_1$ and $\tau_2=\tau_1$. 
Recall that every $p^n_{\sigma,\tau} \in\mathcal{P}_n$ lies on the  horizontal line
$$ 
y_\tau^n=\frac{2\tau+1}{2^{n+1}}.  
$$
It follows that if $n_2\neq n_1$ or $\tau_2 \neq \tau_1$, then  
\begin{equation}\label{minimalr}
    |y^{n_1}_{\tau_1}-y^{n_2}_{\tau_2}|=\left|\frac{2\tau_1+1}{2^{n_1+1}}-\frac{2\tau_2+1}{2^{n_2+1}}\right|\ge \frac{1}{2^{n_0-3}}> 4r,
\end{equation}
where the last inequality follows from \eqref{n0def}. On the other hand, since both $p^{n_1}_{\sigma_1,\tau_1}$ and $p^{n_2}_{\sigma_2,\tau_2}$ intersect $\mathbb{D}(z_0,r)$, we have
$$
|y^{n_1}_{\tau_1}-y^{n_2}_{\tau_2}|\le\dist(p^{n_1}_{\sigma_1,\tau_1},p^{n_2}_{\sigma_2,\tau_2})\le\diam(\mathbb{D}(z_0,r))=2r, 
$$ 
which contradicts \eqref{minimalr}. The proof is complete. 
\end{proof}

It remains to prove the packing condition \eqref{eq:packingcondition} for $\Omega'=\varphi(\Omega)$. Let $w_0 \in \mathbb{C}$ and $s>0$. Since  $\mathcal C_N(\Omega)$ is locally finite in $\mathbb{C}$, it suffices to show that $\card \mathcal{P}^{\Omega'}_{w_0,s} \leq 10^{5}$, where 
\begin{eqnarray*}
\mathcal{P}^{\Omega'}_{w_0,s} &:=& \{\varphi(p): \, p \in \mathcal{C}', \diam(\varphi(p)) \geq s,\, \varphi(p) \cap \mathbb{D}(w_0,s) \neq \emptyset\}, \\ 
\mathcal{C}'&:=& \{p \in \mathcal{C}(\Omega): \, p \subset \mathbb{U}_k \text{ for some } k \geq 10^5\}. 
\end{eqnarray*}

Note that $\diam(\varphi(p))<10^{-10}$ for every $p \in \mathcal{C}'$, so we may assume that $s<10^{-10}$. Moreover, if $|w_0|<100s$ and if $p \in \mathcal{C}'$ is such that $\varphi(p)$ intersects $\mathbb{D}(w_0,s)$, then a computation involving $\varphi$ shows that $\diam(\varphi(p))<s$. It follows that $\mathcal{P}^{\Omega'}_{w_0,s} = \emptyset$ in this case. 

Now suppose that $|w_0| \geq 100s$ and $\varphi(p) \in \mathcal{P}^{\Omega'}_{w_0,s}$. If 
$w_1,w_2 \in \mathbb{D}(w_0,2s)$ are such that $\frac{s}{2} \leq |w_1-w_2| \leq 2s$, then applying the identity 
$$
\left|\frac{1}{w_1}  - \frac{1}{w_2} \right|=\frac{|w_1-w_2|}{|w_1||w_2|} 
$$
gives 
$$
\frac{s|z_0|^2}{3} \leq |z_1-z_2| \leq 3s|z_0|^2 \quad \text{where} \quad z_l:=w_l^{-1} \,\,\text{for}\,\, l=0,1,2.   
$$ 
It follows that if we denote $r:=s|z_0|^2$, then 
\begin{equation} \label{seima}
\diam(p) \geq \frac{r}{3}\quad \text{and} \quad 
\varphi^{-1}(\mathbb{D}(w_0,s)) \subset \mathbb{D}(z_0,3r). 
\end{equation}
We cover $\mathbb{D}(z_0,3r)$ with $1000$ disks $\mathbb{D}_m$ of radius $r/3$. 
Then \eqref{seima} yields the desired estimate 
$$
\card \mathcal{P}^{\Omega'}_{w_0,s} \leq \sum_{m=1}^{1000} \card \mathcal{P}^\Omega_m \leq 10^5, 
$$
where 
$$
\mathcal{P}^\Omega_m :=\{p \in \mathcal{C}': \, \diam(p) \geq r/3, \, p \cap \mathbb{D}_m \neq \emptyset\}, 
$$
and the last inequality follows by applying \eqref{packingslits} to each $\mathbb{D}_m$. We have proved condition \eqref{eq:packingcondition} for $\Omega'$. The proof of Theorem \ref{thm:counterexample2} is complete.


\section{Cofat and cospread domains}\label{sec:cofatspread} 
In this section, we show that the class of domains that satisfy the assumptions of Theorem \ref{thm:mainthm} contains the cofat domains, as well as the domains that satisfy the quasitripod condition \eqref{eq:quasitripod} and the packing condition \eqref{eq:packingcondition}. 

We first consider the cofat domains. Recall that $p \subset \mathbb{C}$ is \emph{$\tau$-fat}, $\tau>0$, if 
\begin{equation} \label{eq:cofat'}
\operatorname{area}(p\cap \mathbb{D}(z,r))\geq \tau\,r^2 \quad 
\text{for every } z \in p \text{ and } 0<r \leq \diam(p). 
\end{equation}
Moreover, a domain $\Omega \subset \hatc$ is \emph{cofat}, if there is $\tau>0$ such that every $p \in \mathcal{C}_N(\Omega)$ is $\tau$-fat. 

\begin{lemma} \label{lem:cofatlemma} 
Every $\tau$-cofat domain $\Omega \subset \hatc$, $\infty \in \Omega$, is $(\tau/8)$-UNC and satisfies the packing condition \eqref{eq:packingcondition} with constant $N=4\pi/\tau$. 
\end{lemma}

\begin{proof}
The packing condition \eqref{eq:packingcondition} is a straightforward consequence of \eqref{eq:cofat'} and the additivity of the Lebesgue measure. Indeed, if $p \in \mathcal{C}(\Omega)$ intersects $\D(z_0,r)$ 
and has diameter at least $r$, then 
\begin{equation}\label{eq:aksi}
\area(p \cap \mathbb{D}(z_0,2r)) \geq \tau r^2 
\end{equation}
by \eqref{eq:cofat'}. Since the elements of $\mathcal{C}(\Omega)$ are pairwise disjoint, we conclude that there can be at most $4\pi/\tau$ elements $p$ that satisfy \eqref{eq:aksi}. 

It remains to show that every $p \in \mathcal{C}_N(\Omega)$ satisfies the UNC-condition \eqref{eq:non-collineardef}. We fix $a_p,b_p \in p$ so that $|a_p-b_p|=\diam(p)$, and a third point $c_p \in p$ so that if $L$ is the line passing through $a_p$ and $b_p$, then 
$$
\dist(c_p,L)=\max_{z\in p}\dist(z,L)=:h. 
$$
We apply the fatness condition \eqref{eq:cofat'} for $z_0=a_p$ and $r=\diam(p)$ to conclude that 
$$
4h \diam(p) \geq \area(p \cap \D(a_p,\diam(p))) \geq \tau \diam(p)^2, 
$$
which implies $h \geq \tau \diam(p)/4$. In particular, if $\triangle_p$ is the triangle with vertices $a_p,b_p,c_p$, then 
\[
\operatorname{area}(\triangle_p) = \tfrac12 \diam(p) \cdot h \geq \frac{\tau \diam(p)^2}{8}. 
\] 
We conclude that the UNC-condition \eqref{eq:non-collineardef} holds with 
$\varepsilon=\tau/8$. 
\end{proof}

We now consider quasitripods, and first give the precise definitions of the concepts used in \eqref{eq:quasitripod}. 
Recall that a homeomorphism $\phi \colon  X \to Y$ between subsets of 
$\mathbb{C}$ is $\eta$-\emph{quasisymmetric} if there exists a homeomorphism $\eta:[0,\infty)\to [0,\infty)$ such that for every triple of distinct points  $z_1,z_2,z_3 \in X$, 
$$
\frac{|\phi(z_2)-\phi(z_1)|}{|\phi(z_3)-\phi(z_1)|} \leq 
\eta \left(\frac{|z_2-z_1|}{|z_3-z_1|}\right). 
$$ 
The \emph{standard tripod} $T_0 \subset \mathbb{C}$ is the union of segments $[0,e^{i \cdot 2j\pi / 3}]$, $j=0,1,2$. We call $T \subset \mathbb{C}$ an \emph{$\eta$-quasitripod} if there exists an $\eta$-quasisymmetric homeomorphism $\phi \colon T_0 \to T$. For quasisymmetric homeomorphisms defined on $T_0$, we necessarily have $\eta(1)\ge1$.

Note that the quasitripod condition given in \cite{esmayli2024conformal} is quantitatively equivalent to \eqref{eq:quasitripod} by \emph{Väisälä's theorem} (\cite{heinonen2001lectures}*{Corollary 10.22}). We now show that \eqref{eq:quasitripod} implies the UNC-condition \eqref{eq:non-collineardef}.

\begin{lemma}
Let $\Omega \subset \hatc$ be a domain that contains $\infty$. Suppose that there exists $\eta$ such that every $p \in \mathcal{C}_N(\Omega)$ contains an $\eta$-quasitripod $T$ with
\begin{equation}\label{eq:diamtdiamp}
    \diam(T) \geq \diam(p)/\eta(1). 
\end{equation}
Then $\Omega$ is $16^{-1}\eta(1)^{-5}$-UNC.   
\end{lemma} 
\begin{proof}
We fix $p\in\mathcal C_N(\Omega)$ and $T=\phi(T_0)\subset p$ such that \eqref{eq:diamtdiamp} holds. We denote the vertices of $T_0$ by $z_0,z_1,z_2$, and let $a_p:=\phi(0)\in T$. By the $\eta$-quasisymmetry of $\phi$, there are $0<\nu_j \leq 1$, $j=0,1,2$, such that  
\begin{equation}\label{defofR}
    \phi(w_j) \in \mathbb{S}(a_p,R), \quad \text{where} \quad w_j:=\nu_jz_j \quad \text{and} \quad R:=\tfrac{\diam(T)}{2\eta(1)}.
\end{equation}
We assume without loss of generality that 
\begin{equation} \label{eq:fare}
|\phi(w_1)-\phi(w_2)| \leq |\phi(w_j)-\phi(w_k)|, \quad j \neq k, 
\end{equation}
and denote $b_p:=\phi(w_1)$ and $c_p:=\phi(w_2)$. 

Note that $|w_1| \leq |w_1-w_2|$. Therefore, by quasisymmetry we have 
\begin{equation} \label{eq:loso}
\tfrac{R}{|c_p-b_p|} \leq \eta\left(\frac{|w_1|}{|w_1-w_2|} \right) \leq \eta(1). 
\end{equation} 

Let $\triangle_p$ be the isosceles triangle with vertices $a_p,b_p,c_p$. By \eqref{eq:fare}, the angle between $a_pb_p$ and $a_pc_p$ is bounded from above by $2\pi/3$. Therefore, \eqref{eq:loso}  and elementary trigonometry show that 
\begin{equation} \label{eq:kwak}
\area(\triangle_p) \geq \tfrac{R^2}{4\eta(1)}. 
\end{equation}
Since $a_p,b_p,c_p \in p$ by assumption, the desired UNC-condition \eqref{eq:non-collineardef} follows by combining \eqref{eq:diamtdiamp}, \eqref{defofR} and \eqref{eq:kwak}. 
\end{proof}

\medskip 

\medskip
\textbf{Declaration of generative AI use:} During the preparation of this manuscript, the authors used ChatGPT to assist with the preparation of the figures, routine computations,  language editing, identifying the reference in the proof of Lemma \ref{lem:sameside}, and proofreading.

	\bibliographystyle{plain}
	\bibliography{ref1.bib}

\end{document}